\documentclass[11pt]{amsart}
\usepackage{amsmath,amsthm,amsfonts,amscd,amssymb,eucal,latexsym,mathrsfs}
\usepackage[all]{xy}

\newtheorem{theorem}{Theorem}[section]
\newtheorem{corollary}[theorem]{Corollary}
\newtheorem{lemma}[theorem]{Lemma}
\newtheorem{proposition}[theorem]{Proposition}

\theoremstyle{definition}
\newtheorem{definition}[theorem]{Definition}
\newtheorem{remark}[theorem]{Remark}

\newtheorem{example}[theorem]{Example}

\newcommand{\htopol}{{\text{\rm h}}_{\text{\rm top}}}

\newcommand{\spn}{{\rm span}}

\newcommand{\re}{{\rm re}}
\newcommand{\im}{{\rm im}}
\newcommand{\Aut}{{\rm Aut}}
\newcommand{\interior}{{\rm int}}
\newcommand{\CA}{{\rm CA}}

\newcommand{\Ad}{{\rm Ad}\,}

\newcommand{\id}{{\rm id}}

\newcommand{\cb}{{\rm cb}}

\newcommand{\cB}{{\mathcal B}}
\newcommand{\cH}{{\mathcal H}}
\newcommand{\cK}{{\mathcal K}}
\newcommand{\cW}{{\mathcal W}}

\newcommand{\cC}{{\mathcal C}}
\newcommand{\cU}{{\mathcal U}}

\newcommand{\cT}{{\mathcal T}}
\newcommand{\cQ}{{\mathcal Q}}
\newcommand{\cP}{{\mathcal P}}
\newcommand{\cS}{{\mathcal S}}

\newcommand{\Cb}{{\mathbb C}}
\newcommand{\Zb}{{\mathbb Z}}

\newcommand{\Rb}{{\mathbb R}}
\newcommand{\Nb}{{\mathbb N}}

\newcommand{\even}{{\rm e}}

\newcommand{\diam}{{\rm diam}}

\newcommand{\supp}{{\rm supp}}

\newcommand{\IE}{{\rm IE}}

\newcommand{\IT}{{\rm IT}}
\newcommand{\IN}{{\rm IN}}

\newcommand{\maxtensor}{{\rm max}}
\newcommand{\Ind}{{\rm Ind}}

\newcommand{\Ab}{{\rm Ab}}
\newcommand{\Con}{{\rm Con}}
\newcommand{\NE}{{\mathcal{N}}}
\newcommand{\Inf}{{\mathcal{I}}}

\newcommand{\TSyn}{{\mathcal{TS}}}
\newcommand{\Exp}{{\rm Exp}}

\newcommand{\oA}{{\boldsymbol{A}}}
\newcommand{\oU}{{\boldsymbol{U}}}
\newcommand{\ox}{{\boldsymbol{x}}}
\newcommand{\hcb}{{\text{\rm h}_\text{\rm c}}}

\newcommand{\M}{M}
\newcommand{\LY}{{\rm LY}}

\newcommand{\WAP}{{\rm WAP}}
\newcommand{\HNS}{{\rm HNS}}

\newcommand{\RP}{{\rm RP}}

\newcommand{\CP}{{\rm CP}}
\newcommand{\WM}{{\rm WM}}

\newcommand{\op}{{\rm op}}
\newcommand{\comp}{{\rm c}}
\newcommand{\eq}{{\rm eq}}

\allowdisplaybreaks

\begin{document}

\title[Independence in topological and $C^*$-dynamics]{Independence in
topological and $C^*$-dynamics}

\author{David Kerr}
\author{Hanfeng Li}
\address{\hskip-\parindent
David Kerr, Department of Mathematics, Texas A{\&}M University,
College Station TX 77843-3368, U.S.A.}
\email{kerr@math.tamu.edu}

\address{\hskip-\parindent
Hanfeng Li, Department of Mathematics, SUNY at Buffalo,
Buffalo NY 14260-2900, U.S.A.}
\email{hfli@math.buffalo.edu}

\date{April 4, 2007}

\begin{abstract}
We develop a systematic approach to the study of independence in 
topological dynamics with an emphasis on combinatorial methods.
One of our principal aims is to
combinatorialize the local analysis of topological entropy and
related mixing properties. We also 
reframe our theory of dynamical independence in terms of tensor
products and thereby expand its scope to $C^*$-dynamics.
\end{abstract}

\maketitle

\section{Introduction}

The probabilistic notion of independence underlies several key concepts
in ergodic theory as means for expressing randomness or indeterminism.
Strong mixing, weak mixing, and ergodicity all capture the idea of asymptotic
independence, the first in a strict sense and the second two in a mean sense.
Moreover, positive entropy in $\Zb$-systems
is reflected via the Shannon-McMillan-Breiman theorem in independent behaviour
along positive density subsets of iterates (see for example Section~3 of \cite{GW}).

One can also speak of independence in topological dynamics, in which case the issue
is not the size of certain intersections as in 
the probabilistic context but rather the
simple combinatorial fact of their nonemptiness. Although the topological
analogues of strong mixing, weak mixing, and ergodicity and their relatives
form the subject of an extensive body of research (stemming in large
part from Furstenberg's seminal work on disjointness \cite{Fur,CDSRP}),
a systematic approach to independence as a unifying concept for expressing and
analyzing recurrence and mixing properties seems to be absent in the
topological dynamics literature. The present paper aims to establish such an
approach, with a particular emphasis on combinatorial arguments. At the same
time we propose a recasting of the theory in terms of tensor products.
This alternative formulation has the advantage of being applicable to
general $C^*$-dynamical systems, in line with the principle in operator
space theory that it is the tensor product viewpoint which typically
enables the quantization of concepts from Banach space theory \cite{ER,Pis}.

In fact it is from within the theory of Banach spaces that the
inspiration for our combinatorial approach to dynamical independence originates.
In the Banach space context, independence at the dual level is
associated with $\ell_1$ structure, as was strikingly
demonstrated by Rosenthal in the proof of his characterization of Banach
spaces containing $\ell_1$ isomorphically \cite{Ros}. Rosenthal's
groundbreaking $\ell_1$ theorem initiated a line of research based in Ramsey
methods that led to the work of Bourgain, Fremlin, and
Talagrand on pointwise compact sets of Baire class one functions \cite{BFT}
(see \cite{Gow,Tod} for general references).
The transfer of these ideas to the dynamical realm was initiated by K{\"o}hler,
who applied the Bourgain-Fremlin-Talagrand dichotomy for spaces of
Baire class one functions to obtain a corresponding statement for enveloping
semigroups of continuous interval maps \cite{K}.
This dynamical Bourgain-Fremlin-Talagrand dichotomy was recently extended by Glasner
and Megrelishvili to general metrizable systems \cite{GM} (see also
\cite{tame}). The dichotomy hinges on the notion
of tameness, which was introduced by K{\"o}hler (under the term regularity) and refers
to the absence of orbits of continuous functions which contain an infinite subset
equivalent to the standard basis of $\ell_1$.

The link between topological entropy and $\ell_1$ structure via coordinate
density was
discovered by Glasner and Weiss, who proved using techniques from the local
theory of Banach spaces that if a homeomorphism of a compact metric space $X$ has
zero entropy then so does the induced weak$^*$ homeomorphism of the space of
probability measures on $X$ \cite{GW}. This connection to Banach space geometry
was pursued in \cite{EID} with applications to Voiculescu-Brown entropy in
$C^*$-dynamics and then in \cite{DEBS} within a general Banach space framework.
It was shown in \cite{DEBS} that functions in the topological Pinsker
algebra can be described by the property that their orbits do not admit a positive
density subset equivalent to the standard basis of $\ell_1$. While tameness
is concerned with infinite subsets of orbits with no extra structure and
thereby calls for the application of Ramsey methods, the positive density
condition in the case of entropy reflects a tie to quantitative results in
the local theory of Banach spaces involving the Sauer-Perles-Shelah lemma and
Hilbertian geometry. What is common to both cases
is the dynamical appearance of $\ell_1$ as a manifestation of
combinatorial independence. How this link between linear geometry and
combinatorial structure plays out in the
study of local dynamical behaviour is in general not so straightforward
however, and a major goal of this paper is to understand the local situation
from a combinatorial standpoint.

Over the last ten years a substantial local theory of topological entropy
for $\Zb$-systems has unfolded around the concept of entropy pair introduced by
Blanchard in \cite{Disj}. Remarkably, every significant result to date involving entropy
pairs has been obtained using measure-dynamical techniques by way of a
variational principle. This has raised the question of whether more direct 
topological-combinatorial arguments can be found
(see for example \cite{mu}). Applying a local variational principle,
Huang and Ye have recently obtained a characterization of 
entropy pairs, and more generally of entropy tuples,
in terms of an independence property \cite{LVRA}. We will give a
combinatorial proof of this result in Section~\ref{S-entropy comb} with
a key coordinate density lemma inspired by work of Mendelson and Vershynin
\cite{MV}. Our argument has the major advantage
of portability and provides the basis for a versatile combinatorial
approach to the local analysis of entropy.
It works equally well for noninvertible surjective continuous maps and actions 
of discrete amenable groups, applies to sequence entropy where no variational 
principle exists (Section~\ref{S-null}), and is potentially of use  
in the study of Banach space geometry. The tuples of points enjoying
the independence property relevant to entropy we call IE-tuples, and in analogy
we also define IN-tuples and IT-tuples as tools for the local study of
sequence entropy and tameness, respectively. While positive entropy is
keyed to independent behaviour along positive density subsets of iterates,
what matters for positive sequence entropy and untameness are independence
along arbitrarily large finite subsets and independence along infinite subsets,
respectively (see Sections~\ref{S-null} and \ref{S-tame}). We investigate how
these local independence properties are interrelated at the global level
(Section~\ref{S-I-indep}), as well as how
their quantizations are connected to various types of asymptotic Abelianness in
$C^*$-dynamical systems (Section~\ref{S-NC}).

We begin the main body of the paper in Section~\ref{S-prelim} by laying down
the general notation and definitions that will form the
groundwork for subsequent discussions. Our setting will be that of
topological semigroups with identity acting continuously
on compact Hausdorff spaces by surjective continuous maps, or strongly
continuously on unital $C^*$-algebras
by injective unital $^*$-endomorphisms, with notable specializations to
singly generated systems in Section~\ref{S-entropy comb} (except for the last
part on actions of discrete amenable groups) and
Section~\ref{S-entropy tensor} and to actions of groups in Section~\ref{S-minimal} 
(where the main results are in fact
for the Abelian case), the second half of Section~\ref{S-NC}, and Section~\ref{S-UHF}.
In Section~\ref{S-entropy comb}
we introduce the notion of IE-tuple and 
establish several basic properties in parallel with those of entropy tuples,
including behaviour under taking products, which we prove by measure and density
arguments in a $0$--$1$ product space. We then argue that
entropy tuples are IE-tuples by applying the key coordinate density result which
appears as Lemma~\ref{L-key}.
In particular,
we recover the result of Glasner on entropy pairs in products \cite{mu}
without having invoked a variational principle. 
Using IE-tuples we give an alternative proof of the fact due to Blanchard, Glasner,
Kolyada, and Maass \cite{BGKM} that positive entropy implies Li-Yorke chaos.
Our arguments show moreover that the set of entropy pairs 
contains a dense subset consisting of Li-Yorke pairs. 
To conclude Section~\ref{S-entropy comb} we discuss how the theory 
readily extends to actions of discrete amenable groups. Note in contrast that 
the measure-dynamical approach as it has been developed for $\Zb$-systems does not 
directly extend to the general amenable case, as one needs for example
to find a substitute for the procedure of taking powers of a
single automorphism (see \cite{LVRA} and Section~19.3 of \cite{ETJ}).

We continue our discussion on entropy in Section~\ref{S-entropy tensor} by
shifting to the tensor product perspective and determining how concepts like
uniformly positive entropy translate, with a view towards formulating
the notion of a $C^*$-algebraic K-system.

In Section~\ref{S-null} we define IN-tuples as the analogue of IE-tuples for
topological sequence entropy. We give a local description of 
nullness (i.e., the vanishing of sequence entropy over all sequences)
at the Banach space level in terms of IN-pairs and show that nondiagonal
IN-pairs are the same as sequence entropy pairs as defined in \cite{NSSEP}.
Here our combinatorial approach is essential, as there is no variational
principle for sequence entropy.

Section~\ref{S-tame} concerns independence in the context of tameness.
We define IT-tuples and establish several properties in relation to untameness.
While nullness implies tameness, we illustrate with a WAP subshift example
that the converse is false. Section~\ref{S-minimal}
investigates tame extensions of minimal systems. Under the hypothesis of
an Abelian acting group, we establish in Theorem~\ref{T-decomposition} a
proximal-equicontinuous decomposition that answers a question of
Glasner from \cite{tame} and show in Theorem~\ref{T-uniquely ergodic} that
tame minimal systems are uniquely ergodic.
These theorems generalize results from \cite{NSSEP} which cover the
metrizable null case.

In Section~\ref{S-I-indep} we define I-independence as a tensor product property
that may be thought of as a $C^*$-dynamical analogue of
measure-theoretic weak mixing. Theorem~\ref{T-equiv-indep-prod} asserts that,
for systems on compact Hausdorff spaces,
I-independence is equivalent to uniform untameness (resp.\ uniform
nonnullness) of all orders and the weak mixing (resp.\ transitivity)
of all of the $n$-fold product systems, and, in the case of an
Abelian acting group,
to untameness, nonnullness, and weak mixing. We also demonstrate that, for
general $C^*$-dynamical systems, I-independence implies complete untameness.

Section~\ref{S-NC} focuses on independence in the noncommutative context.
For dynamics on simple unital nuclear $C^*$-algebras, we show that independence 
essentially amounts to Abelianness.
It is also observed that in certain situations the existence of a faithful weakly mixing
state implies independence along a thickly syndetic set.
In the opposite direction, we
prove in Section~\ref{S-UHF} that I-independence in the setting of a UHF
algebra implies weak mixing for the unique tracial state. Moreover, for Bogoliubov
actions on the even CAR algebra, I-independence is actually equivalent
to weak mixing for the unique tracial state. Continuing with the
theme of UHF algebras, we round out Section~\ref{S-UHF} by showing that, in the
type $d^\infty$ case, I-independence for $^*$-automorphisms is point-norm generic.

In the final two sections we construct an example of a tame
nonnull Toeplitz subshift (Section~\ref{S-Toeplitz}) and prove that
the action of a convergence group is null (Section~\ref{S-convergence}).

After this paper was finished Wen Huang informed us that he has shown
that every tame minimal action of an Abelian group on a compact metrizable space 
is a highly proximal extension of an equicontinuous system and
is uniquely ergodic \cite{TSSP}. Corollary~\ref{C-HP over eq} and 
Theorem~\ref{T-uniquely ergodic} in our paper strengthen these results 
from the perspective of our geometric formulation of tameness. 
We also remark that our paper answers
the last three questions in \cite{TSSP}, the first negatively and the
second two positively. 
\medskip

\noindent{\it Acknowledgements.} The initial stages of this work were carried
out while the first author was visiting the University of Tokyo under a JSPS
Postdoctoral Fellowship and the second author was at the University of Toronto. 
The first author is grateful to JSPS for support and to Yasuyuki Kawahigashi 
for hosting his stay at the University of Tokyo over the 2004--2005 academic year.
We thank Eli Glasner for some useful comments.


\section{General notation and basic notions}\label{S-prelim}

By a {\it dynamical system} we mean a pair $(X,G)$ where $X$ is
a compact Hausdorff space and $G$ is a topological semigroup with
identity with a continuous action $(s,x)\mapsto sx$ on $X$ by surjective
continuous maps.
By a {\it $C^*$-dynamical system} we mean a triple $(A,G,\alpha )$
where $A$ is a unital $C^*$-algebra, $G$ is a topological semigroup
with identity, and
$\alpha$ is an action $(s,a)\mapsto \alpha_s (a)$ of $G$ on
$A$ by injective unital $^*$-endomorphisms. The identity of $G$ will
be written $e$. We denote by $G_0$ the set $G\setminus \{ e \}$.

A dynamical system $(X,G)$ gives rise to an action $\alpha$
of the opposite semigroup $G^\op$ on $C(X)$ defined by
$\alpha_s (f)(x) = (sf)(x) = f(sx)$ for all $s\in G^\op$, $f\in C(X)$, and $x\in X$,
using the same notation for corresponding elements of $G$ and $G^\op$.
Whenever we define a property for $C^*$-dynamical systems we will also speak
of the property for dynamical systems and surjective continuous maps by
applying the definition to this associated $C^*$-dynamical system.
If a $C^*$-dynamical system is defined by a single $^*$-endomorphism
then we will talk about properties of the system as properties of
the $^*$-endomorphism whenever convenient, with a similar comment applying
to singly generated dynamical systems.

For a semigroup with identity, we write $\NE$ for the collection of nonempty
subsets of $G_0$, $\Inf$ for the collection of infinite subsets of $G_0$,
and $\TSyn$ for the collection of thickly syndetic subsets
of $G_0$. Recall that a subset $K$ of $G$ is {\it syndetic} if there is
a finite subset $F$ of $G$ such that $FK = G$ and {\it thickly syndetic}
if for every finite subset $F$ of $G$ the set $\bigcap_{s\in F} sK$ is
syndetic. When $G=\Zb$ we say that a subset $I\subseteq G$ 
has {\it positive density} if the limit
\[ \lim_{n\to\infty} \frac{|I\cap \{ -n,-n+1, \dots ,n \} |}{2n+1} \]
exists and is nonzero. When $G=\Zb_{\geq 0}$ we say that a subset 
$I\subseteq G$ has {\it positive density} if the limit
\[ \lim_{n\to\infty} \frac{|I\cap \{ 0,1, \dots ,n \} |}{n+1} \]
exists and is nonzero.

Given a $C^*$-dynamical system $(A,G,\alpha )$ and an element $a\in A$, a subset
$I\subseteq G$ is said to be an {\it $\ell_1$-isomorphism set} for $a$ if
the set $\{ \alpha_s (a) \}_{s\in I}$ is equivalent to the standard basis of
$\ell_1^I$. If this a $\lambda$-equivalence for some
$\lambda\geq 1$ then we also refer to $I$ as an 
{\it $\ell_1$-$\lambda$-isomorphism set}. These definitions also make sense more
generally for actions on Banach spaces by isometric endomorphisms.
For a dynamical system $(X,G)$ we will speak of
$\ell_1$-isomorphism sets for elements of $C(X)$ in reference to the induced action
of $G^\op$ described above.

For a dynamical system $(X,G)$ we denote by $\M (X)$ the weak$^*$ compact
convex set of Borel probability measures on $X$ and by $\M (X,G)$ the
weak$^*$ closed convex subcollection of $G$-invariant Borel probability measures.

For subsets $Y$ and $Z$ of a metric space $(X,d)$ and $\varepsilon > 0$ we write
$Y\subseteq_\varepsilon Z$ and say that $Z$ approximately includes $Y$ to
within $\varepsilon$ if for every $y\in Y$ there is a $z\in Z$ with
$d(z,y) < \varepsilon$.
For a set $X$ and an $m\in\Nb$ we write
$\Delta_m (X)=\{(x, \dots , x)\in X^m:x\in X\}$.
An element of $X^m$ which is not contained in $\Delta_m (X)$ is said to be
{\it nondiagonal}.
The linear span of a set $\Omega$ of elements in a linear space will often
be written $[\Omega ]$.



A partition of unity $\{ f_1 , \dots , f_n \}$ of a compact
Hausdorff space $X$ is said to be
{\it effective} if $\max_{x\in X} f_i (x) = 1$ for
each $i=1,\dots ,n$, which is equivalent to the existence of
$x_1 , \dots , x_n \in X$ such that $f_i (x_j ) = \delta_{ij}$.
In this case the linear map $\spn \{ f_1 , \dots ,f_n \} \to \Cb^n$
given by evaluation
at the points $x_1 , \dots , x_n$ is an isometric order isomorphism.
See Section~8 of \cite{GILQ} for more information on order
structure and finite-dimensional approximation in commutative
$C^*$-algebras.

An {\it operator space} is a closed subspace of a $C^*$-algebra, or equivalently
of $\cB (\cH )$ for some Hilbert space $\cH$. The distinguishing  
characteristic of an operator space is its collection of matrix norms, in terms 
of which one can formulate an abstract definition, the
equivalence of which with the concrete definition is a theorem of Ruan.
A linear map $\varphi : V\to W$ between operator spaces is said to be
{\it completely bounded} if 
$\sup_{n\in\Nb} \| \id_{M_n} \otimes\varphi \| < \infty$, in which case  
we refer to this supremum as the c.b.\ (completely bounded) norm,
written $\| \varphi \|_\cb$.
We say that a map $\varphi : V\to W$ between operator spaces is
a {\it $\lambda$-c.b.-isomorphism} if it is invertible and the
c.b.\ norms of $\varphi$ and $\varphi^{-1}$ satisfy
$\| \varphi \|_\cb \| \varphi^{-1} \|_\cb \leq\lambda$.
The minimal tensor product of operator spaces $V\subseteq \cB (\cH )$ and 
$W\subseteq \cB (\cK )$, written $V\otimes W$, is the closure of the 
algebraic tensor product of $V$ and $W$ under its canonical embedding into 
$\cB (\cH\otimes_2 \cK )$. When applied to closed subspaces of commutative 
$C^*$-algebras, the minimal operator space tensor product is the same as the 
Banach space injective tensor product (ignoring the matricial data).

An {\it operator system} is a closed unital self-adjoint subspace of 
a unital $C^*$-algebra.
Let $V$ be an operator system and $I,I'$ nonempty finite sets with
$I \subseteq I'$. We regard $V^{\otimes I}$ as an operator subsystem
of $V^{\otimes I'}$ under the complete order embedding given
by $v\mapsto v\otimes 1\in V^{\otimes I} \otimes
V^{\otimes I' \setminus I} = V^{\otimes I'}$. With respect to such
inclusions we define $V^{\otimes J}$ for any index set $J$ as a direct
limit over the finite subsets of $J$. For general references on operator spaces
and operator systems see \cite{ER,Pis}.

A collection $\{ (A_{i,0} , A_{i,1} )\}_{i\in I}$ of pairs of disjoint
subsets of a set $X$ is said to be {\it independent} if for every finite set
$F\subseteq I$ and $\sigma\in \{ 0,1 \}^F$ we have
$\bigcap_{i\in F} A_{i,\sigma (i)} \neq\emptyset$.

Our basic concept of dynamical independence and its quantization are given by the
following definitions.

\begin{definition}\label{D-indep}
Let $(X,G)$ be a dynamical system.
For a tuple $\oA = (A_1, A_2, \dots, A_k )$ of subsets of $X$, we say
that a set $J\subseteq G$ is an {\it independence set} for $\oA$
if for every nonempty finite subset $I\subseteq J$ and function
$\sigma : I\to \{1,2, \dots, k\}$ we have
$\bigcap_{s\in J} s^{-1} A_{\sigma(s)} \neq \emptyset$, where $s^{-1} A$
for a set $A\subseteq X$ refers to the inverse image
$\{ x\in X : sx \in A \}$.
\end{definition}

\begin{definition}\label{D-tensorindep}
Let $(A,G,\alpha )$ be a $C^*$-dynamical system, and let $V$ be a
finite-dimensional operator subsystem of $A$. Associated to every tuple
$(s_1 , \dots , s_k )$ of elements of $G$ is the linear dynamical multiplication 
map $V^{\otimes [1,k]} \to A$ determined on elementary tensors by
$a_1 \otimes\cdots\otimes a_k \mapsto \alpha_{s_1} (a_1 )
\cdots \alpha_{s_k} (a_k )$.
For $\lambda\geq 1$, we say that a tuple of elements of $G$ is a
{\it $\lambda$-contraction tuple}
for $V$ if the associated multiplication map
has c.b.\ norm at most $\lambda$, a {\it $\lambda$-expansion tuple}
for $V$ if the multiplication map
has an inverse with c.b.\ norm at most $\lambda$, and a
{\it $\lambda$-independence tuple} for $V$ if the multiplication map is a
$\lambda$-c.b.-isomorphism onto its image. A subset $J\subseteq G$ 
is said to be a {\it $\lambda$-contraction set}, {\it $\lambda$-expansion set}, or
{\it $\lambda$-independence set} if every tuple of distinct elements in $J$ 
is of the corresponding type.
\end{definition}

\begin{remark}\label{R-comm}
If $A$ is a commutative $C^*$-algebra, then for each of the linear
maps in Definition~\ref{D-tensorindep} the norm and
c.b.\ norm coincide. Moreover the linear
map $V^{\otimes [1,k]} \to A$
is contractive and thus $\varphi$ is a $\lambda$-isomorphism onto its image
if and only if it has a bounded inverse of norm at most $\lambda$.
Thus in this case a $\lambda$-expansion tuple (resp.\ set) is the same as a
$\lambda$-independence tuple (resp.\ set).
\end{remark}

Notice that if a tuple $(s_1 , \dots , s_k )$ of elements of $G$ is
a $\lambda$-expansion tuple for a finite-dimensional operator subsystem
$V\subseteq A$ with $\lambda$ close to one, then the associated multiplication
map $\varphi$ gives rise to an operator space matrix norm
on $V^{\otimes [1,k]}$ which is close to being a cross norm in the sense that for
all $a_1 , \dots , a_k \in V$ we have
\begin{align*}
\lambda^{-1} \| a_1 \| \cdots \| a_k \| &=
\lambda^{-1} \| a_1 \otimes\cdots\otimes a_k \| \leq
\| \varphi (a_1 \otimes\cdots\otimes a_k ) \| \\
&= \| \alpha_{s_1} (a_1 ) \cdots \alpha_{s_k} (a_k ) \| \leq \| a_1 \| \cdots
\| a_k \| .
\end{align*}

\begin{definition}\label{D-contr}
Let $(A,G,\alpha )$ be a $C^*$-dynamical system, and let $V$ be a
finite-dimensional operator subsystem of $A$. For $\varepsilon\geq 0$,
we define $\Con (\alpha , V, \varepsilon )$ to be the set of all
$s\in G_0$ such that $(e,s)$ is a $(1+\varepsilon )$-contraction tuple,
$\Exp (\alpha , V, \varepsilon )$ to be the set of all
$s\in G_0$ such that $(e,s)$ is a $(1+\varepsilon )$-expansion tuple,
and $\Ind (\alpha , V, \varepsilon )$ to be the set of all $s\in G_0$ such
that $(e,s)$ is a $(1+\varepsilon )$-independence tuple.
For a collection $\cC$ of subsets of $G_0$ which is closed under taking
supersets, we say that
the system $(A,G,\alpha )$ or the action $\alpha$ is
{\it $\cC$-contractive} if for every finite-dimensional operator subsystem
$V\subseteq A$ and $\varepsilon > 0$ the set
$\Con (\alpha , V, \varepsilon )$ is a member of $\cC$,
{\it $\cC$-expansive} if the same criterion holds with respect to the set
$\Exp (\alpha , V, \varepsilon )$, and {\it $\cC$-independent} if the
criterion holds with respect to the set $\Ind (\alpha , V, \varepsilon )$.
\end{definition}

We will mainly be interested in applying Definition~\ref{D-contr} to the
collections $\NE$, $\Inf$, and $\TSyn$ as defined above.
Note that, by Remark~\ref{R-comm}, for dynamical systems $\cC$-contractivity is
automatic and $\cC$-expansivity and $\cC$-independence amount
to the same thing.

To verify $\cC$-contractivity, $\cC$-expansivity, or $\cC$-independence it
suffices to check the condition in question over a collection of operator
subsystems with dense union. This fact is recorded in Proposition~\ref{P-dense}
and rests on the following perturbation lemma, which is a slight variation on
Lemma~2.13.2 of \cite{Pis} with essentially the same proof.

\begin{lemma}\label{L-perturb}
Let $V$ be a finite-dimensional operator space with Auerbach system\linebreak
$\{ (v_i , f_i ) \}_{i=1}^n$ and let $\varepsilon > 0$ be such that
$\varepsilon (1+\varepsilon ) < 1$.
Let $W$ be an operator space, $\rho : V\to W$ a linear map which is an
isomorphism onto its image with
$\max (\| \rho \|_\cb , \| \rho^{-1} \|_\cb ) < 1+\varepsilon$,
and $w_1 , \dots ,w_n$ elements of $W$ such that
$\| \rho (v_i ) - w_i \| < \dim (V)^{-1} \varepsilon$ for each
$i=1,\dots , n$. Then the linear map
$\varphi : V\to W$ determined by $\varphi (v_i ) = w_i$ for $i=1,\dots ,n$
is an isomorphism onto its image with
$\| \varphi \|_\cb < 1 + 2\varepsilon$ and
$\| \varphi^{-1} \|_\cb <
(1+\varepsilon )(1 - \varepsilon (1+\varepsilon ))^{-1}$.
\end{lemma}

\begin{proof}
Define the linear map $\delta : V\to W$ by $\delta (v) =
\sum_{i=1}^n f_i (v) (w_i - \rho (v_i ))$ for all $v\in V$. Since the
norm and c.b.\ norm of a rank-one linear map coincide, we have
$\| \delta \|_\cb \leq \sum_{i=1}^n \| f_i \| \| w_i - \rho (v_i ) \|
< \varepsilon$, and thus since $\varphi = \rho + \delta$ we see that
$\| \varphi \|_\cb < 1 + 2\varepsilon$. The bound $\| \varphi^{-1} \|_\cb
< (1+\varepsilon )(1 - \varepsilon (1+\varepsilon ))^{-1}$ follows
by applying Lemma~2.13.1 of \cite{Pis}.
\end{proof}

\begin{proposition}\label{P-dense}
Let $(A,G,\alpha )$ be a $C^*$-dynamical system. Let
$\mathfrak{S}$ be a collection of finite-dimensional operator subsystems of
$A$ with the property that for every finite set $\Omega\subseteq A$ and
$\varepsilon > 0$ there is a $V\in\mathfrak{S}$ such that
$\Omega\subseteq_\varepsilon V$. Let $\cC$ be a collection of subsets of $G_0$
which is closed under taking supersets.
Then $\alpha$ is $\cC$-contractive if and only if for every $V\in\mathfrak{S}$
and $\varepsilon > 0$ the set $\Con (\alpha , V, \varepsilon )$ is a member
of $\cC$, $\cC$-expansive if and only if for every $V\in\mathfrak{S}$
and $\varepsilon > 0$ the set $\Exp (\alpha , V, \varepsilon )$ is a member
of $\cC$, and $\cC$-independent if and only if for every $V\in\mathfrak{S}$
and $\varepsilon > 0$ the set $\Ind (\alpha , V, \varepsilon )$ is a member
of $\cC$.
\end{proposition}

\begin{proof}
We will show the third equivalence, the first two involving similar
perturbation arguments.
For the nontrivial direction, suppose that for every $V\in\mathfrak{S}$
and $\varepsilon > 0$ the set $\Ind (\alpha , V, \varepsilon )$ is a member
of $\cC$. Let $V$ be a finite-dimensional operator subsystem of
$A$ and let $\varepsilon > 0$. Let $\delta$ be a positive real number to
be further specified below. By assumption we can find a
$W\in\mathfrak{S}$ such that the set $\Ind (\alpha , W, \delta )$ is
a member of $\cC$ and the unit ball of $V$ is approximately included to
within $\delta$ in $W$. By Lemma~\ref{L-perturb},
if $\delta$ is sufficiently small then, taking an Auerbach basis
$\cS = \{ v_i \}_{i=1}^r$ for $V$ and choosing $w_1, \dots , w_r \in W$ with
$\| w_i - v_i \| < \delta$ for each $i=1, \dots r$, the linear map
$\rho : V\to W$ determined on $\cS$ by $\rho (v_i ) = w_i$
is an isomorphism onto its image with
$\| \rho \|_\cb \| \rho^{-1} \|_\cb < \sqrt[4]{1 + \varepsilon}$.
Now let $s\in \Ind (\alpha , W, \delta )$.
Let $\varphi : W\otimes W \to A$  be the dynamical multiplication map
determined on elementary tensors by
$w_1 \otimes w_2 \mapsto w_1 \alpha_s (w_2 )$.
Consider the linear map $\theta : [\rho (V)\alpha_s (\rho (V))] \to
[V\alpha_s (V)]$ determined by
$\rho (v_1 ) \alpha_s (\rho (v_2 )) \mapsto v_1 \alpha_s (v_2 )$ for
$v_1 , v_2 \in V$, which is well defined by our choice of $\rho$ and $s$.
Another application of Lemma~\ref{L-perturb} shows that if $\delta$
is small enough then the composition $\theta\circ\varphi$ is an isomorphism
onto its image satisfying
$\| \theta\circ\varphi \|_\cb \| (\theta\circ\varphi )^{-1} \|_\cb
< \sqrt{1 + \varepsilon}$. Notice now that the dynamical multiplication map
$\psi : V\otimes V \to A$ determined on elementary tensors
by $v_1 \otimes v_2 \mapsto v_1 \alpha_s (v_2 )$ factors as
$\theta\circ\varphi\circ (\rho\otimes\rho )$, and hence
\[ \| \psi \|_\cb \| \psi^{-1} \|_\cb \leq
\| \theta\circ\varphi \|_\cb \| (\theta\circ\varphi )^{-1} \|_\cb
\| \rho \|_\cb^2 \| \rho^{-1} \|_\cb^2 \leq 1+\varepsilon . \]
Thus $s\in \Ind (\alpha , V, \varepsilon )$, and we conclude that $\alpha$ is
$\cC$-independent.
\end{proof}


\section{Topological entropy and combinatorial independence}\label{S-entropy comb}

The local theory of topological entropy based on entropy pairs is
developed in the literature for $\Zb$-systems, but here will we
consider general continuous surjective maps. In fact one of the novel features
of our combinatorial approach is that it applies not only to singly
generated systems but also to actions of any discrete amenable group, as we will 
indicate in the last part of the section. Thus, with the exception of the last part 
of the section, $G$ will be one of the additive semigroups $\Zb$ and $\Zb_{\geq 0}$
and we will denote the generating surjective endomorphism of $X$ by
$T$. 

Recall that the topological entropy $\htopol (T,\cU )$ of an open
cover $\cU$ of $X$ with respect to $T$ is defined as
$\lim_{n\to\infty} \frac{1}{n}\ln N(\cU \vee T^{-1} \cU
\vee\cdots\vee T^{-n+1} \cU )$, where $N(\cdot )$ denotes the
minimal cardinality of a subcover. The topological entropy $\htopol
(T)$ of $T$ is the supremum of $\htopol (T,\cU )$ over all open
covers $\cU$ of $X$.
A pair $(x_1 ,x_2 )\in X^2\setminus \Delta_2(X)$ is said to be an
{\it entropy pair} if whenever $U_1$ and $U_2$ are closed disjoint
subsets of $X$ with $x_1 \in\interior (U_1 )$ and $x_2 \in\interior
(U_2 )$, the open cover $\{ U_1^\comp , U_2^\comp \}$ has positive
topological entropy. More generally, following \cite{TEE} we call a
tuple $\ox = (x_1 , \dots, x_k)\in X^k\setminus \Delta_k(X)$ an {\it
entropy tuple} if whenever $U_1, \dots , U_l$ are closed pairwise
disjoint neighbourhoods of the distinct points in the list 
$x_1 , \dots , x_k$, the open cover 
$\{ U_1^\comp , \dots , U_l^\comp \}$ has positive topological entropy.

\begin{definition}\label{D-IE-pair}
We call a tuple $\ox = (x_1 , \dots, x_k )\in X^k$ an {\it IE-tuple}
(or an {\it IE-pair} in the case $k=2$) if for every product
neighbourhood $U_1 \times\cdots\times U_k$ of $\ox$ the tuple $(U_1,
\dots, U_k)$ has an independence set of positive density. We denote
the set of $\IE$-tuples of length $k$ by $\IE_k (X,G)$.
\end{definition}

The argument in the second paragraph of the proof of Theorem~3.2 in
\cite{GW} shows the following lemma, which will be repeatedly useful
for converting finitary density statements to infinitary ones.

\begin{lemma}\label{L-density}
A tuple $\oA = (A_1 , \dots , A_k )$ of subsets of $X$ has an
independence set of positive density if and only if there exists a
$d>0$ such that for any $M>0$ we can find an interval $I$ in $G$
with $|I|\geq M$ and an independence set $J$ for $\oA$ contained in
$I$ for which $|J|\geq d|I|$.
\end{lemma}

On page 684 of \cite{GW} a pair $(x_1 , x_2 )\in X^2\setminus
\Delta_2(X)$ is defined to be an {\it E-pair} if for any
neighbourhoods $U_1$ and $U_2$ of $x_1$ and $x_2$, respectively,
there exists a $\delta > 0$ and a $k_0$ such that for every $k\geq
k_0$ there exists a sequence $0\leq n_1 < n_2 < \cdots < n_k <
k/\delta$ such that $\bigcap_{j=1}^k T^{-n_j} (U_{\sigma(j)}
)\neq\emptyset$ for every $\sigma\in\{ 1,2 \}^k$. From
Lemma~\ref{L-density} we see that E-pairs are the same as
nondiagonal IE-pairs.

We now proceed to establish some facts concerning the set of
IE-pairs as captured in Propositions~\ref{P-basic E}, \ref{P-f to
pair E},
and Theorem~\ref{T-product E}. For $\Zb$-systems these can be proved
via a local variational principle route by combining the known
analogues for entropy pairs (see \cite{Disj,ETJ,DEBS}) with Huang
and Ye's characterization (in different terminology) of entropy
pairs as IE-pairs \cite{LVRA}, which itself will be reproved and
extended to cover noninvertible surjective continuous maps in
Theorem~\ref{T-E vs IE} (see also the end of the section for actions
of discrete amenable groups).

Let $k\ge 2$ and let $Z$ be a nonempty finite set. Let $\mathcal{U}$
be the cover of $\{0,1, \dots , k\}^Z=\prod_{z\in Z}\{0,1, \dots ,
k\}$ consisting of subsets of the form $\prod_{z\in
Z}\{i_z\}^\comp$, where $1\le i_z\le k$ for each $z\in Z$. For
$S\subseteq \{0,1,\dots , k\}^Z$ we write $F_S$ to denote the
minimal number of sets in $\mathcal{U}$ one needs to cover $S$.

The following result plays a key role in our combinatorial approach
to the study of $\IE$-tuples (and $\IN$-tuples in
Section~\ref{S-null}). The idea of considering the property (i) in
its proof comes from the proof of Theorem~4 in \cite{MV}.

\begin{lemma}\label{L-key}
Let $k\ge 2$ and let $b>0$ be a constant. There exists a constant
$c>0$ depending only on $k$ and $b$ such that for every finite set
$Z$ and $S\subseteq \{0, 1, \dots , k\}^{Z}$ with $F_S\ge k^{b|Z|}$
there exists a $W\subseteq Z$ with $|W|\ge c|Z|$ and $S|_W\supseteq
\{1, \dots , k\}^W$.
\end{lemma}

\begin{proof}
Pick a constant $0<\lambda<\frac{1}{3}$ such that
$b_1:=b+\log_k(1-\lambda)>0$. Set $b_2:= \log_k \!\big(
\frac{1-\lambda}{\lambda} \big) >0$ and $t=(2b_2 )^{-1} b_1\log_2
\!\big( \frac{k+1}{k} \big)$.

Denote by $H_S$ the number of non-empty subsets $W$ of $Z$ such that
$S|_W\supseteq \{1, \dots , k\}^W$. By Stirling's formula there is a
constant $c>0$ depending only on $t$ (and hence depending only on
$k$ and $b$) such that $\sum_{1\le j\le cn}\binom{n}{j}<2^{tn}$ for
all $n$ large enough. If $H_S\ge 2^{t|Z|}$ and $|Z|$ is large
enough, then there exists a $W\subseteq Z$ for which $|W|\ge c|Z|$
and $S|_W \supseteq \{1, \dots , k\}^W$. It thus suffices to show
that $H_S\ge 2^{t|Z|}$.

Set $S_0=S$ and $Z_0=Z$. We shall construct $Z_0\supseteq
Z_1\supseteq \dots\supseteq Z_m$ for $m:=\lceil
t|Z|/\log_2\frac{k+1}{k} \rceil$ and $S_j\subseteq \{0, 1, \dots ,
k\}^{Z_j}$ for all $1\le j\le m$ with the following properties:
\begin{enumerate}
\item[(i)] $H_{S_{j-1}}\ge \frac{k+1}{k}H_{S_j}$ for all $1\le j\le m$,

\item[(ii)] $F_{S_j}\ge k^{b|Z|}(1-\lambda)^{|Z\setminus Z_j|-j}\lambda^j$ for all
$0\le j\le m$.
\end{enumerate}

Suppose that we have constructed $Z_0, \dots , Z_j$ and $S_0, \dots
, S_j$ with the above properties for some $0\le j<m$. If we have a
$Q\subseteq Z_j$ and a $\sigma\in \{1,\dots , k\}^{Z_j\setminus Q}$
such that $F_{S_{j, \sigma}}\ge (1-\lambda)^{|Z_j\setminus
Q|}F_{S_j}$, where $S_{j, \sigma}$ is the restriction of $\{f\in
S_j:f(x)\neq \sigma(x) \mbox { for all } x\in Z_j\setminus Q\}$ on
$Q$, then
\begin{align*}
|Q|
&\ge \log_k(F_{S_{j, \sigma}}) \\
&\ge (|Z_j\setminus Q|)\log_k(1-\lambda)+\log_k(F_{S_j}) \\
&\ge (|Z_j\setminus Q|+|Z\setminus Z_j|-j)
\log_k(1-\lambda)+b|Z|+j\log_k{\lambda} \\
&= (|Z|-|Q|-j)(b+\log_k(1-\lambda))+b(|Q|+j)+j\log_k{\lambda} \\
&= (|Z|-|Q|)b_1+b|Q|- jb_2 \\
&\ge (|Z|-|Q|)b_1+b|Q|-b_2 \Big( 1+t|Z|/\log_2 \!\Big( \frac{k+1}{k}
\Big)
\Big) \\
&= (|Z|-|Q|)b_1+b|Q|-b_2-\frac{b_1}{2}|Z|
\end{align*}
and hence $|Q|\ge \frac{|Z|b_1/2-b_2}{1+b_1-b}\ge 2$ when $|Z|$ is
large enough. Take $Q$ and $\sigma$ as above such that $|Q|$ is
minimal. Then $|Q|\ge 2$. Pick a $z\in Q$, and set $S_{j,i}$ to be
the restriction of $\{f\in S_{j, \sigma}:f(z)=i\}$ to
$Z_{j+1}:=Q\setminus \{z\}$ for $i=1,\dots , k$. Then
\[F_{S_{j, i}}\ge \lambda(1-\lambda)^{|Z_j\setminus Q|}F_{S_j}\ge
k^{b|Z|}(1-\lambda)^{|Z\setminus Z_{j+1}|-(j+1)}\lambda^{j+1}\] for
$i=1,\dots , k$ (here one needs the fact that $|Q|\ge 2$). Now take
$S_{j+1}$ to be one of the sets among $S_{j, 1}, \dots , S_{j, k}$
with minimal $H$-value, say $S_{j, l}$. For each $1\le i\le k$
denote by $B_i$ the set of nonempty subsets $W\subseteq Z_{j+1}$
such that $S_{j, i}|_W\supseteq \{1,\dots , k\}^W$. Note that
$H_{S_j}\ge |\bigcup^k_{i=1} B_i|+|\bigcap^k_{i=1} B_i|$. If
$|\bigcup^k_{i=1} B_i|\ge \frac{k+1}{k}|B_l|$, then $H_{S_j}\ge
\frac{k+1}{k}|B_l|=\frac{k+1}{k}H_{S_{j+1}}$. Suppose that
$|\bigcup^k_{i=1} B_i|< \frac{k+1}{k}|B_l|$. Note that
\begin{gather*}
\big| {\textstyle\bigcap}^k_{i=1} B_i \big|\cdot k+ \big( \big|
{\textstyle\bigcup}^k_{i=1} B_i \big| - \big|
{\textstyle\bigcap}^k_{i=1} B_i \big| \big) (k-1) \ge
\sum^k_{i=1}|B_i| \ge k|B_l|.
\end{gather*}
Thus
\[ \big| {\textstyle\bigcap}^k_{i=1} B_i \big| \ge k|B_l|-
(k-1)\big| {\textstyle\bigcup}^k_{i=1} B_i \big| \ge
k|B_l|-(k-1)\cdot \frac{k+1}{k}|B_l|=\frac{1}{k}|B_l|.\] Therefore
$H_{S_j}\ge |\bigcup^k_{i=1} B_i|+|\bigcap^k_{i=1} B_i|\ge
|B_l|+\frac{1}{k}|B_l| = \frac{k+1}{k}H_{S_{j+1}}$. Hence the
properties (i) and (ii) are also satisfied for $j+1$.

A simple calculation shows that $k^{b|Z|}(1-\lambda)^{|Z\setminus
Z_m|-m}\lambda^m\ge k^{b|Z|}(1-\lambda)^{|Z|-m}\lambda^m>1$ when
$|Z|$ is large enough. Thus $F_{S_m}>1$ according to property (ii)
and hence $H_{S_m}\ge 1$. By property (i) we have $H_S\ge
(\frac{k+1}{k})^mH_{S_m}\ge (\frac{k+1}{k})^m\ge 2^{t|Z|}$. This
completes the proof of the proposition.
\end{proof}

\noindent For a cover $\cU$ of $X$ we denote by $\hcb(T, \cU)$ the
combinatorial entropy of $\cU$ with respect to $T$, which is defined
using the same formula as for the topological entropy of open
covers.

\begin{lemma}\label{L-E vs IE}
Let $k\ge 2$. Let $U_1, \dots , U_k$ be pairwise disjoint subsets of $X$ and
set $\cU=\{U^{\comp}_1, \dots , U^{\comp}_k\}$. Then $\oU:=(U_1, \dots , U_k)$
has an independence set of positive density if and only if $\hcb(T,
\cU)>0$.
\end{lemma}

\begin{proof} The ``only if'' part is trivial. For the ``if'' part,
set $b:=\hcb(X, \mathcal{U})$ and consider the map $\varphi_n:
X\rightarrow \{0, 1, \dots , k\}^{\{1, \dots, n\}}$ defined by
\[ (\varphi_n(x))(j) = \left\{ \begin{array}{l@{\hspace*{7mm}}l} i,
& \mbox{if } T^{j}(x)\in U_i \mbox{ for some } 1\le i\le k , \\
0, & \mbox{otherwise.} \end{array} \right. \] Then $N\big(
\bigvee^n_{i=1}T^{-i}\mathcal{U} \big) =F_{\varphi_n(X)}$, and so
$F_{\varphi_n(X)}>e^{\frac{b}{2}n}$ for all large enough $n$. By
Lemma~\ref{L-key} there exists a constant $c>0$ depending only on
$k$ and $b$ such that $\varphi_n(X)|_W\supseteq \{1,\dots , k\}^W$
for some $W\subseteq \{1, \dots, n\}$ with $|W|\ge cn$ when $n$ is
large enough. Then $W$ is an independence set for the tuple $\oU$.
Thus by Lemma~\ref{L-density} $\oU$ has an independence set of
positive density.
\end{proof}

We will need the following consequence of Karpovsky and Milman's
generalization of the Sauer-Perles-Shelah lemma
\cite{Sauer,Shelah,KM}. It also follows directly from
Lemma~\ref{L-key}.

\begin{lemma}[\cite{KM}]\label{L-KM}
Given $k\geq 2$ and $\lambda > 1$ there is a constant $c>0$ such
that, for all $n\in \Nb$, if $S\subseteq \{1, 2, \dots , k \}^{\{1,
2, \dots , n\}}$ satisfies $|S|\geq ((k-1)\lambda )^n$ then there is
an $I\subseteq \{1, 2, \dots , n\}$ with $|I|\geq cn$ and $S|_I =
\{1, 2,\dots , k \}^I$.
\end{lemma}

The case $|Z|=1$ of the following lemma appeared in \cite{Pajor83}.

\begin{lemma}\label{L-decomposition E combinatorics}
Let $Z$ be a finite set such that $Z\cap \{1, 2, 3\} = \emptyset$.
There exists a constant $c>0$ depending only on $|Z|$ such that, for
all $n\in \Nb$, if $S\subseteq (Z\cup \{ 1, 2\} )^{\{1, 2, \dots , n
\}}$ is such that $\Gamma_n |_S : S\to (Z\cup \{3\})^{\{1, 2, \dots
, n\}}$ is bijective, where $\Gamma_n : (Z\cup \{ 1, 2\})^{\{ 1, 2,
\dots , n\}}\to (Z\cup \{3\})^{\{1, 2, \dots , n\}}$ converts the
coordinate values $1$ and $2$ to $3$, then there is some $I\subseteq
\{ 1, 2, \dots , n\}$ with $|I|\ge cn$ and either $S|_I \supseteq
(Z\cup \{1\})^I$ or $S|_I\supseteq (Z\cup \{ 2\})^I$.
\end{lemma}

\begin{proof} The case $Z=\emptyset$ is trivial. So we assume that
$Z$ is nonempty. Fix a (small) constant $0<t<\frac{1}{8}$ which we
shall determine later. Denote by $S'$ the elements of $S$ taking
value in $Z$ on at least $(1-4t)n$ many coordinates in $\{1, 2,
\dots , n\}$. Then $|S'|\ge \binom{n}{3tn}|Z|^{(1-4t)n}$ when $n$ is
large enough. Note that each $\sigma\in S'$ takes values $1$ or $2$
on at most $4tn$ many coordinates in $\{1, 2, \dots , n\}$. For
$i=1, 2$, set $S'_i$ to be the elements in $S'$ taking value $i$ on
at most $2tn$ many coordinates in $\{1, 2, \dots , n\}$. Then
$\max(|S'_1|, |S'_2|)\ge \frac{1}{2}|S'|\ge
\frac{1}{2}\binom{n}{3tn}|Z|^{(1-4t)n}$ when $n$ is large enough.
Without loss of generality, we may assume that $|S'_1|\ge
\frac{1}{2}\binom{n}{3tn}|Z|^{(1-4t)n}$. For each $\beta\subseteq
\{1, 2, \dots , n\}$ with $|\beta|\le 2tn$ denote by $S^{\beta}$ the
set of elements in $S'_1$ taking value $1$ exactly on $\beta$. The
number of different $\beta$ is
\[ \sum_{0\le m\le 2tn}\binom{n}{m}\le (2tn+1)\binom{n}{2tn} . \]
By Stirling's formula we can find  $M_1, M_2 > 0$ such that for all
$n\in \Nb$ we have
\[ \binom{n}{2tn}\le \frac{M_1}{\sqrt{tn}}
\bigg( \frac{1}{(1-2t)^{1-2t}(2t)^{2t}} \bigg)^n \] and
\[\binom{n}{3tn}\ge \frac{M_2}{\sqrt{tn}}
\bigg( \frac{1}{(1-3t)^{1-3t}(3t)^{3t}} \bigg)^n . \] Therefore,
when $n$ is large enough we can find some $\beta$ such that
\[|S^{\beta}|\ge
\frac{\frac{1}{2}\binom{n}{3tn}|Z|^{(1-4t)n}}{(2tn+1)\binom{n}{2tn}}
\ge M|Z|^{(1-4t)n}(f(t))^n\frac{1}{2tn+1} , \] where
$M:=\frac{M_2}{2M_1}>0$ and
$f(t):=\frac{(1-2t)^{1-2t}(2t)^{2t}}{(1-3t)^{1-3t}(3t)^{3t}}$. Note
that $\lim_{t\to 0^+}t^{-1} \ln{f(t)} = \infty$. Fix $t$ such that
$f(t)\ge(2|Z|)^{4t}$. Then there is some $n_0 > 0$ such that
\[|S^{\beta}|\ge M (|Z|2^{4t})^n\frac{1}{2tn+1}\ge (|Z|2^t)^n\]
for all $n\ge n_0$. By Lemma~\ref{L-KM} there exists a constant
$c>0$ depending only on $|Z|$ such that for all $n\ge n_0$ we can
find an $I\subseteq \{1, 2, \dots , n\}\setminus \beta$ for which
$|I|\ge c|\{1, 2, \dots , n\}\setminus \beta|\ge c(1-2t)n$ and
$S^{\beta}|_I = (Z\cup \{2\})^I$. Now we may reset $c$ to be
$\min(c(1-2t), 1/n_0)$.
\end{proof}

As an immediate consequence of Lemma~\ref{L-decomposition E
combinatorics} we have:

\begin{lemma}\label{L-decomposition indep}
Let $c$ be as in Lemma~\ref{L-decomposition E combinatorics} for
$|Z|=k-1$. Let $\oA=(A_1, \dots, A_k)$ be a $k$-tuple of subsets of
$X$ and suppose $A_1 = A_{1, 1}\cup A_{1,2}$. If $H$ is a finite
independence set for $\oA$, then there exists some $I\subseteq H$
such that $|I|\ge c|H|$ and $I$ is an independence set for
$(A_{1,1}, \dots, A_k)$ or $(A_{1,2}, \dots, A_k)$.
\end{lemma}

The next lemma follows directly from Lemmas~\ref{L-density} and
\ref{L-decomposition indep}.

\begin{lemma}\label{L-decomposition E}
Let $\oA= (A_1, \dots, A_k )$ be a $k$-tuple of subsets of $X$ which has
an independence set of positive density. Suppose that $A_1 = A_{1,
1}\cup A_{1, 2}$. Then at least one of the tuples $(A_{1,1},\dots,
A_k)$ and $(A_{1,2}, \dots, A_k)$ has an independence set of
positive density.
\end{lemma}

\begin{proposition}\label{P-basic E}
The following are true:
\begin{enumerate}
\item Let $(A_1, \dots , A_k )$ be a tuple of closed subsets of $X$ which has an 
independence set of positive density. Then there exists an IE-tuple $(x_1,\dots ,
x_k)$ with $x_j\in A_j$ for all $1\le j\le k$.

\item $\IE_2(X, T)\setminus \Delta_2(X)$ is nonempty if and only if $\htopol(T)>0$.

\item $\IE_k(X,T)$ is a closed $T\times\cdots\times T$-invariant
subset of $X^k$.

\item Let $\pi:(X, T)\rightarrow (Y,S)$ be a factor map. Then
$(\pi\times\cdots\times \pi )(\IE_k(X, T))=\IE_k(Y, S)$.

\item Suppose that $Z$ is a closed $T$-invariant subset of $X$. Then
$\IE_k(Z, T|_Z )\subseteq \IE_k(X, T)$.
\end{enumerate}
\end{proposition}

\begin{proof}
Assertion (1) follows from Lemma~\ref{L-decomposition E} and a
simple compactness argument. One can easily show that $\htopol(T)>0$
if and only if there is some two-element open cover $\cU=\{U_1,
U_2\}$ of $X$ with positive topological entropy (see for instance
the proof of Proposition~1 in \cite{Disj}). Then assertion (2)
follows directly from assertion (1) and Lemma~\ref{L-E vs IE}.
Assertions (3)--(5) either are trivial or follow directly from
assertion (1).
\end{proof}

We remark that (2) and (4) of Proposition~\ref{P-basic E} show that
the topological Pinsker factor (i.e., the largest zero-entropy factor)
of $(X,T)$ is obtained from the
closed invariant equivalence relation on $X$ generated by the set of
IE-pairs (cf.\ \cite{ZEF}).

In \cite{DEBS} we introduced the notion of $\CA$ entropy for
isometric automorphisms of a Banach space. One can easily extend the
definition to isometric endomorphisms of Banach spaces and check
that Theorem~3.5 in \cite{DEBS} holds in this general setting. In
particular, $\htopol(T)>0$ if and only if there exists an $f\in
C(X)$ with an $\ell_1$-isomorphism set (with respect to the induced
$^*$-endomorphism $f\mapsto f\circ T$) of positive density. Thus if
$f$ is a function in $C(X)$ with an $\ell_1$-isomorphism set of
positive density then the dynamical factor of $(X,T)$ spectrally
generated by $f$ has positive entropy and hence by
Proposition~\ref{P-basic E}(2) has a nondiagonal $\IE$-pair, from
which we infer using Proposition~\ref{P-basic E}(4) that $(X,T)$ has
an $\IE$-pair $(x,y)$ such that $f(x) \neq f(y)$. This yields one
direction of the following proposition. The other direction follows
by a standard argument which appears in the proof of the
Rosenthal-Dor $\ell_1$ theorem \cite{Dor}.

\begin{proposition}\label{P-f to pair E}
Let $f\in C(X)$. Then $f$ has an $\ell_1$-isomorphism set of
positive density if and only if there is an IE-pair $(x, y)$ with
$f(x)\neq f(y)$.
\end{proposition}

We next describe in Proposition~\ref{P-Xm} the set $\IE_1(X, T)$,
which can also be identified with $\IE_k(X, T)\cap \Delta_k(X)$ for
each $k$. The following lemma is mentioned on page 35 of \cite{BGH}.
For completeness we provide a proof here.

\begin{lemma}\label{L-closed positive density}
Let $A$ be a closed subset of $X$. Then $A$ has an independence set
of positive density if and only if there exists a $\mu \in \M(X, T)$
with $\mu(A)>0$.
\end{lemma}

\begin{proof}
Suppose $\mu(A)>0$ for some $\mu \in \M(X, T)$. Let $H$ be a finite
subset of $G$. Denote by $M$ the maximum over all $x\in X$ of the
cardinality of the set of $s\in H$ such that $x\in s^{-1}A$. Then
$\sum_{s\in H}1_{s^{-1}A}\le M$. Thus
\[ M\ge \int \sum_{s\in H}1_{s^{-1}A}\, d\mu=|H|\mu(A).\]
Therefore we can find a subset $I\subseteq H$ such that $|I|\ge
|H|\mu(A)$ and $I$ is an independence set for $A$. By
Lemma~\ref{L-density} $A$ has an independence set of positive
density. This proves the ``if'' part.

Conversely, suppose that $A$ has an independence set $H$ with
density $d>0$. We consider the case $G=\Zb$. The case $G=\Zb_{\ge
0}$ is dealt with similarly. Let $y\in \bigcap_{s\in H}s^{-1}A$.
Then
\[ \liminf_{n\to \infty} \bigg(
\frac{1}{2n+1}\sum^n_{j=-n}\delta_{T^j(y)} \bigg) (A)\ge d. \] Take
an accumulation point $\mu$ of the sequence $\{(2n+1)^{-1}
\sum^n_{j=-n}\delta_{T^j(y)})\}_{n\in \Nb}$ in $\M(X)$. Then $\mu\in
\M(X, T)$ and $\mu(A)\ge d$. This proves the ``only if'' part.
\end{proof}

As a consequence of Lemma~\ref{L-closed positive density} we have:

\begin{proposition}\label{P-Xm}
The set $\IE_1(X, T)$ is the closure of the union of $\supp(\mu)$
over all $\mu \in \M(X, T)$.
\end{proposition}

We describe next in Theorem~\ref{T-product E} the IE-tuples of a
product system. The corresponding statement for entropy pairs was
proved by Glasner in \cite[Theorem 3.(7)(9)]{mu} (see also
\cite[Theorem 19.24]{ETJ}) for a product of metrizable systems using
a local variational principle.

For any tuple $\oA=(A_1, \dots, A_k)$ of subsets of $X$, denote by
$\cP_{\oA}$ the set of all independence sets for $\oA$. Identifying
subsets of $G$ with elements of $\Omega_2:=\{0, 1\}^{G}$ by taking
indicator functions, we may think of $\cP_{\oA}$ as a subset of
$\Omega_2$. Endow $\Omega_2$ with the shift induced from
addition by $1$ on $G$. Clearly $\cP_{\oA}$ is closed and
shift-invariant (i.e., the image of $\cP_{\oA}$ under the shift
coincides with $\cP_{\oA}$). We say a closed shift-invariant subset
$\cP\subseteq \Omega_2$ has {\it positive density} if it has an
element with positive density. Then by definition $\oA$ has an
independence set of positive density exactly when $\cP_{\oA}$ has
positive density. We also say $\cP$ is {\it hereditary} if any
subset of any element in $\cP$ is an element in $\cP$. Note that
$\cP_{\oA}$ is hereditary.

For $s\in G$ we denote by $[s]$ the set of subsets of $G$ containing
$s$.

\begin{lemma}\label{L-positive density subsets}
Let $\cP$ be a closed shift-invariant subset of $\Omega_2$. Then
$\cP$ has positive density if and only if $\mu(\cP\cap [0])>0$ for
some shift-invariant Borel probability measure $\mu $ on $\cP$.
\end{lemma}

\begin{proof}
Notice that $\cP$ has positive density if and only if $\cP\cap [0]$
has an independence set of positive density. The lemma then follows
from Lemma~\ref{L-closed positive density}.
%
\end{proof}

\begin{lemma}\label{L-intersection of positive density}
Let $\cP$ and $\cQ$ be hereditary closed shift-invariant subsets of
$\Omega_2$ with positive density. Then $\cP\cap \cQ$ is also a
hereditary closed shift-invariant subset of $\Omega_2$ with positive
density.
\end{lemma}

\begin{proof}
Clearly $\cP\cap \cQ$ is a hereditary closed shift-invariant subset
of $\Omega_2$. By Lemma~\ref{L-positive density subsets} there is a
shift-invariant Borel probability measure $\mu$ (resp.\ $\nu$) on
$\cP$ (resp.\ $\cQ$) such that $\mu (\cP \cap [0])>0$ (resp.\ $\nu
(\cQ\cap [0])>0$).
Then $(\mu\times \nu )(U)>0$ for $U:=(\cP\cap [0])\times (\cQ\cap
[0])$. By Lemma~\ref{L-closed positive density}, $U$ has an
independence set $H$ of positive density. Take a pair
\[(x', y')\in \bigcap\nolimits_{s\in
H}s^{-1}U= (\cP\times \cQ)\cap \Big( \Big( \bigcap\nolimits_{s\in
H}[s]\Big)\times\Big( \bigcap\nolimits_{s\in H}[s] \Big) \Big) .\]
Then $H\subseteq x', y'$. Since both $\cP$ and $\cQ$ are hereditary,
$H\in \cP\cap \cQ$. This finishes the proof.
\end{proof}

\begin{theorem}\label{T-product E}
Let $(X, T)$ and $(Y, S)$ be dynamical systems. Then
\[ \IE_k(X\times Y, T\times S)=\IE_k(X, T)\times \IE_k(Y, S).\]
\end{theorem}

\begin{proof}
By Proposition~\ref{P-basic E}.(4) the left hand side is included
in the right hand side. The other direction follows from
Lemma~\ref{L-intersection of positive density}.
\end{proof}

We next record the fact that, as an immediate consequence of
Lemma~\ref{L-E vs IE}, nondiagonal $\IE$-tuples are the same as
entropy tuples. This was established for $\Zb$-systems in
\cite{LVRA} using a local variational principle.

\begin{theorem}\label{T-E vs IE}
Let $(x_1, \dots , x_k )$ be a tuple in $X^k\setminus \Delta_k(X)$
with $k\ge 2$. Then $(x_1, \dots , x_k)$ is an entropy tuple if and
only if it is an IE-tuple.
\end{theorem}


In \cite[Corollary 2.4(2)]{BGKM} Blanchard, Glasner, Kolyada, and
Maass showed using the variational principle that positive entropy
implies Li-Yorke chaos. We will next give
an alternative proof of this fact using IE-tuples. The notion of
Li-Yorke chaos was introduced in \cite{BGKM} and is based on ideas
from \cite{LY}. In the case that $X$ is a metric space with metric
$\rho$, a pair $(x_1, x_2)$ of points in $X$ is said to be a {\it
Li-Yorke pair} (with modulus $\delta$) if 
\[ \limsup_{n\to \infty}\rho(T^nx_1, T^nx_2)=\delta>0 \, \, \mbox{ and }
\, \, \liminf_{n\to \infty}\rho(T^nx_1, T^nx_2)=0 . \] A set
$Z\subseteq X$ is said to be {\it scrambled} if all nondiagonal
pairs of points in $Z$ are Li-Yorke. The system $(X, T)$ is said to
be {\it Li-Yorke chaotic} if $X$ contains an uncountable scrambled
set.

We begin by establishing the following lemma, in which we use the
notation established just before Lemma~\ref{L-positive density subsets}. 
For any subset $\cP$ of 
$\Omega_2$, we say that a finite
subset $J\subset G$ has {\it positive density with respect to $\cP$}
if there exists a $K\subseteq G$ with positive density such that
$(K-K)\cap (J-J)=\{0\}$ and $K+J\in \cP$. We say that a subset
$J\subseteq G$ has {\it positive density with respect to $\cP$} if
every finite subset of $J$ has positive density with respect to
$\cP$.

\begin{lemma}\label{L-positive density double}
Let $\cP$ be a hereditary closed shift-invariant subset of
$\Omega_2$ with positive density. Then there exists a $J\subseteq
\Zb_{\ge 0}$ with positive density which also has positive density with
respect to $\cP$.
\end{lemma}

\begin{proof}
Denote by $\cQ$ the set of subsets of $\Zb_{\ge 0}$ which have positive
density with respect to $\cP$. Then $\cQ$ is a hereditary closed
shift-invariant subset of $\{ 0,1 \}^{\Zb_{\ge 0}}$. Take $H\in \cP$
with density $d>0$. Fix $0<d'<\frac{d}{3}$. We claim that if $n$
is large enough so that $0<\frac{d+d'}{2d' n}-\frac{2}{n+2}=:b$
and $d'n\ge 2$, then there exists an $S\in \cQ$ such that
$|S|=c_n:=\lfloor d' n\rfloor$ and $S\subseteq [0, n]$. When $m$ is
large enough, we have $|[-m, m]\cap H|>\frac{d+d'}{2}m$. Arrange the
elements of $[-m, m]\cap H$ in increasing order as $a_1<a_2< \dots
<a_k$, where $k>\frac{d+d'}{2}m$. Consider the numbers
$a_{jc_n}-a_{(j-1)c_n+1}$ for $1\le j\le d_{n, m}:=\Big\lfloor
\frac{\frac{d+d'}{2}m}{c_n}\Big\rfloor$. Set $M_{n, m}=|\{1\le j\le
d_{n, m}: a_{jc_n}-a_{(j-1)c_n+1}\le n\}|.$ Then
\[ 2m+1\ge a_k-(a_1-1)\ge (d_{n, m}-M_{n, m})(n+2).\]
Thus
\[M_{n, m}\ge d_{n, m}-\frac{2m+1}{n+2}\ge
\frac{\frac{d+d'}{2}m}{d' n}-1-\frac{2m+1}{n+2}>\frac{b}{2}m \] when
$m$ is large enough. Note that if $a_{jc_n}-a_{(j-1)c_n+1}\le n$,
then $0<a_{(j-1)c_n+2}-a_{(j-1)c_n+1}<\dots
<a_{jc_n}-a_{(j-1)c_n+1}$ are contained in $[0, n]$. Consequently,
when $m$ is large enough, there exist an $S_m\subseteq [0, n]$ with
$|S_m|=c_n$ and a $W_m\subseteq [-m, m]\cap H$ with $|W_m|\ge
\frac{b}{2\binom{n}{c_n-1}}m$ such that $(W_m-W_m)\cap
(S_m-S_m)=\{0\}$ and $W_m+S_m\subseteq H$. Since $\cP$ is
hereditary, $W_m+S_m\in \cP$. Then there exists an $S\subseteq [0,
n]$ which coincides with infinitely many of the $S_m$. Note that the
collection $\cW_S$ consisting of the sets $W\subseteq G$ such that $(W-W)\cap
(S-S)=\{0\}$ and $W+S\in \cP$ is a closed shift-invariant subset of
$\Omega_2$. By Lemma~\ref{L-density} $\cW_S\cap [0]$ has an
independence set of positive density. Thus $S\in \cQ$.
By Lemma~\ref{L-density} again we see that $\cQ\cap [0]$ has an
independence set of positive density. In other words, $\cQ$ has
positive density. This finishes the proof.
\end{proof}

\begin{theorem}\label{T-LY}
Suppose that $X$ is metrizable with a metric $\rho$.  Suppose that
$k\ge 2$ and $\ox=(x_1, \dots, x_k)$ is an $\IE$-tuple in
$X^k\setminus \Delta_k(X)$. For each $1\le j\le k$, let $A_j$ be a
neighbourhood of $x_j$. Then there exist a $\delta>0$ and a Cantor
set $Z_j\subseteq A_j$ for each $j=1,\dots ,k$ such that the
following hold:
\begin{enumerate}
\item every nonempty finite
tuple of points in $Z:=\bigcup_jZ_j$ is an $\IE$-tuple;

\item for all $m\in \Nb$, distinct $y_1, \dots, y_m \in Z$,
and $y'_1, \dots, y'_m \in Z$ one has
\[ \liminf_{n\to \infty}\max_{1\le i\le m} \rho(T^ny_i, y'_i)=0.\]
\end{enumerate}
\end{theorem}

\begin{proof}
We may assume that the $A_j$ are closed and pairwise disjoint.
We shall construct, via induction on $m$, closed nonempty subsets
$A_{\sigma}$ for $\sigma\in \Sigma_m:=\{1, 2, \dots, k\}^{\{1, 2,
\dots, m\}}$ with the following properties:
\begin{enumerate}
\item[(a)] when $m=1$, $A_\sigma = A_{\sigma (1)}$ for all $\sigma\in \Sigma_m$,

\item[(b)] when $m\ge 2$, $A_{\sigma}\subseteq A_{\sigma |_{\{1, 2, \dots, m-1\}}}$
for all $\sigma\in \Sigma_m$,

\item[(c)] when $m\ge 2$, for every map $\gamma:\Sigma_m\rightarrow \Sigma_{m-1}$ there
exists an $h_{\gamma}\in G$ with $h_\gamma \geq m$ such that
$h_{\gamma}(A_\sigma )\subseteq A_{\gamma(\sigma )}$ for all
$\sigma\in \Sigma_m$,

%

\item[(d)] when $m\ge 2$, $\diam(A_{\sigma})\le 2^{-m}$ for all $\sigma\in \Sigma_m$,

\item[(e)] for every $m$, the collection $\{A_\sigma : \sigma\in \Sigma_m\}$,
ordered into a tuple, has an independence set of positive density.
\end{enumerate}
Suppose that we have constructed such $A_\sigma$ over all $m$. Then
the $A_\sigma$ for all $\sigma$ in a given $\Sigma_m$ are pairwise
disjoint because of property (c). Set $\Sigma = \{1, 2, \dots,
k\}^{\Nb}$. Properties (b) and (d) imply that for each $\sigma \in
\Sigma$ we have $\bigcap_m A_{\sigma |_{\{1, 2, \dots, m\}}}=
\{z_\sigma \}$ for some $z_\sigma \in X$ and that $Z_j=\{z_\sigma
:\sigma\in \Sigma \text{ and } \, \sigma (1)=j\}$ is a Cantor set
for each $j=1,\dots ,k$. Property (a) implies that $Z_j\subseteq
A_j$. Condition (1) follows from properties (d) and (e). Condition
(2) follows from properties (c) and (d).

We now construct the $A_\sigma$. Define $A_\sigma$ for $\sigma\in
\Sigma_1$ according to property (a). By assumption property (e) is
satisfied for $m=1$. Assume that we have constructed $A_\sigma$ for
all $\sigma\in \Sigma_j$ and $j=1,\dots ,m$ with the above
properties. Set $\oA_m$ to be $\{A_\sigma : \sigma\in \Sigma_m\}$
ordered into a tuple. By Lemma~\ref{L-positive density double} there
exist an $H\subseteq \Zb_{\ge 0}$ with positive density which also
has positive density with respect to $\cP_{\oA_m}$. Then for any
nonempty finite subset $J$ of $H$, the sets $A_{J,
\omega}:=\bigcap_{h\in J}h^{-1} A_{\omega (h)}$ for all $\omega\in
(\Sigma_m)^J$, taken together as a tuple, have an independence set
of positive density. Replacing $H$ by $H-h$ for the smallest $h\in
H$ we may assume that $0\in H$ and hence may require that $0\in J$.
For each $\gamma\in (\Sigma_m )^{\Sigma_{m+1}}$ take an
$h_{\gamma}\in J$ with $h_\gamma \geq m+1$. As we can take $|J|$ to
be arbitrarily large, we may assume that $h_{\gamma}\neq
h_{\gamma'}$ for $\gamma\neq \gamma'$ in $(\Sigma_m
)^{\Sigma_{m+1}}$.
Take a map $f:\Sigma_{m+1}\rightarrow (\Sigma_m)^J$ such that
$(f(\sigma ))(0)=\sigma |_{\{1, \dots, m\}}$ and $(f(\sigma
))(h_{\gamma})=\gamma(\sigma )$ for all $\sigma\in \Sigma_{m+1}$ and
$\gamma\in (\Sigma_m )^{\Sigma_{m+1}}$. Set $A_\sigma =A_{J,
f(\sigma )}$ for all $\sigma\in \Sigma_{m+1}$. Then properties (b),
(c), and (e) hold for $m+1$. For each $\sigma\in\Sigma_{m+1}$ write
$A_\sigma$ as the union of finitely many closed subsets each with
diameter no bigger than $2^{-(m+1)}$. Using
Lemma~\ref{L-decomposition E} we may replace $A_\sigma$ by one of
these subsets. Consequently, property (d) is also satisfied for
$m+1$. This completes the induction procedure and hence the proof of
the theorem.
\end{proof}

The set $Z$ in Theorem~\ref{T-LY} is clearly scrambled. As a
consequence of Proposition~\ref{P-basic E}(2) and Theorem~\ref{T-LY}
we obtain the following corollary.

\begin{corollary}\label{C-LY}
Suppose that $X$ is metrizable. If $\htopol(T)>0$, then $(X, T)$ is
Li-Yorke chaotic.
\end{corollary}

As mentioned above, Corollary~\ref{C-LY} was proved in
\cite[Corollary 2.4.(2)]{BGKM} using measure-dynamical techniques.

Denote by $\LY(X, T)$ the set of Li-Yorke pairs in $X\times X$.
Employing a local variational principle, Glasner showed in 
\cite[Theorem 4.(3)]{mu} (see also \cite[Theorem 19.27]{ETJ}) that 
for $\Zb$-systems the set of proximal entropy pairs is dense in the
set of entropy pairs. As a
consequence of Theorem~\ref{T-LY} we have the following improvement:

\begin{corollary}\label{C-proximal}
Suppose that $X$ is metrizable. Then $\LY(X, T)\cap \IE_2(X, T)$ is
dense in $\IE_2(X, T)\setminus \Delta_2(X)$.
\end{corollary}

The next corollary is a direct consequence of Theorem~\ref{T-LY} in
the case $X$ is metrizable and follows from the proof of
Theorem~\ref{T-LY} in the general case.

\begin{corollary}\label{C-no isolated}
Let $W$ be a neighbourhood of an $\IE$-pair $(x_1, x_2)$ in
$X^2\setminus \Delta_2(X)$. Then $W\cap \IE_2(X, T)$ is not of the
form $\{x_1\}\times Y$ for any subset $Y$ of $X$. In particular,
$\IE_2(X, T)$ does not have isolated points.
\end{corollary}

Corollary~\ref{C-no isolated}, as stated for entropy pairs, was
proved by Blanchard, Glasner, and Host in \cite[Theorem 6]{BGH}.

We round out this section by explaining how the facts about IE-tuples
captured in Propositions~\ref{P-basic E}, \ref{P-f to pair E}, and \ref{P-Xm} and
Theorems~\ref{T-product E} and \ref{T-E vs IE}
apply to actions of any discrete amenable group. So for the remainder of the
section $G$ will be a discrete amenable group. 
For a finite $K\subseteq G$ and $\delta>0$ we denote by $M(K, \delta)$ 
the set of all nonempty finite subsets $F$ of $G$ which are $(K, \delta)$-invariant
in the sense that  
\[ |\{s\in F: Ks\subseteq F\}|\ge (1-\delta )|F| . \]
According to the F{\o}lner characterization of 
amenability, $M(K, \delta)$ is nonempty for every finite set $K\subseteq G$ and 
$\delta > 0$. This is equivalent to the existence of a 
{\it F{\o}lner net}, i.e., a net $\{ F_\gamma \}_\gamma$ of nonempty finite 
subsets of $G$ such that 
$\lim_\gamma | sF_\gamma \Delta F_\gamma | / | F_\gamma | = 0$ for all $s\in G$.
For countable $G$ we may take this net to be a sequence,
in which case we speak of a {\it F{\o}lner sequence}. 

For countable $G$,
a sequence $\{ F_n \}_{n\in\Nb}$ of nonempty finite subsets of $G$ is said
to be {\it tempered} if for some $c > 0$ we have
$| \bigcup_{k=1}^{n-1} F_k^{-1} F_n | \leq c | F_n |$ for all $n\in\Nb$. 
By \cite[Prop.\ 1.4]{PET}, every F{\o}lner sequence has a tempered subsequence.
Below we will make use of
the pointwise ergodic theorem of Lindenstrauss \cite{PET}, which applies to 
tempered F{\o}lner sequences.

The following subadditivity 
result was established by Lindenstrauss and Weiss for countable 
$G$ \cite[Theorem 6.1]{LW}. We will reduce the general case to their result.
Note that there is also a version for left invariance.

\begin{proposition}\label{P-subadditive} 
If $\varphi$ is a real-valued function 
which is defined on the set of finite subsets of $G$ and satisfies 
\begin{enumerate} 
\item $0\le \varphi(A)<+\infty$ and $\varphi(\emptyset)=0$,  

\item $\varphi(A)\le \varphi(B)$ for all $A\subseteq B$,  

\item $\varphi(As)=\varphi(A)$ for all finite $A\subseteq G$ and $s\in G$,   

\item $\varphi(A\cup B)\le \varphi(A)+\varphi(B)$ if $A\cap B=\emptyset$,  
\end{enumerate} 
then $\frac{1}{|F|} \varphi(F)$ converges to some limit $b$ 
as the set $F$ becomes  
more and more invariant in the sense that 
for every $\varepsilon>0$ there exist a 
finite set $K\subseteq G$ and a $\delta>0$ such that 
$\big| \frac{1}{|F|} \varphi(F) -b \big| <\varepsilon$ 
for all $(K, \delta)$-invariant 
sets $F\subseteq G$.  
\end{proposition}

\begin{proof} 
For a finite $K\subseteq G$ and $\delta>0$ we
denote by $Z(K, \delta)$ the closure of 
$\big\{ \frac{1}{|F|}\varphi(F): F\in M(K, \delta)\big\}$. 
The collection of sets $Z(K, \delta)$ has the finite intersection property 
and clearly the conclusion of the proposition is equivalent to the set
$Z:=\bigcap_{(K, \delta)}Z(K, \delta)$ containing only one point. 
Suppose that $x$ and $y$ are 
distinct points of $Z$. Then one can find a sequence $\{ (A_n, B_n) \}_{n\in\Nb}$ of 
pairs of finite subsets of $G$ such that $\{A_n\}_{n\in \Nb}$ and 
$\{B_n\}_{n\in \Nb}$ are both F{\o}lner sequences for the subgroup $H$ of $G$ 
generated by $\bigcup_{n\in \Nb} (A_n\cup B_n)$ and   
$\max\! \big( \big| x-\frac{1}{|A_n|}\varphi(A_n)\big| , 
\big| y-\frac{1}{|B_n|}\varphi(B_n) \big| \big) <1/n$   
for each $n\in \Nb$. Since subgroups of discrete amenable groups are amenable  
\cite[Proposition 0.16]{Paterson}, $H$ is amenable. This contradicts the 
result of Lindenstrauss and Weiss. Thus $Z$ contains only one point.  
\end{proof} 

Let $(X, G)$ be a dynamical system. We define the topological entropy of $(X,G)$
by first defining the topological entropy of a finite open cover of $X$ using
Proposition~\ref{P-subadditive} and then taking a supremum
over all finite open covers
(this was originally introduced in \cite{JMO} without the
subadditivity result).
For a finite tuple $\oA=(A_1, \dots , A_k)$ 
of subsets of $X$, Proposition~\ref{P-subadditive} also applies to the 
function $\varphi_{\oA}$ given by 
\[ \varphi_{\oA}(F)=\max\{|F\cap J|: J \mbox{ is an 
independence set for } \oA\} . \] 
This permits us to define the 
{\it independence density} $I(\oA )$ of $\oA$ as the limit of 
$\frac{1}{|F|}\varphi_{\oA}(F)$ as $F$ becomes more and more invariant,
providing a numerical measure of the dynamical independence of $\oA$.

\begin{proposition}\label{P-ps} 
Let $(X, G)$ be a dynamical system. Let 
$\oA=(A_1, \dots , A_k)$ be a tuple of subsets of $X$. Let $c>0$. Then the following 
are equivalent: 
\begin{enumerate} 
\item $I(\oA ) \geq c$,

\item for every $\varepsilon>0$ there exist a finite set $K\subseteq G$ and a 
$\delta>0$ such that for every $F\in M(K, \delta)$ there is an independence set $J$ 
for $\oA$ with $|J\cap F|\ge (c-\varepsilon)|F|$.  

\item for every finite set
$K\subseteq G$ and $\varepsilon>0$ there exist an $F\in M(K, \varepsilon)$ 
and an independence set $J$ for $\oA$ such that $|J\cap F|\ge (c-\varepsilon)|F|$.  
\end{enumerate} 
When $G$ is countable, these conditions are also equivalent to: 
\begin{enumerate} 
\item[4.] for every tempered F{\o}lner sequence $\{F_n\}_{n\in \Nb}$ of $G$ there  
is an independence set $J$ for $\oA$ such that  
$\lim_{n\to \infty}\frac{|F_n\cap J|}{|F_n|}\ge c$.  

\item[5.] there are a tempered F{\o}lner sequence $\{F_n\}_{n\in \Nb}$ of $G$ 
and an independence set $J$ for $\oA$ such that  
$\lim_{n\to \infty}\frac{|F_n\cap J|}{|F_n|}\ge c$. 
\end{enumerate} 
\end{proposition}  

\begin{proof} 
The equivalences (1)$\Leftrightarrow$(2)$\Leftrightarrow$(3) follow from 
Proposition~\ref{P-subadditive}. Assume now that $G$ is countable. 
Then the implications (4)$\Rightarrow$(5)$\Rightarrow$(3) are trivial.  
Suppose that (3) holds and let us show (4). Then one 
can easily show that there is a $G$-invariant Borel probability measure $\mu$ on 
$\cP_{\oA}\subseteq \{0, 1\}^G$ with $\mu([e]\cap \cP_{\oA})\ge c$, as in the proof 
of Lemma~\ref{L-closed positive density}. Here $\cP_{\oA}$ is defined as before
Lemma~\ref{L-positive density subsets} and
$\{0, 1\}^G$ is 
equipped with the shift given by $sx(t)=x(ts)$ 
for all $x\in \{0, 1\}^G$ and $s, t\in G$.  Replacing $\mu$ by a suitable ergodic 
$G$-invariant Borel probability measure in the ergodic decomposition of $\mu$, we 
may assume that $\mu$ is ergodic. Let $\{F_n\}_{n\in \Nb}$ be a tempered F{\o}lner 
sequence for $G$. The pointwise ergodic theorem \cite[Theorem 1.2]{PET} asserts 
that $\lim_{n\to \infty} \frac{1}{| F_n |}\sum_{s\in F_n} f(sx)=\int f\, d\mu$ 
$\mu$-a.e.\ for every $f\in L^1(\mu)$. Setting $f$ to be the characteristic 
function of $[e]\cap \cP_{\oA}$ and taking $J$ to be some $x$ satisfying the 
above equation, we get (4). 
\end{proof} 

Effectively extending Definition~\ref{D-IE-pair},
we call a tuple $\ox = (x_1 , \dots, x_k )\in X^k$ an {\it IE-tuple}
(or an {\it IE-pair} in the case $k=2$) if for every tuple
$\oU = (U_1,\dots, U_k)$ associated to a product
neighbourhood $U_1 \times\cdots\times U_k$ of $\ox$
the independence density $I(\oU )$ is nonzero. Denoting
the set of $\IE$-tuples of length $k$ by $\IE_k (X,G)$ and replacing everywhere
the existence of positive density independence sets for a tuple $\oA$ by
the nonvanishing of the independence density $I(\oA )$
in our earlier discussion for singly generated systems, we see 
that Propositions~\ref{P-basic E} and \ref{P-Xm} 
and Theorem~\ref{T-product E} 
continue 
to hold in our current setting. 
In particular, the topological Pinsker factor (i.e., the largest 
zero-entropy factor) of $(X,G)$ arises from the closed invariant 
equivalence relation on $X$ generated by the set of IE-pairs.

Given an action $\alpha$ of $G$ by isometric automorphisms of a 
Banach space $V$ and an element $v\in V$, we can apply the left invariance version
of Proposition~\ref{P-subadditive} to the function $\varphi_{v,\lambda}$ given by 
\[ \varphi_{v,\lambda} (F)=\max\{|F\cap J|: J \mbox{ is an }\ell_1 \mbox{-}
\lambda \mbox{-isomorphism set for } v\} \] 
and define the {\it $\ell_1$-$\lambda$-isomorphism density} $I(v,\lambda )$ of $v$ 
as the limit of $\frac{1}{|F|}\varphi_{v,\lambda} (F)$ as $F$ becomes 
more and more invariant.
Defining the CA entropy of $\alpha$ by taking a limit supremum of averages 
as in Section~2 of \cite{DEBS} but this time along a F{\o}lner net, one can
check that the analogue of Theorem~3.5 of \cite{DEBS} holds. In particular,
the topological entropy of $(X,G)$ is nonzero if and only if there exists an $f\in
C(X)$ with nonvanishing $\ell_1$-$\lambda$-isomorphism density
for some $\lambda\geq 1$.
Then we obtain the analogue of Proposition~\ref{P-f to pair E}, i.e.,
a function $f\in C(X)$ has nonvanishing $\ell_1$-$\lambda$-isomorphism density
for some $\lambda\geq 1$ if and only if there is an $\IE$-pair $(x,y)$ in
$X\times X$ with $f(x) \neq f(y)$.

Finally, if we define entropy tuples in the same way as for singly
generated systems, then Theorem~\ref{T-E vs IE} still holds.


\section{Topological entropy and tensor product
independence}\label{S-entropy tensor}

In this section we will see how combinatorial independence in the context of
entropy translates into the
language of tensor products, with a hint at how the theory might
thereby be extended to noncommutative $C^*$-dynamical systems.
As in the previous section, our dynamical systems here will have
acting semigroup $\Zb$ or $\Zb_{\geq 0}$ with generating endomorphism $T$.

To start with, we remark that, given a dynamical system $(X,T)$ and
denoting by $\alpha_T$ the induced $^*$-endomorphism $f\mapsto f\circ T$ of
$C(X)$, the following conditions are equivalent:
\begin{enumerate}
\item $\htopol (T) > 0$,

\item $\IE_2 (X,T)\setminus\Delta_2 (X) \neq \emptyset$,

\item there is an $f\in C(X)$ with an $\ell_1$-isomorphism set of positive
density,

\item there is a 2-dimensional operator subsystem $V\subseteq C(X)$ and a
$\lambda\geq 1$ such that $V$ has a $\lambda$-independence set of positive density
for $\alpha_T$.
\end{enumerate}
The equivalence (1)$\Leftrightarrow$(3) was established in \cite{DEBS}
for $\Zb$-systems and is similarly seen to be valid for endomorphisms,
and (1)$\Leftrightarrow$(2) is Proposition~\ref{P-basic E}(2)
of Section~\ref{S-entropy comb}. To show (2)$\Rightarrow$(4) we simply need
to take a pair $(A,B)$ of disjoint closed subsets of $X$ with an independence
set of positive density and consider the operator system generated by
a norm-one self-adjoint function in $C(X)$ taking the constant values $1$
and $-1$ on $A$ and $B$, respectively. Finally, the implication
(4)$\Rightarrow$(3) is a consequence of the following lemma, which
expresses in a form suited to our context a well-known phenomenon
observed by Rosenthal in the proof of his $\ell_1$ theorem \cite{Ros,ell1},
namely that $\ell_1$ geometry ensues in a natural way from independence.

\begin{lemma}\label{L-tensor-indep}
Let $V$ be an operator system and let $v$ be a nonscalar 
element of $V$. For each $j\in\Nb$ set $v_j = 1\otimes v\otimes 1
\in V^{\otimes [1,j-1]} \otimes
V\otimes V^{\otimes [j+1,\infty )} = V^{\otimes\Nb}$. 
Set $Z = \{ \sigma (v) : \sigma\in S(V) \}$ where $S(V)$ is the state
space of $V$. Let $0 < \eta < \frac12 \diam (Z)$. Then
for all $n\in\Nb$ and complex scalars $c_1 , \dots , c_n$ we have
\[ \frac{\eta}{4} \sum_{j=1}^n | c_j | \leq
\bigg\| \sum_{j=1}^n c_j v_j \bigg\| . \]
\end{lemma}

\begin{proof}
Choose points $b_1 , b_2 \in Z$ such that 
$| b_1 - b_2 | = \diam (Z)$, and take disks $D_1$ and $D_2$ in the
complex plane centred at $b_1$ and $b_2$, respectively, with
common diameter $d$ and at distance greater than
$\max (2\eta ,2d)$ from each other. For each $j\in\Nb$ define the
subsets
\begin{align*}
U_j &= \{ \sigma\in S(V^{\otimes\Nb} ) : \sigma (v_j ) \in D_1 \} ,\\
V_j &= \{ \sigma\in S(V^{\otimes\Nb} ) : \sigma (v_j ) \in D_2 \}
\end{align*}
of the state space $S(V^{\otimes\Nb} )$.
Then the collection of pairs $(U_j , V_j )$ for 
$j\in\Nb$ is independent, and so we obtain the result by the proof
in \cite{Dor}.
\end{proof}

What can be said about the link between $\ell_1$ structure and
independence in connection with the global picture of
entropy production?
This question is tied to the divergence in topological dynamics between
the notions of completely positive entropy and uniformly positive
entropy.

Recall that the dynamical system $(X,T)$ is said to have
{\it completely positive entropy} or {\it c.p.e.} if each of its nontrivial
factors has positive topological entropy, i.e., if its Pinsker
factor is trivial \cite{UPE}. The functions in
$C(X)$ which lie in the topological Pinsker algebra (the $C^*$-algebraic
manifestation of the Pinsker factor) are characterized by the fact that
they lack an $\ell_1$-isomorphism set of positive density for the induced
$^*$-endomorphism of $C(X)$ \cite{DEBS}. Thus $(X,T)$ has
c.p.e.\ precisely when every nonscalar function in $C(X)$ has
an $\ell_1$-isomorphism set of positive density.

The system $(X,T)$ is said to have {\it uniformly positive entropy} or
{\it u.p.e.} if every nondiagonal pair in
$X\times X$ is an entropy pair \cite{UPE}.
More generally, following \cite{TEE} we say that $(X,T)$ has
{\it u.p.e.\ of order $n$} if every tuple in $X^n \setminus\Delta_n (X)$ is an
entropy tuple (see the beginning of Section~\ref{S-entropy comb}).
By Theorem~\ref{T-E vs IE}, this is equivalent to every $n$-tuple of nonempty
open subsets of $X$ having an independence set of positive density.
Finally, we say that $(X,T)$ has {\it u.p.e.\ of all orders} if it has
u.p.e.\ of order $n$ for each $n\geq 2$. U.p.e.\ implies c.p.e., but the
converse is false. Also, in \cite{LVRA} it is shown that,
for every $n\geq 2$, u.p.e.\ of order $n$ does not imply u.p.e.\ of order
$n+1$.

The following propositions supply functional-analytic characterizations
of u.p.e.\ and u.p.e.\ of all orders, complementing the
characterization of c.p.e.\ in terms of $\ell_1$-isomorphism sets.


\begin{proposition}\label{P-equiv-D-indep}
Let $X$ be a compact Hausdorff space and $T:X\to X$ a surjective continuous map.
Then the following are equivalent:
\begin{enumerate}
\item $(X,T)$ has u.p.e.\ of all orders,

\item for every finite set $\Omega\subseteq C(X)$ and $\delta > 0$
there is a finite-dimensional operator subsystem
$V\subseteq C(X)$ which approximately includes $\Omega$ to within $\delta$
and has a $1$-independence set of positive density,

\item for every finite set $\Omega\subseteq C(X)$ and $\delta > 0$
there is a finite-dimensional operator subsystem
$V\subseteq C(X)$ which approximately includes $\Omega$ to within $\delta$
and has a $\lambda$-independence set of positive density for some
$\lambda\geq 1$.
\end{enumerate}
\end{proposition}

\begin{proof}
(1)$\Rightarrow$(2). Let $\Omega$ be a finite subset of $A$ and let
$\delta > 0$. Then we can construct a partition of unity $\{ g_1 , \dots , g_k \}$
in $C(X)$ such that the operator system it spans approximately includes
$\Omega$ to within $\delta$ and for each $i=1,\dots ,k$ we have $g_i (x) = 1$
for all $x$ in some nonempty open set $U_i$. By (1) and Theorem~\ref{T-E vs IE},
the tuple $(U_1 , \dots , U_k )$ has an independence set $I$ of positive density.
Let $(s_1 , \dots , s_n )$ be a tuple of distinct elements of $I$ and let
$\varphi : V^{\otimes [1,n]} \to C(X)$ be the
contractive linear map determined on elementary tensors by
$f_1 \otimes\cdots\otimes f_n \mapsto
(f_1 \circ T^{s_1} ) \cdots (f_n \circ T^{s_n} )$. Then the collection
$\{ (g_{\sigma (1)} \circ T^{s_1} ) \cdots 
(g_{\sigma (n)} \circ T^{s_n} ) : \sigma\in \{ 1,\dots ,k\}^{\{ 1,\dots ,n\}} \}$ is
an effective $k^n$-element partition of unity of $X$ and hence is isometrically
equivalent to the standard basis of $\ell_\infty^{k^n}$. Since the subset
$\{ g_{\sigma (1)} \otimes\cdots\otimes g_{\sigma (n)}
: \sigma\in \{ 1,\dots ,k\}^{\{ 1,\dots ,n\}} \}$
of $V^{\otimes [1,n]}$ is also isometrically equivalent
to the standard basis of $\ell_\infty^{k^n}$, we conclude that $I$ is
a $1$-independence set for $V$, yielding (2).

(2)$\Rightarrow$(3). Trivial.

(3)$\Rightarrow$(1). Let $k\geq 2$ and let $(U_1 , \dots , U_k )$ be a $k$-tuple of
nonempty open subsets of $X$ with pairwise disjoint closures. To obtain (1)
it suffices to show that $(U_1 , \dots , U_k )$ has an independence set of
positive density. Since the
sets $U_i$ have pairwise disjoint closures we can
construct positive norm-one functions $g_1 , \dots , g_k \in C(X)$ such
that, for each $i$, the set on which $g_i$ takes the value $1$ is a closed subset of
$U_i$. By (3), given a $\delta > 0$ we can
find an operator subsystem $V\subseteq C(X)$ such that $V$ has
a $\lambda$-independence set $J$ of positive density for some $\lambda\geq 1$
and there exist
$f_1 , \dots , f_k \in V$ for which $\| f_i - g_i \| < \delta$ for
each $i=1, \dots ,k$. Since for each $i$ continuity implies that
$g_i^{-1} ((\rho , 1]) \subseteq U_i$ for some $\rho\in (0,1)$, by taking
$\delta$ small enough we may ensure that the norm-one self-adjoint elements
$h_i = \| (f_i + f_i^* )/2 \|^{-1} (f_i + f_i^* )/2 \in V$ for $i=1, \dots , k$
are defined and for some $\theta\in (0,1)$ satisfy
$h_i^{-1} ((\theta ,1]) \subseteq U_i$ for each $i$.

Choose an $r\in\Nb$ large enough so that $\theta^r < \lambda^{-1}$ and a
$b>0$ small enough so that $\frac{k}{k-1} 2^{-2b} > 1$. By Stirling's
formula there is a $c\in (0,1/2)$ such that $\binom{n}{cn} \leq 2^{bn}$
for all $n\in\Nb$. By Lemma~\ref{L-KM} there is a $d>0$ such that,
for all $n\in\Nb$, if $\Gamma\subseteq\{ 1,\dots ,k \}^{\{ 1,\dots ,n \}}$ and
$| \Gamma | \geq (k2^{-2b} )^n = \big( (k-1)\big( \frac{k}{k-1} 2^{-2b} \big)
\big)^n$ then
there exists a set $I\subseteq \{ 1,\dots , n\}$ such that $|I| \geq dn$ and
$\Gamma |_I = \{ 1,\dots ,k \}^I$.

By shifting $J$ we may assume that it contains $0$, and so by positive density
there is an $a>0$ such that
for each $m\in\Zb$ the set $J_m := J \cap \{ -m, -m+1 , \dots , m\}$
has cardinality at least $am$.
Now suppose we are given an $m\in\Nb$ with $m\geq \max (r/ab , r/ac)$.
Enumerate the elements of $J_m$ as $j_1 , \dots , j_n$.
For each $\sigma\in\{ 1,\dots ,k \}^{\{ 1,\dots ,n \}}$ we define a set
$K_\sigma \subseteq I$ as follows. Pick an $x\in X$ such
that $h_{\sigma (1)} (T^{j_1} x)\cdots h_{\sigma (n)} (T^{j_n} x) =
\| (h_{\sigma (1)} \circ T^{j_1} )\cdots
(h_{\sigma (n)} \circ T^{j_n} ) \|$. Since
$\| (h_{\sigma (1)} \circ T^{j_1} )\cdots
(h_{\sigma (n)} \circ T^{j_n} ) \| \geq \lambda^{-1} > \theta^r$, there exists
a $K_\sigma \subseteq \{ 1, \dots , n\}$ with
$|K_\sigma | = n-r$ such that
$h_{\sigma (i)} (T^{j_i} x) > \theta$ for each $i\in K_\sigma$.
Hence $x\in\bigcap_{i\in K_\sigma} T^{-j_i} U_{\sigma (i)}$, so
that $\bigcap_{i\in K_\sigma} T^{-j_i} U_{\sigma (i)}$ is nonempty.

Now since $r\leq cam \leq cn$, we have
\[ \big| \big\{ K_\sigma : \sigma\in\{ 1,\dots ,k \}^{\{ 1,\dots ,n \}}
\big\} \big| \leq \binom{n}{r} \leq 2^{bn} . \]
We can thus find a
$K\in \big\{ K_\sigma : \sigma\in \{ 1,\dots ,k \}^{\{ 1,\dots ,n \}} \big\}$
such that the set $\mathcal{R}$ of all
$\sigma\in\{ 1,\dots ,k \}^{\{ 1,\dots ,n \}}$ for
which $K_\sigma = K$ has cardinality at least
$k^n / 2^{bn}$. Then the set of all
restrictions of elements of $\mathcal{R}$ to $K$ has cardinality
at least $| \mathcal{R} | / 2^{n-|K|} \geq 2^{-r} (k2^{-b} )^n
\geq (k2^{-2b} )^n$.
It follows that there is a set $I_m \subseteq J_m$ with
$|I_m | \geq d|J_m | \geq dam$
such that the set of all restrictions of elements of $\mathcal{R}$
to $I_m$ is $\{ 1,\dots ,k \}^{I_m}$. Since
$\bigcap_{i\in I_m} T^{-j_i} U_{\sigma (i)} \neq\emptyset$ for every
$\sigma\in\{ 1,\dots ,k \}^{I_m}$, $I_m$ is an independence set for
$(U_1 , \dots , U_k )$. By Lemma~\ref{L-density} we conclude that
$(U_1 , \dots , U_k )$ has an independence set of positive density, finishing
the proof.
\end{proof}

We denote by $\cS_2 (X)$
the collection of $2$-dimensional operator
subsystems of $C(X)$ equipped with the metric given by the Hausdorff
distance between unit balls.

\begin{proposition}\label{P-upe}
For a $\Zb$-dynamical system $(X,T)$ the following are
equivalent:
\begin{enumerate}
\item $(X,T)$ has u.p.e.,

\item the collection of $2$-dimensional operator
subsystems of $C(X)$ which have a
$1$-independence set of positive density is dense in $\cS_2 (X)$,

\item the collection of $2$-dimensional operator
subsystems of $C(X)$ which have a
$\lambda$-independence set of positive density for some
$\lambda\geq 1$ is dense in $\cS_2 (X)$.
\end{enumerate}
\end{proposition}

\begin{proof}
(1)$\Rightarrow$(2). Let $V$ be a
$2$-dimensional operator subsystem of $C(X)$. Then $V$ has a linear basis
of the form $\{ 1,g \}$ for some nonscalar $g\in C(X)$, which we may
assume to be self-adjoint by replacing it with $g+g^*$ if necessary. By
scaling and scalar-translating $g$, we may furthermore assume that the
spectrum of $g$ is a subset of $[0,1]$ containing $0$ and $1$.
Given $\delta > 0$, by a simple perturbation argument we can construct a positive
norm-one $h\in C(X)$ such that $\| h-g \| < \delta$, $h(x) = 0$ for all
$x$ in some nonempty open set $U_0$, and $h(x) = 1$ for all $x$ in some
nonempty open set
$U_1$. By taking $\delta$ small enough we can make the unit ball of the
operator system $W = \spn \{ 1,h \}$ as close as we wish to the unit ball of $V$,
and so we will obtain (2)
once we show that $W$ has a $1$-independence set of positive density, and this can
be done as in the proof of the corresponding implication in
Proposition~\ref{P-equiv-D-indep} using the partition of unity $\{ h , 1-h \}$.

(2)$\Rightarrow$(3). Trivial.

(3)$\Rightarrow$(1). Argue as in the proof of the corresponding implication in
Proposition~\ref{P-equiv-D-indep}.
\end{proof}

In \cite{LVRA} Huang and Ye proposed u.p.e.\ of all orders as the most suitable
topological analogue of a K-system. In view of
Proposition~\ref{P-equiv-D-indep} we might then define a unital $^*$-endomorphism
$\alpha$ of a unital $C^*$-algebra $A$ to be a $C^*$-algebraic K-system if for every
finite set $\Omega\subseteq A$ and $\delta > 0$ there is a finite-dimensional
operator subsystem $V\subseteq A$ which approximately includes $\Omega$ to within
$\delta$ and has a $(1+\varepsilon )$-independence set of positive density for
every $\varepsilon > 0$. This property holds prototypically for the shift on the
infinite minimal tensor product $A^{\otimes\Zb}$ for any unital $C^*$-algebra $A$,
and it implies completely positive
Voiculescu-Brown entropy \cite{toral}, as can be see from Remark~3.10
of \cite{EID} and Lemma~\ref{L-tensor-indep}. In fact to deduce completely positive
Voiculescu-Brown entropy all we need is that the collection
of $2$-dimensional operator subsystems of $A$ which have a $\lambda$-independence
set of positive density for some $\lambda\geq 1$ is dense in the 
collection of all $2$-dimensional
operator subsystems of $A$ with respect to the metric given by Hausdorff
distance between unit balls. 

\begin{remark}
Following up on the discussion in the last part of Section~\ref{S-entropy comb},
we point out that the results of this section hold more generally for any
dynamical system $(X,G)$ with $G$ discrete and amenable 
if the relevant terms are interpreted or reformulated as follows. 
U.p.e., u.p.e\ of order $n$, and u.p.e.\ of all orders are defined in the 
same way as for singly generated systems. 
For a finite-dimensional operator subsystem $V\subseteq C(X)$
and a $\lambda\geq 1$, the left invariance version of 
Proposition~\ref{P-subadditive} applies to the 
function $\varphi_{V,\lambda}$ given by 
$\varphi_{V,\lambda} (F)=\max\{|F\cap J|: J \mbox{ is a }
\lambda\mbox{-independence set for } V\}$, so that we may
define the {\it $\lambda$-independence density} $I(V,\lambda )$ of $V$ as 
the limit of $\frac{1}{|F|}\varphi_{V,\lambda} (F)$ 
as $F$ becomes more and more invariant.
We can then replace ``$\lambda$-independence set of positive density'' 
everywhere above by ``nonvanishing $\lambda$-independence density''.
\end{remark}


\section{Topological sequence entropy, nullness, and independence}\label{S-null}

In this section we will examine the local theory of topological
sequence entropy and nullness from the viewpoint of independence.
Sequence entropy is developed in the literature for single
continuous maps, but here we will work in the framework of a general
dynamical system $(X,G)$.

Following Goodman \cite{SE}, for a sequence $\mathfrak{s} = \{ s_n
\}_{n\in\Nb}$ in $G$ we define the topological sequence entropy of
$(X,G)$ with respect to $\mathfrak{s}$ and a finite open cover $\cU$
of $X$ by
\[ \htopol (X,\cU ; \mathfrak{s} ) = \limsup_{n\to\infty} \frac1n \log
N\bigg( \bigvee_{i=1}^n s_i^{-1} \cU \bigg) \] where $N(\cdot )$
denotes the minimal cardinality of a subcover. Following Huang, Li,
Shao, and Ye \cite{NSSEP}, we call a nondiagonal pair $(x,y)\in
X\times X$ a {\it sequence entropy pair} if for any disjoint closed
neighbourhoods $U$ and $V$ of $x$ and $y$, respectively, there
exists a sequence $\mathfrak{s}$ in $G$ such that $\htopol(X, \{U^{\comp},
V^{\comp} \};\mathfrak{s} )>0$.  More generally, following 
Huang, Maass, and Ye \cite{HMY} we call a
tuple $\ox = (x_1 , \dots, x_k)\in X^k\setminus \Delta_k(X)$ a {\it sequence
entropy tuple} if whenever $U_1, \dots , U_l$ are closed pairwise
disjoint neighbourhoods of the distinct points in the list 
$x_1 , \dots , x_k$, the open cover
$\{ U_1^\comp , \dots , U_l^\comp \}$ has positive topological sequence
entropy with respect to some sequence in $G$. 
We say that $(X, G)$ is {\it null} if $\htopol (X,
\cU; \mathfrak{s})=0$ for all open covers $\cU$ of $X$ and all
sequences $\mathfrak{s}$ in $G$. We say that $(X, G)$ is {\it
nonnull} otherwise. Then the basic facts recorded in
\cite[Proposition 2.1]{NSSEP} and \cite[Proposition 3.2]{HMY}
also hold in our general setting.


\begin{definition}\label{D-IN-pair}
We call a tuple $\ox = (x_1 , \dots, x_k )\in X^k$ an {\it IN-tuple}
(or an {\it IN-pair} in the case $k=2$) if for any product
neighbourhood $U_1 \times\cdots\times U_k$ of $\ox$ the tuple $(U_1,
\dots, U_k)$ has arbitrarily large finite independence sets. We
denote the set of $\IN$-tuples of length $k$ by $\IN_k (X,G)$.
\end{definition}

We will show in Proposition~\ref{P-SE vs IN} that sequence entropy tuples are
exactly nondiagonal $\IN$-tuples.

First we record some basic facts pertaining to IN-tuples.
For a cover $\cU$ of $X$ and a sequence $\mathfrak{s}$ in $G$ we
denote by $\hcb(T, \cU; \mathfrak{s})$ the combinatorial sequence
entropy of $\cU$ with respect to $\mathfrak{s}$, which is defined
using the same formula as for the topological sequence entropy of
open covers. From Lemma~\ref{L-key} we infer the following analogue
of Lemma~\ref{L-E vs IE}.

\begin{lemma}\label{L-SE vs IN}
Let $k\ge 2$. Let $U_1, \dots , U_k$ be disjoint subsets of $X$ and
set $\cU=\{U^{\comp}_1, \dots , U^{\comp}_k\}$. Then $\oU:=(U_1, \dots , U_k)$
has arbitrarily large finite independence sets if and only if\linebreak
$\hcb(T, \cU; \mathfrak{s})>0$ for some sequence $\mathfrak{s}$ in
$G$.
\end{lemma}

From Lemma~\ref{L-decomposition indep} we infer the following
analogue of Lemma~\ref{L-decomposition E}.

\begin{lemma}\label{L-decomposition N}
Let $\oA= (A_1, \dots, A_k )$ be a $k$-tuple of subsets of $X$ with
arbitrarily large finite independence sets. Suppose that $A_1 =
A_{1, 1}\cup A_{1, 2}$. Then at least one of the tuples
$(A_{1,1},\dots, A_k)$ and $(A_{1,2}, \dots, A_k)$ has arbitrarily
large finite independence sets.
\end{lemma}

Using Lemmas~\ref{L-SE vs IN} and ~\ref{L-decomposition N} we obtain
the following analogue of Proposition~\ref{P-basic E}.

\begin{proposition}\label{P-basic N}
The following are true:
\begin{enumerate}
\item Let $(A_1, \dots , A_k )$ be a tuple of closed subsets of $X$ which has
arbitrarily large finite independence
sets. Then there exists an IN-tuple $(x_1,\dots , x_k)$ with $x_j\in
A_j$ for all $1\le j\le k$.

\item $\IN_2(X, G)\setminus \Delta_2(X)$ is nonempty if and only if $(X,
G)$ is nonnull.

\item $\IN_k(X,G)$ is a closed $G$-invariant
subset of $X^k$.

\item Let $\pi:(X, G)\rightarrow (Y, G)$ be a factor map. Then
$(\pi\times\cdots\times \pi )(\IN_k(X, G))=\IN_k(Y, G)$.

\item Suppose that $Z$ is a closed $G$-invariant subset of $X$. Then
$\IN_k(Z, G)\subseteq \IN_k(X, G)$.
\end{enumerate}
\end{proposition}
%
%
%
%
%

We remark that (2) and (4) of Proposition~\ref{P-basic N} show that
the largest null factor of $(X,G)$ is obtained from the closed
invariant equivalence relation on $X$ generated by the set of
IN-pairs.

\begin{corollary}\label{C-permanence N}
For a fixed $G$, the class of null $G$-systems is preserved under
taking factors, subsystems and products.
\end{corollary}

\begin{corollary}\label{C-se=1}
Suppose that $(X,G)$ is nonnull. Then there exist an open cover
$\mathcal{U}=\{U, V\}$ of $X$ and a sequence $\mathfrak{s}$ in $G$
such that $\htopol (X, \mathcal{U};\mathfrak{s} )=\log{2}$. In the
case $G=\Zb$ the sequence can be taken in $\Nb$.
\end{corollary}

\begin{definition}\label{D-null}
We say that a function $f\in C(X)$ is {\it null} if there does not
exist a $\lambda\geq 1$ such that $f$ has arbitrarily large finite
$\ell_1$-$\lambda$-isomorphism sets. 
Otherwise $f$ is said to be {\it nonnull}.
\end{definition}

The proof of the following proposition is similar to that of
Proposition~\ref{P-f to pair E}, where this time we use
Theorem~5.8 in \cite{DEBS} as formulated for the more general
context of isometric endomorphisms of Banach spaces.

\begin{proposition}\label{P-f to pair N}
Let $f\in C(X)$. Then $f$ is nonnull if and only if there is an
IN-pair $(x, y)$ with $f(x)\neq f(y)$.
\end{proposition}

In analogy with the situation for entropy, nondiagonal IN-tuples
turn out to be the same as sequence entropy tuples, as follows from
Lemma~\ref{L-SE vs IN}.

\begin{theorem}\label{P-SE vs IN}
Let $(x_1, \dots , x_k )$ be a tuple in $X^k\setminus \Delta_k(X)$
with $k\ge 2$. Then $(x_1, \dots , x_k)$ is a sequence entropy tuple
if and only if it is an IN-tuple.
\end{theorem}


In parallel with Blanchard's definition of uniform positive entropy
\cite{UPE}, we say that $(X, G)$ is {\it uniformly nonnull} if
$\IN_2(X, G)=X\times X$, or, equivalently, if any pair of nonempty
open subsets of $X$ has arbitrarily large finite independence sets.
By Theorem~\ref{P-SE vs IN} this is the same as s.u.p.e.\ as defined 
on page 1507 of \cite{NSSEP}.
We say that $(X, G)$ is {\it completely nonnull} if the maximal null
factor of $(X, G)$ is trivial. Uniform nonnullness implies complete
nonnullness, but the converse is false. Blanchard showed in
\cite[Example 8]{UPE} that the shift action on $X=\{a, b\}^{\Zb}\cup
\{a, c\}^{\Zb}$ has completely positive entropy but is not
transitive. Hence this action is completely nonnull but by
Theorem~\ref{T-equiv-indep-prod} fails to be uniformly nonnull.

We now briefly examine the tensor product viewpoint. In analogy with
the case of entropy (see the beginning of Section~\ref{S-entropy
tensor}), it can be shown that the system $(X,G)$ is nonnull if and
only if there is a 2-dimensional operator subsystem $V\subseteq
C(X)$ and a $\lambda\geq 1$ such that $V$ has arbitrarily large
finite $\lambda$-independence sets for the induced $C^*$-dynamical system.
Using arguments similar to those in the proof of
Proposition~\ref{P-upe} (with $\cS_2 (X)$ defined as in the
discussion there) one can show:

\begin{proposition}\label{P-unn}
For a dynamical system $(X,G)$ the following are equivalent:
\begin{enumerate}
\item $(X,G)$ is uniformly nonnull,

\item the collection of $2$-dimensional operator
subsystems of $C(X)$ which have arbitrarily large finite $1$-independence
sets is dense in $\cS_2 (X)$,


\item for every $2$-dimensional operator subsystem $V\subseteq C(X)$
there is a $\lambda\geq 1$ such that $V$ has arbitrarily large
finite $\lambda$-independence sets.
\end{enumerate}
\end{proposition}

For $n\geq 2$ we say that $(X, G)$ is {\it uniformly nonnull of
order $n$} if $\IN_n(X, G)=X^n$, or, equivalently, if every
$n$-tuple of nonempty open subsets of $X$ has arbitrarily large
finite independence sets. We say that $(X, G)$ is {\it uniformly
nonnull of all orders} if it is uniformly nonnull of order $n$ for
every $n\geq 2$. The analogue of Proposition~\ref{P-unn} for uniform
nonnullness of all orders is also valid and will be subsumed as part of
Proposition~\ref{P-equiv-indep} and Theorem~\ref{T-equiv-indep-prod} 
in connection with I-independence, to which uniform nonnullness of 
all orders is equivalent.


\section{Tameness and independence}\label{S-tame}

Let $(X, G)$ be a dynamical system. We
say that $(X, G)$ is {\it tame} if no element $f\in C(X)$ has an
infinite $\ell_1$-isomorphism set and {\it untame} otherwise. The concept 
of tameness was introduced by K\"{o}hler in \cite{K} under the term regularity. 
Here we are following the terminology of \cite{tame}. Actually in \cite{tame} 
the system $(X, G)$ is defined to be tame if its enveloping semigroup is separable 
and Fr\'{e}chet, which is equivalent to our geometric definition when 
$X$ is metrizable.

\begin{definition}\label{D-IT-pair}
We call a tuple $\ox = (x_1 , \dots, x_k )\in X^k$ an {\it IT-tuple}
(or an {\it IT-pair} in the case $k=2$) if for any product
neighbourhood $U_1 \times\cdots\times U_k$ of $\ox$ the tuple $(U_1,
\dots, U_k)$ has an infinite independence set. We denote the set of
$\IT$-tuples of length $k$ by $\IT_k (X,G)$.
\end{definition}
%

In contrast to the density conditions in the context
of entropy and sequence entropy, we are interested
here in the existence of infinite sets along which independence (or,
in the definition of tameness, equivalence to the standard
$\ell_1$ basis) occurs. This places us in the realm of Rosenthal's
$\ell_1$ theorem \cite{Ros} and the Ramsey methods involved in its
proof (see \cite{Gow}). Indeed we begin by observing the following
immediate consequence of \cite[Theorem 2.2]{Ros}. Note that though
Theorem 2.2 of \cite{Ros} is stated for sequences of pairs of
disjoint subsets, the proof there works for sequences of tuples of
(not necessarily disjoint) subsets.

\begin{lemma}\label{L-infinite indep}
Let $\oA=(A_1, \dots , A_k)$ be a tuple of closed
subsets of $X$. Then the pair $\oA$ has an infinite independence set
if and only if there is an infinite set $H\subseteq G$ such that
for every infinite set $W\subseteq H$ there exists an $x\in X$ for
which the sets $(x, A_j, W)^{\perp}:=\{s\in W: sx\in A_j\}$ for
$j=1, \dots , k$ are all infinite.
\end{lemma}

\begin{lemma}\label{L-decomposition R}
Let $\oA= (A_1, \dots, A_k )$ be a $k$-tuple of 
closed subsets of $X$ with an infinite independence set. Suppose
that $A_1 = A_{1, 1}\cup A_{1, 2}$ and that both $A_{1, 1}$ and $A_{1,
2}$ are closed. Then at least one of the tuples $(A_{1,1}, \dots,
A_k)$ and $(A_{1,2}, \dots, A_k)$ has an infinite independence set.
\end{lemma}

\begin{proof} Let $H\subseteq G$
be as in Lemma~\ref{L-infinite indep} for $\oA$. Suppose that
neither of $(A_{1,1}, \dots, A_k)$ and $(A_{1, 2}, \dots, A_k)$ has
an infinite independence set. Then by Lemma~\ref{L-infinite indep}
we can find infinite subsets $W_1\supseteq W_2$ of $H$ such that for
every $x\in X$ and $i=1,2$ at least one of the sets $(x, A_{1, i},
W_i)^{\perp}$ and $(x, A_j, W_i)^{\perp}$ for $2\le j\le k$ is
finite. But there exists an $x\in X$ such that $(x, A_j,
W_2)^{\perp}$ is infinite for all $1\le j\le k$. Thus $(x, A_{1, 1},
W_1)^{\perp}$ and $(x, A_{1, 2}, W_2)^{\perp}$ are finite. Since
$(x, A_1, W_2)^{\perp}\subseteq (x, A_{1, 1}, W_2)^{\perp}\cup (x,
A_{1, 2}, W_2)^{\perp}$, we obtain a contradiction. Therefore at
least one of the tuples $(A_{1,1}, \dots, A_k)$ and $(A_{1,2},
\dots, A_k)$ has an infinite independence set.
\end{proof}

\begin{proposition}\label{P-basic R}
The following are true:
\begin{enumerate}
\item Let $(A_1, \dots , A_k )$ be a tuple of closed subsets of $X$ 
which has an infinite independence set.
Then there exists an IT-tuple $(x_1,\dots , x_k)$ with $x_j\in A_j$
for all $1\le j\le k$.

\item $\IT_2(X, T)\setminus \Delta_2(X)$ is nonempty if and only if $(X, G)$ is untame.

\item $\IT_k(X,G)$ is a closed $G$-invariant
subset of $X^k$.

\item Let $\pi:(X, G)\rightarrow (Y,G)$ be a factor map. Then
$(\pi\times\cdots\times \pi )(\IT_k(X, G))=\IT_k(Y, G)$.

\item Suppose that $Z$ is a closed $G$-invariant subset of $X$. Then
$\IT_k(Z, G)\subseteq \IT_k(X, G)$.
\end{enumerate}
\end{proposition}
%
%
%
%

The proof of Proposition~\ref{P-basic R} is similar to that of
Proposition~\ref{P-basic E}, with (1) following from
Lemma~\ref{L-decomposition R} and (2) following from
Proposition~\ref{P-f to pair R} below.

We remark that (2) and (4) of Proposition~\ref{P-basic R} show that
the largest tame factor of $(X, G)$ is obtained from the closed
invariant equivalence relation on $X$ generated by the set of
IT-pairs.

\begin{corollary}\label{C-permanence R}
For a fixed $G$, the class of tame $G$-systems is preserved under
taking factors, subsystems, and products.
\end{corollary}

\begin{proposition}\label{P-f to pair R}
Let $f\in C(X)$. Then $f$ has an infinite $\ell_1$-isomorphism set
if and only if there is an $\IT$-pair $(x, y)$ with $f(x)\neq f(y)$.
\end{proposition}

\begin{proof}
The ``if'' part follows as in the analogous situation for entropy by
the well-known Rosenthal-Dor argument. For the ``only if'' part, by
Proposition~\ref{P-basic R}(1) it suffices to show  the existence of
a pair of disjoint closed subsets $A$ and $B$ of $X$ which have an
infinite independence set and satisfy $f(A)\cap f(B)=\emptyset$. Let
$H=\{s_j:j\in \Nb\}\subseteq G^\op$ be an $\ell_1$-isomorphism set
for $f$. Then the sequence $\{s_j f \}_{j\in \Nb}$ has no weakly
convergent subsequence. 
Using Lebesgue's theorem one sees that
$\{s_j f \}_{j\in \Nb}$ has no pointwise convergent subsequence. In
Gowers' proof of Rosenthal's $\ell_1$ theorem
\cite[page 1079-1080]{Gow} it is shown that
there exist disjoint closed subsets $Z_1, Z_2\subseteq \Cb$ and a
subsequence $\{s_{n_k} f\}_{k\in \Nb}$ such that the sequence of
pairs $\{((s_{n_k} f)^{-1}(Z_1), (s_{n_k} f)^{-1}(Z_2))\}_{k\in
\Nb}$ is independent. Therefore $\{s_{n_k} : k\in \Nb\}$ is an
infinite independence set for the pair $(f^{-1}(Z_1) ,
f^{-1}(Z_2))$, yielding the proposition.
\end{proof}


In parallel with the cases of entropy and nullness, we say that $(X,
G)$ is {\it uniformly untame} if $\IT_2(X, G)=X\times X$, or,
equivalently, if any pair of nonempty open subsets of $X$ 
has an infinite independence set. We say that $(X, G)$ is
{\it completely untame} if the maximal tame factor of $(X,
G)$ is trivial. Complete untameness is strictly weaker than
uniform untameness, as illustrated by Blanchard's example
mentioned in the paragraph following Proposition~\ref{P-SE vs IN}.

We demonstrate in the following example that untame systems need
not be null by constructing a $\WAP$ (weakly almost periodic)
nonnull subshift. The proof of \cite[Corollary 5.7]{DEBS} shows that
if $(X, G)$ is $\HNS$ (hereditarily nonsensitive 
\cite[Definition 9.1]{GM}) then it
is tame. Since $\WAP$ systems are $\HNS$ \cite[Section 9]{GM}, 
our example is tame.

\begin{example}\label{E-tame nonnull}
Let $0<m_1<m_2<\cdots$ be a sequence in $\Nb$ with $m_j-m_i>m_i-m_k$
for all $j>i>k$. Let $\{k_n\}_{n\in \Nb}$ be an unbounded sequence
in $\Nb$ and let $S_j=\sum^j_{i=1}k_i$ for $j\in\Nb$ be the partial
sums of $\{k_n\}_{n\in \Nb}$. Denote by $A_j$ the set of all
elements in $\{0, 1\}^{\Zb}$ whose support is contained in
$\{m_k:S_{j-1}<k\le S_j\}$. One checks easily that the union $X$ of
the orbits of $\bigcup_j A_j$ under the shift $T$ is closed. Thus
$(X, \Zb)$ is a subshift. Denote by $Z$ the set of elements in $X$
supported exactly at one point. It is easy to see that $\IN_2(X,
\Zb)\setminus \Delta_2(X)=(Z\times \{0\})\cup (\{0\}\times Z)$. Since the set of elements
in $C(X)$ whose $\Zb$-orbit is precompact in the weak topology of
$C(X)$ is a closed $\Zb$-invariant algebra of $C(X)$, to see that $(X,
\Zb)$ is $\WAP$ it suffices to check that the orbit of $f$ is
precompact in the weak topology, where $f\in C(X)$ is defined by
$f(x)=x(0)$ for $x\in X$. However, the union of the zero function
and the orbit of $f$ is evidently compact in $\Cb^X$ and thus is
compact in the weak topology by a result of Grothendieck \cite{Gr} 
\cite[Theorem 1.43.1]{ETJ}.
\end{example}

In Section~\ref{S-Toeplitz} we will see that there exist minimal
tame nonnull $\Zb$-systems.

Turning finally to tensor products, in analogy with entropy and
sequence entropy one can show that the system $(X,G)$ is untame
if and only if there is a 2-dimensional operator subsystem
$V\subseteq C(X)$ with an infinite independence set for the induced
$C^*$-dynamical system. With $\cS_2 (X)$ defined as in the
paragraph preceding Proposition~\ref{P-upe}, one has the following
characterizations of uniform untameness (see the proof of
Proposition~\ref{P-equiv-indep}).

\begin{proposition}\label{P-unr}
For a dynamical system $(X,G)$ the following are equivalent:
\begin{enumerate}
\item $(X,G)$ is uniformly untame,

\item the collection of $2$-dimensional operator
subsystems of $C(X)$ which have an infinite $1$-independence set is
dense in $\cS_2 (X)$,


\item the collection of $2$-dimensional operator
subsystems of $C(X)$ which have an infinite $\lambda$-independence
set for some $\lambda\geq 1$ is dense in $\cS_2 (X)$.
\end{enumerate}
\end{proposition}

For $n\geq 2$ we say that $(X,G)$ is {\it uniformly untame of
order $n$} if $\IT_n(X, G)=X^n$, or, equivalently, if every
$n$-tuple of nonempty open subsets of $X$ has an infinite
independence set. If $(X,G)$ is uniformly untame of order $n$
for each $n\geq 2$ then we say that it is {\it uniformly untame
of all orders}. For the analogue of Proposition~\ref{P-unr} for
uniform untameness of all orders see Proposition~\ref{P-equiv-indep} 
and Theorem~\ref{T-equiv-indep-prod}.


\section{Tame extensions of minimal systems}\label{S-minimal}

An extension $\pi:X\rightarrow Y$ of dynamical systems with acting
group $G$ is said to be {\it tame} if
$\IT_{\pi}\setminus \Delta_2(X)=\emptyset$, 
where $R_{\pi}:=\{(x_1, x_2)\in X\times X:\pi(x_1)=\pi(x_2)\}$
and $\IT_{\pi}:=\IT_2(X, G)\cap R_{\pi}$.
In this section we will analyze the structure of tame extensions of minimal
systems.
Throughout $G$ will be a group, and we will frequently
specialize to the Abelian case, to which the main results
(Theorems~\ref{T-metric tame minimal to highly proximal},
\ref{T-decomposition}, and \ref{T-uniquely ergodic}) apply.
Theorems~\ref{T-metric tame minimal to highly proximal}
and \ref{T-decomposition} address the relation between proximality and
equicontinuity within the frame of tame extensions, while
Theorem~\ref{T-uniquely ergodic} asserts that tame minimal
systems are uniquely ergodic. 

We refer the reader to \cite{Aus,Vri} for
general information on extensions of dynamical systems. 

As a direct consequence of Proposition~\ref{P-basic R}(4) we have:

\begin{proposition}\label{P-basic R extension}
The following are true:
\begin{enumerate}
\item Let $\psi:X\rightarrow Y$ and $\varphi:Y\rightarrow Z$ be extensions.
Then $\varphi\circ \psi$ is tame if and only if both $\varphi$ and $\psi$ are tame.

\item Let $\pi_j:X_j\rightarrow Y_j$ for $j\in J$ be extensions. 
If every $\pi_j$ is tame, then
$\prod_{j\in J} \pi_j:\prod_{j\in J}X_j\rightarrow \prod_{j\in J}Y_j$
is tame.

\item Let $\{\pi_{\alpha}:X_{\alpha}\rightarrow Y \}_{\alpha<\nu}$
be an inverse system of extensions, where $\nu$ is an ordinal.
Then $\pi:=\varprojlim \pi_{\alpha}$ is tame if and only if
every $\pi_{\alpha}$ is tame.

\item For every commutative diagram
\begin{gather*}
\xymatrix{
{X'} \ar[d]_{\pi'} \ar[r]^{\sigma} & {X}\ar[d]^{\pi} \\
Y'  \ar [r]^{\tau} & Y
}
\end{gather*}
of extensions, if $X'$ is a subsystem of $Y'\times X$, $\pi'$ and $\sigma$ are the
restrictions to $X'$ of the coordinate projections of $Y'\times X$ onto $Y'$ and $X$,
respectively, and $\pi$ is tame, then $\pi'$ is tame.
\end{enumerate}
\end{proposition}

A continuous map $f:A\rightarrow B$ between topological spaces is said to be 
{\it semi-open}
if the image of every nonempty open subset of $A$ under $f$ has nonempty interior.
For a dynamical system $(X, G)$, denote by $\RP(X, G)$ the {\it regionally proximal relation} of $X$,
that is, $\RP(X, G)=\bigcap_{U\in \mathcal{U}_{X}}\overline{GU}$, where
$\mathcal{U}_{X}$ is the collection of open neighbourhoods
of the diagonal $\Delta_2(X)$ in $X^2$.
Following \cite{NSSEP} we call a pair $(x_1,x_2)\in X^2\setminus\Delta_2(X)$
a {\it weakly mixing pair} if for any neighbourhoods $U_1$ and $U_2$ of
$x_1$ and $x_2$, respectively,
there exists an $s\in G$ such that $sU_1\cap U_1\neq \emptyset$ and
$sU_1\cap U_2\neq \emptyset$. We denote the set of weakly mixing pairs by $\WM(X, G)$.
When $(X, G)$ is minimal, $\WM(X, G)=\RP(X, G)\setminus \Delta_2(X)$ 
\cite[Theorem 2.4(1)]{NSSEP}.
The proofs of Lemma~3.1 and Corollary~4.1 in \cite{NSSEP} yield the following two results.

\begin{lemma}\label{L-semiopen to IT}
Suppose that $G$ is Abelian. Let $(x_1, x_2)\in X^2\setminus \Delta_2(X)$ and
$A=\overline{G(x_1, x_2)}$. If $\pi_1:A\rightarrow X$ is semi-open, where $\pi_1$ is the
projection to the first coordinate, then the following are equivalent:
\begin{enumerate}
\item $(x_1, x_2)\in \IT_2(X, G)$,

\item $(x_1, x_2)\in \IN_2(X, G)$,

\item $(x_1, x_2)\in \WM(X, G)$.
\end{enumerate}
\end{lemma}

\begin{lemma}\label{L-distal ext}
Suppose that $G$ is Abelian.
Let $\pi:(X, G)\rightarrow (Y, G)$ be a distal extension of minimal systems.
Let $(x_1, x_2)\in R_{\pi} \setminus \Delta_2(X)$. Then the following are equivalent:
\begin{enumerate}
\item $(x_1, x_2)\in \IT_2(X, G)$,

\item $(x_1, x_2)\in \IN_2(X, G)$,

\item $(x_1, x_2)\in \WM(X, G)$,

\item $(x_1, x_2)\in \RP(X, G)$.
\end{enumerate}
\end{lemma}

We call an extension $\pi:X\rightarrow Y$ a {\it Bron{\v s}tein extension} 
if $R_{\pi}$ has a
dense subset of almost periodic points.
Recall that $\pi$ is said to be {\it highly proximal} if
for every nonempty open subset $U$ of $X$ and every point $z\in Y$ there exists an
$s\in G$ such that $\pi^{-1}(z)\subseteq sU$. 
In this case $(X, G)$ and $(Y, G)$ are necessarily
minimal and $\pi$ is proximal. We say that
$\pi$ is {\it strictly PI} if it can be obtained
by a transfinite succession of proximal and equicontinuous extensions,
and {\it PI} if there exists a proximal extension $\psi:X'\rightarrow X$
such that $\pi\circ \psi$ is strictly PI.

If in these definitions proximality is replaced by high proximality, then we obtain the notions
of {\it strictly HPI} extensions and {\it HPI} extensions.

\begin{lemma}\label{L-B+tame}
Suppose that $G$ is Abelian. Let $\pi:X\rightarrow Y$ be a tame Bron{\v s}tein extension of minimal
systems. Suppose that $\pi$ has no nontrivial equicontinuous factors. Then $\pi$ is an isomorphism.
\end{lemma}

\begin{proof}
For an extension $\psi:X'\rightarrow Y'$
denote by $\RP_{\psi}$ the {\it relative regionally proximal relation} of $\psi$,
that is, $\RP_{\psi}=\bigcap_{U\in \mathcal{U}_{X'}}\overline{GU\cap R_{\psi}}$, where
$\mathcal{U}_{X'}$ is the collection of open neighbourhoods
of the diagonal $\Delta_2(X')$ in $X'\times X'$.
As $X$ is minimal, if we denote by $S_{\pi}$ the smallest closed 
$G$-invariant equivalence relation on $X$ containing $\RP_{\pi}$,
then the induced extension $X/S_{\pi}\rightarrow Y$ is an equicontinuous factor of $\pi$ \cite{Ell}
\cite[Theorem V.2.21]{Vri}.
Since $\pi$ has no nontrivial equicontinuous factors, $S_{\pi}=R_{\pi}$.
It is a theorem of Ellis that for any Bron{\v s}tein extension $\psi$ of minimal systems, $\RP_{\psi}$
is an equivalence relation \cite[Theorem 2.6.2]{Veech} \cite[Theorem VI.3.20]{Vri}.
Therefore $R_{\pi}=S_{\pi}=\RP_{\pi}\subseteq \RP(X, G)$.
Since $X$ is minimal, $\WM(X, G)=\RP(X, G)\setminus \Delta_2(X)$ 
\cite[Theorem 2.4(1)]{NSSEP}.
Thus $R_{\pi}\setminus \Delta_2(X) \subseteq \WM(X, G)$. 
Suppose that $\pi$ is not an isomorphism.
Since $\pi$ is a Bron{\v s}tein extension, 
we can find an almost periodic point $(x_1, x_2)$ in
the nonempty open subset $R_{\pi}\setminus \Delta_2(X)$ of $R_{\pi}$.
As extensions of minimal systems are semi-open, we conclude that 
$(x_1, x_2)\in \IT_2(X, G)$
by Lemma~\ref{L-semiopen to IT}. This is in contradiction to the tameness of $\pi$. Thus
$\pi$ is an isomorphism.
\end{proof}

\begin{lemma}\label{L-tame to PI}
Suppose that $G$ is Abelian. Then any tame extension $\pi:X\rightarrow Y$ of minimal
systems is PI.
\end{lemma}

\begin{proof} 
Consider the canonical PI tower of $\pi$ \cite{EGS} \cite[Theorem VI.4.20]{Vri}.
This is a commuting diagram of extensions of minimal systems of the form
displayed in Proposition~\ref{P-basic R extension}(4),
where $\sigma$ is proximal, $\tau$ is strictly PI,
$\pi'$ is RIC (relatively incontractible), 
and $\pi'$ has no nontrivial equicontinuous factors.
Furthermore, $X'$ is a subsystem of $Y'\times X$, and
$\pi'$ and $\sigma$ are the restrictions to
$X'$ of the coordinate projections of $Y'\times X$ onto $Y'$ and $X$, respectively
(see for example \cite[VI.4.22]{Vri}).
Since $\pi$ is tame, by Proposition~\ref{P-basic R extension}(4) so is $\pi'$.
As every RIC-extension is a Bron{\v s}tein extension
\cite[Corollary 5.12]{EGS} \cite[Corollary VI.2.8]{Vri}, $\pi'$ is a Bron{\v s}tein extension.
Then $\pi'$ is an isomorphism by Lemma~\ref{L-B+tame}. Hence $\pi$ is PI.
\end{proof}

Lemma~\ref{L-tame to PI} generalizes a result of Glasner \cite[Theorem 2.3]{tame}, who
proved it for tame metrizable $(X, G)$.

A theorem of van der Woude asserts that an extension $\pi:X\rightarrow Y$ of minimal systems is HPI
if and only if every topologically transitive subsystem $W$ of
$R_{\pi}$ for which both of the coordinate projections $\pi_i:W\rightarrow X$ are semi-open
is minimal \cite[Theorem 4.8]{HPI}.

\begin{lemma}\label{L-tame metrizable to HPI}
Suppose that $G$ is Abelian. Then any tame extension $\pi:X\rightarrow Y$ of minimal metrizable systems
is HPI.
\end{lemma}

\begin{proof}
Let $S_{\pi}$ be as in the proof of Lemma~\ref{L-B+tame}.
Then we have the natural extensions $\psi:X\rightarrow X/S_{\pi}$ and
$\varphi: X/S_{\pi}\rightarrow Y$ with $\pi=\varphi\circ \psi$, and
$\varphi$ is equicontinuous. To show that $\pi$ is HPI, it suffices to show
that $\psi$ is HPI.
Let $W$ be a topologically transitive subsystem of $R_{\psi}=S_{\pi}$ for which both of the
coordinate projections $\pi_i:W\rightarrow X$ are semi-open. 
By the theorem of van der Woude,
the proof will be complete once we show that $W$ is minimal.
Since $G$ is Abelian, $X$ has a $G$-invariant Borel probability measure by
a result of Markov and Kakutani \cite[Theorem VII.2.1]{Dav}.
Then $\RP(X, G)$ is an equivalence relation by \cite[Theorem 9.8]{Aus} 
\cite[Theorem V.1.17]{Vri}.
Thus $S_{\pi}\subseteq \RP(X, G)$.
Since $X$ is minimal, $\WM(X, G)=\RP(X, G)\setminus \Delta_2(X)$ 
\cite[Theorem 2.4(1)]{NSSEP}.
Then $R_{\psi}\setminus \Delta_2(X)=S_{\pi}\setminus \Delta_2(X) \subseteq \WM(X, G)$.
As $X$ is metrizable and $W$ is topologically transitive, we can find a transitive point
$(x_1, x_2)$ of $W$. Since $\pi$ is tame,
by Lemma~\ref{L-semiopen to IT} we must have $x_1=x_2$. 
Therefore $W=\Delta_2(X)$ is minimal.
\end{proof}

\begin{lemma}\label{L-open tame to HPI}
Suppose that $G$ is Abelian. Let $\pi:X\rightarrow Y$ be an open
tame extension of minimal systems. If $Y$ is metrizable, then $\pi$
is HPI.
\end{lemma}

To prove Lemma~\ref{L-open tame to HPI} we need a variation of
the ``reduction-to-the-metric-case construction'' in \cite{MW}, which is in turn
a relativization of a method of Ellis \cite{Ellis}.

Let $\pi:(X, G)\rightarrow (Y, G)$ be an open extension of $G$-systems.
Denote by $\CP(X)$ the set of all continuous pseudometrics on $X$.
For a
$\rho\in \CP(X)$
and a countable subgroup
$H$ of $G$, define $C_{\rho, H}:=\{(x_1, x_2)\in X\times X : \rho(sx_1, sx_2)=0 \mbox{ for all } s\in H\}$.
Then $C_{\rho, H}$ is a closed $H$-invariant equivalence relation on $X$.
Denote $X/C_{\rho, H}$ by $X_{\rho, H}$, and denote the quotient map $X\rightarrow X_{\rho, H}$ by
$\psi_{\rho, H}$. Say $H=\{s_1, s_2, \dots\}$ with $s_1 = e$.
Define a pseudometric $d_{\rho, H}$ on $X$ by
\begin{gather}\label{E-d_H}
d_{\rho, H}(x_1, x_2):=\sum^{\infty}_{i=1}2^{-i}\rho(s_ix_1, s_ix_2)
\end{gather}
for all $x_1, x_2\in X$.
Then $X_{\rho, H}$ is a metrizable space with a metric $d'_{\rho, H}$ defined by
\begin{gather}\label{E-d'_H}
d'_{\rho, H}(\psi_{\rho, H}(x_1), \psi_{\rho, H}(x_2))=d_{\rho, H}(x_1, x_2)
\end{gather}
for all $x_1, x_2\in X$.
Let $C'_{\rho, H}=\{(y_1, y_2)\in Y\times Y : \psi_{\rho, H}(\pi^{-1}(y_1))
=\psi_{\rho, H}(\pi^{-1}(y_2))\}$.
Then $C'_{\rho, H}$ is a closed $H$-invariant equivalence relation on $Y$.
Denote $Y/C'_{\rho, H}$ by $Y_{\rho, H}$, and denote the quotient map $Y\rightarrow Y_{\rho, H}$ by
$\tau_{\rho, H}$. Define a map $\sigma_{\rho, H}:X\rightarrow  X_{\rho, H}\times Y_{\rho, H}$
by $\sigma_{\rho, H}(x)=(\psi_{\rho, H}(x), \tau_{\rho, H}\circ \pi(x))$, and write $\sigma_{\rho, H}(X)$
as $X^*_{\rho, H}$. Denote by $\pi_{\rho, H}$ the restriction
to $X^*_{\rho, H}$ of the projection of $X_{\rho, H}\times Y_{\rho, H}$ onto $Y_{\rho, H}$. 

\begin{lemma}\cite[Lemma V.3.3]{Vri}\label{L-open to quotient}
For every
$\rho\in \CP(X)$
and countable subgroup
$H$ of $G$ the diagram
\begin{gather}\label{D-open to quotient}
\xymatrix{
{(X, H)} \ar[d]_{\sigma_{\rho, H}} \ar[r]^{\pi} & {(Y, H)}\ar[d]^{\tau_{\rho, H}} \\
(X^*_{\rho, H}, H)  \ar [r]^{\pi_{\rho, H}} & (Y_{\rho, H}, H)
}
\end{gather}
of extensions of $H$-systems commutes and $X^*_{\rho, H}$ is metrizable.
\end{lemma}

For $\rho_1, \rho_2\in \CP(X)$ we write $\rho_1\preceq \rho_2$
if $\rho_1(x_1, x_2)=0$ whenever $\rho_2(x_1, x_2)=0$.

\begin{lemma}\label{L-non-minimal}
Let $W$ be a nonminimal $G$-subsystem of $R_{\pi}$.
Then there exists a $\rho_0\in \CP(X)$
such that
$((\sigma_{\rho, H}\times \sigma_{\rho, H})(W), H)$ is nonminimal for every $\rho\in \CP(X)$ with
$\rho_0\preceq \rho$ and every countable subgroup $H$ of $G$.
\end{lemma}

\begin{proof}
Let $W'$ be a minimal $G$-subsystem of $W$, and let $(x_1, x_2)\in W\setminus W'$.
As in the proof of \cite[Lemma VI.4.41]{Vri}, it can be shown that there exists
a $\rho_0\in \CP(X)$ such that
$\max(\rho_0(x'_1, x_1), \rho_0(x'_2, x_2))>0$ for all $(x'_1, x'_2)\in W'$.
Then for every $\rho\in \CP(X)$ with
$\rho_0\preceq \rho$ and every countable subgroup $H$ of $G$
we have $(\sigma_{\rho, H}\times \sigma_{\rho, H})(x_1 ,x_2 )\notin
(\sigma_{\rho, H}\times \sigma_{\rho, H})(W')$ so that
$((\sigma_{\rho, H}\times \sigma_{\rho, H})(W), H)$ is nonminimal.
\end{proof}

For all $\rho_1\preceq \rho_2$ in $\CP(X)$ and all countable subgroups $H_1\subseteq H_2$ of $G$
it is clear that there exist unique maps $\sigma_{21}:X^*_{\rho_2, H_2}\rightarrow X^*_{\rho_1, H_1}$ and
$\tau_{21}:Y_{\rho_2, H_2}\rightarrow Y_{\rho_1, H_1}$ such that
$\sigma_{21}\circ \sigma_{\rho_2, H_2}=\sigma_{\rho_1, H_1}$ and
$\tau_{21}\circ \tau_{\rho_2, H_2}=\tau_{\rho_1, H_1}$.
For $\rho_1\preceq \rho_2\preceq \cdots$ in $\CP(X)$ and
countable subgroups $H_1\subseteq H_2\subseteq \cdots $ of $G$, we
write $X^*_{\infty}$ and $Y_{\infty}$ for $\varprojlim X^*_{\rho_n, H_n}$ and
$\varprojlim Y_{\rho_n, H_n}$, respectively. Then
we have induced maps $\sigma_{\infty}:X\rightarrow  X^*_{\infty}$,
$\pi_{\infty}:X^*_{\infty}\rightarrow Y_{\infty}$, and
$\tau_{\infty}:Y\rightarrow Y_{\infty}$. The diagram
\begin{gather}\label{D-limit}
\xymatrix{
X \ar[d]_{\sigma_{\infty}} \ar[r]^{\pi} &  Y  \ar[d]^{\tau_{\infty}} \\
X^*_{\infty}   \ar [r]^{\pi_{\infty}} & Y_{\infty}
}
\end{gather}
is easily seen to commute. Moreover, (\ref{D-limit}) can be identified in a natural way
with (\ref{D-open to quotient}) taking
\begin{gather}\label{E-rho H}
\rho:=\sum^{\infty}_{n=1}2^{-n}\rho_n/(\diam(\rho_n)+1) \hspace*{3mm}\mbox{and}\hspace*{3mm}
H=\bigcup^{\infty}_{n=1}H_n .
\end{gather}

\begin{lemma}\label{L-semi-open}
Suppose that $(X, G)$ is minimal.
Let $W$ be a topologically transitive $G$-subsystem of $R_{\pi}$ such that
the coordinate projections $\pi_i:W\rightarrow X$ are both semi-open. Then
for every $\rho_0\in \CP(X)$ and countable subgroup $H_0$ of $G$ there
exist a $\rho\in \CP(X)$ with $\rho_0\preceq \rho$ and a countable
subgroup $H$ of $G$ with $H_0\subseteq H$ such that
\begin{enumerate}
\item[(a)] $(X^*_{\rho, H}, H)$ is a minimal $H$-system,

\item[(b)] $((\sigma_{\rho, H}\times \sigma_{\rho, H})(W), H)$ is a topologically transitive $H$-system,

\item[(c)] the coordinate projections
$\pi^*_i:(\sigma_{\rho, H}\times \sigma_{\rho, H})(W)\rightarrow X^*_{\rho, H}$ for $i=1,2$ are both semi-open.
\end{enumerate}
\end{lemma}

\begin{proof}
We shall show by induction that there exist 
$\rho_0\preceq \rho_1 \preceq \cdots$ in $\CP(X)$,
countable subgroups $H_0\subseteq H_1\subseteq \cdots $ of $G$ and
a sequence $\{(x_n, x'_n)\}_{n\in \Zb_{\ge 0}}$ of elements in $W$ such that for every $n\in \Zb_{\ge 0}$ the following
conditions are satisfied:
\begin{enumerate}
\item[(a')] $H_{n+1}\sigma^{-1}_{\rho_n, H_n}(U)=X$ for any nonempty open subset $U$ of $X^*_{\sigma_n, H_n}$,

\item[(b')] $(\sigma_{\rho_n, H_n}\times\sigma_{\rho_n, H_n})(H_{n+1}(x_n, x'_n))$ is dense in
$(\sigma_{\rho_n, H_n}\times\sigma_{\rho_n, H_n})(W)$,

\item[(c')] for any nonempty open subset $V$ of $(\sigma_{\rho_n, H_n}\times\sigma_{\rho_n, H_n})(W)$ there exist
$t_1, t_2\in X$ and $\delta>0$ such that $B_{\rho_{n+1}}(t_i, \delta)\subseteq
\pi_i((\sigma_{\rho_n, H_n}\times\sigma_{\rho_n, H_n})^{-1}(V)\cap W)$
for $i=1,2$, where $B_{\rho_{n+1}}(t_i, \delta):=\{x\in X : \rho_{n+1}(x, t_i)<\delta\}$.
\end{enumerate}
We first indicate how this can be used to prove the lemma.
Define $\rho$ and $H$ via (\ref{E-rho H}). Then conditions (a) and (b)
follow from (a') and (b'), respectively, as in the proof of \cite[Lemma VI.4.43]{Vri}.
We may identify $(\sigma_{\rho, H}\times\sigma_{\rho, H})(W)$
with $(\sigma_{\infty}\times \sigma_{\infty})(W)=\varprojlim (\sigma_{\rho_n, H_n}\times\sigma_{\rho_n, H_n})(W)$.
Thus for any nonempty open subset $V'\subseteq (\sigma_{\rho, H}\times\sigma_{\rho, H})(W)$ we
can find an $n\in \Nb$ and a nonempty open subset $V$ of $(\sigma_{\rho_n, H_n}\times\sigma_{\rho_n, H_n})(W)$
such that $V'\supseteq (\sigma_{\infty, n}\times \sigma_{\infty, n})^{-1}(V)
\cap (\sigma_{\infty}\times \sigma_{\infty})(W)
=(\sigma_{\infty}\times \sigma_{\infty})((\sigma_{\rho_n, H_n}\times \sigma_{\rho_n, H_n})^{-1}(V)\cap W)$, where
$\sigma_{\infty,n}:X^*_{\infty}\rightarrow X^*_{\sigma_n, H_n}$ is the natural map.
Let $t_i$ and $\delta$ be as in (c').
Then
\begin{align*}
\pi^*_i(V')&\supseteq \pi^*_i((\sigma_{\infty}\times \sigma_{\infty})
((\sigma_{\rho_n, H_n}\times \sigma_{\rho_n, H_n})^{-1}(V)\cap W)) \\
&= \sigma_{\rho, H}(\pi_i((\sigma_{\rho_n, H_n}\times \sigma_{\rho_n, H_n})^{-1}(V)\cap W))\\
&\supseteq \sigma_{\rho, H}(B_{\rho_{n+1}}(t_i, \delta))\supseteq \sigma_{\rho, H}(B_{\rho}(t_i, \delta'))\\
&\supseteq \tilde{\pi}^{-1}_{\rho, H}(B_{d'_{\rho, H}}(\psi_{\rho, H}(t_i), \delta')),
\end{align*}
where $\delta':=2^{-n-1}\delta /(\diam(\rho_{n+1})+1)$, $d'_{\rho, H}$ is the metric on $X_{\rho, H}$ defined in
(\ref{E-d'_H}),
and $\tilde{\pi}_{\rho, H}: X^*_{\rho, H}\rightarrow X_{\rho, H}$ is the coordinate projection.
Therefore $ \pi^*_i(V')$ has nonempty
interior and hence condition (c) follows from (c').

It remains to construct $\{\rho_n\}_{n\ge 0}$, $\{H_n\}_{n \ge 0}$, and $\{(x_n, x'_n)\}_{n\ge 0}$
satisfying (a'), (b') and (c'). Note that (a') and (b') do not depend on the choice of
$\rho_{n+1}$ while (c') does not depend on the choice of $H_{n+1}$ and $(x_n, x'_n)$.
Thus we can choose $H_{n+1}$ and $(x_n, x'_n)$ to satisfy (a') and (b') exactly as in the
proof of \cite[Lemma VI.4.43]{Vri}. Now we explain how to choose $\rho_1$ in
order to satisfy (c').
Repeating the procedure we will obtain the desired $\rho_n$.

Let $\{V_1, V_2, \dots\}$ be a countable nonempty open base for the topology on the metrizable space
$(\sigma_{\rho_0, H_0}\times\sigma_{\rho_0, H_0})(W)$.
Then $Z_{m, i}:=\pi_i((\sigma_{\rho_0, H_0}\times\sigma_{\rho_0, H_0})^{-1}(V_m)\cap W)$ has nonempty interior
for each $m\in \Nb$ and $i=1,2$. Pick a $t_{m,i}$ in the interior of $Z_{m, i}$, and take
a continuous function $f_{m, i}$ on $X$ such that $0\le f_{m, i}\le 1$, $f_{m, i}(t_{m, i})=0$ and
$f_{m, i}|_{Z^{\comp}_{m, i}}=1$. Define $\rho_{m, i}\in \CP(X)$ by $\rho_{m,i}(z_1, z_2)=|f_{m, i}(z_1)-f_{m, i}(z_2)|$
for $z_1, z_2\in X$, and define
$\rho_1\in \CP(X)$ by $\rho_1=\rho_0+\sum^{\infty}_{m=1}2^{-m}(\rho_{m, 1}+\rho_{m, 2})$.
Then $B_{\rho_1}(t_{m, i}, 2^{-m})\subseteq Z_{m, i}$. One sees immediately that
(c') holds for $n=0$. This completes the proof of Lemma~\ref{L-semi-open}.
\end{proof}

\begin{lemma}\label{L-injective}
Suppose that $Y$ is metrizable.
Then there exists a $\rho_1\in \CP(X)$
such that
$\tau_{\rho, H}$ is an isomorphism for every $\rho\in \CP(X)$ with
$\rho_1\preceq \rho$ and every countable subgroup $H$ of $G$.
\end{lemma}

\begin{proof} Let $d$ be a metric on $Y$ inducing its topology.
Simply define $\rho_1\in \CP(X)$ by $\rho_1(x_1, x_2)=d(\pi(x_1), \pi(x_2))$ for
all $x_1, x_2\in X$.
\end{proof}

\begin{proof}[Proof of Lemma~\ref{L-open tame to HPI}]
Suppose that $\pi$ is not HPI. By the theorem of van der Woude there exists
a topologically transitive nonminimal $G$-subsystem $W$ of $R_{\pi}$ such that
the coordinate projections $W\rightarrow X$ are both semi-open.
Let $\rho_0$ and $\rho_1$ be as in Lemmas~\ref{L-non-minimal} and
\ref{L-injective}, respectively. Take $\rho$ and $H$ as in
Lemma~\ref{L-semi-open} such
that $\rho\succeq \rho_0+\rho_1$. Then $\tau_{\rho, H}$ is an isomorphism
and $((\sigma_{\rho, H}\times\sigma_{\rho, H})(W), H)$ is nonminimal.
Since $\pi$ is tame as an extension of $G$-systems, it is also tame as an extension
of $H$-systems. By Proposition~\ref{P-basic R extension}(1), $\pi_{\sigma, H}$ is a tame
extension of $H$-systems.
By Lemma~\ref{L-tame metrizable to HPI} $\pi_{\rho, H}$ is HPI.
Thus any topologically transitive $H$-subsystem of $R_{\pi_{\rho, H}}$
for which both of the coordinate projections to $X^*_{\rho, H}$ are semi-open must be
minimal by the theorem of van der Woude. This is a contradiction to 
Lemma~\ref{L-semi-open}. Therefore $\pi$ is HPI.
\end{proof}

The proof of
the next lemma follows essentially that
of \cite[Theorem 4.3]{NSSEP}.

\begin{lemma}\label{L-PI to proximal}
Suppose that $G$ is Abelian. Let $\pi:X\rightarrow Y$ be a PI (resp.\ HPI) tame extension of 
minimal systems such that every equicontinuous factor of $(X, G)$ factors through $\pi$.
Then $\pi$ is proximal (resp.\ highly proximal).
\end{lemma}

\begin{proof}
Let $\varphi:X'\rightarrow X$ be a proximal (resp.\ highly proximal)
extension such that $\pi\circ \varphi$ is strictly PI (resp.\ strictly HPI).
Replacing $X'$ by a minimal $G$-subsystem we may assume that $X'$ is minimal. 
Denote by $\psi$ the canonical extension
of $(X, G)$ over its maximal equicontinuous factor $(X_\eq , G)$.
Then $\psi$ factors through $\pi$ by our assumption.
Since $\varphi$ is proximal, $(X_\eq , G)$ is also the maximal equicontinuous factor of $(X', G)$.
Now we need:

\begin{lemma}\label{L-no equi}
Let $(Y'_1, G)$ and $(Y'_2, G)$ be systems
with extensions
$X'\rightarrow Y'_1$, $\theta : Y'_1\rightarrow Y'_2$, and $Y'_2\rightarrow Y$ such that the composition
$X'\rightarrow Y'_1\rightarrow Y'_2\rightarrow Y$ is $\pi\circ \varphi$.
Suppose that $\theta$ is distal. Then $\theta$ is an isomorphism.
\end{lemma}

\begin{proof}
Clearly $(X_\eq , G)$ is also the maximal equicontinuous factor of $(Y'_1, G)$.
Suppose that $y_1'$ and $y'_2$ are distinct points in $Y'_1$ with
$\theta(y'_1)=\theta(y'_2)$.
Since $G$ is Abelian, $Y'_1$ has a $G$-invariant Borel probability measure by
a result of Markov and Kakutani \cite[Theorem VII.2.1]{Dav}.
Then $(Y'_1/\RP(Y'_1, G), G)$ is the maximal equicontinuous factor 
of $(Y'_1, G)$
by \cite[Theorem 9.8]{Aus} \cite[Theorem V.1.17]{Vri}.
Thus $(y'_1, y'_2)\in \RP(Y'_1, G)$.
So $(y'_1, y'_2)\in \IT_2(Y'_1, G)$ by Lemma~\ref{L-distal ext}.
By Proposition~\ref{P-basic R} we can find in $X'$ preimages $x'_1$ and $x'_2$
of $y'_1$ and $y'_2$, respectively,
such that $(x'_1, x'_2)\in \IT_2(X', G)$. Then $(\varphi(x'_1), \varphi(x'_2))\in \IT_2(X, G)\cap R_{\pi}=\IT_{\pi}$.
Since $\theta$ is distal, the pair $(y'_1, y'_2)$ is distal.
Then the pair $(x'_1, x'_2)$ is also distal. 
As $\varphi$ is proximal, $\varphi(x'_1)\neq \varphi(x'_2)$.
Thus $\pi$ is not tame,
which contradicts our assumption. Therefore $\theta$ is an isomorphism.
\end{proof}

Back to the proof of Lemma~\ref{L-PI to proximal}. Since equicontinuous extensions are distal, by
Lemma~\ref{L-no equi} we conclude that $\pi\circ \varphi$ can be obtained by
a transfinite succession of proximal (resp.\ highly proximal) extensions.
Since proximal (resp.\ highly proximal) extensions are preserved under transfinite compositions,
$\pi\circ\varphi$ is proximal (resp.\ highly proximal),
and hence $\pi$ is proximal (resp.\ highly proximal).
\end{proof}

\begin{theorem}\label{T-metric tame minimal to highly proximal}
Suppose that $G$ is Abelian. Let $\pi: X\rightarrow Y$ be a tame extension of
minimal systems.
Consider the following conditions:
\begin{enumerate}
\item $\pi$ is highly proximal,

\item $\pi$ is proximal,

\item every equicontinuous factor of $(X, G)$ factors through $\pi$.
\end{enumerate}
Then one has (1)$\Rightarrow$(2)$\Leftrightarrow$(3).
Moreover, if $X$ is metrizable or $\pi$ is open and $Y$ is metrizable, then
conditions (1) to (3) are all equivalent.
\end{theorem}

\begin{proof}
The implications (1)$\Rightarrow$(2)$\Rightarrow$(3) are trivial, and
(3)$\Rightarrow$(2) follows from Lemmas~\ref{L-tame to PI}
and \ref{L-PI to proximal}. When $X$ is metrizable, (3)$\Rightarrow$(1) follows from
Lemmas~\ref{L-tame metrizable to HPI} and \ref{L-PI to proximal}.
When $\pi$ is open and $Y$ is metrizable, (3)$\Rightarrow$(1) follows from
Lemmas~\ref{L-open tame to HPI} and \ref{L-PI to proximal}.
\end{proof}

\begin{theorem}\label{T-decomposition}
Suppose that $G$ is Abelian. Let $\pi: X\rightarrow Y$ be a tame extension of
minimal systems. Then, up to isomorphisms,
$\pi$ has a unique decomposition as $\pi=\varphi\circ \psi$ such that
$\psi$ is proximal and
$\varphi$ is equicontinuous. If furthermore $X$ is metrizable or $\pi$ is
open and $Y$ is metrizable, then $\psi$ is highly proximal.
\end{theorem}

\begin{proof}
When such a decomposition exists, clearly $\varphi$ must be the maximal equicontinuous factor of
$\pi$. This proves uniqueness.
Let $\RP_{\pi}$ and $S_{\pi}$ be as in the proof of Lemma~\ref{L-B+tame}.
Then we have the natural extensions $\psi:X\rightarrow X/S_{\pi}$ and
$\varphi: X/S_{\pi}\rightarrow Y$ with $\pi=\varphi\circ \psi$, and
$\varphi$ is equicontinuous.
For any extension $\theta : X\rightarrow W$ with $(W, G)$ equicontinuous,
since $\RP_{\pi}\subseteq \RP(X, G)\subseteq R_{\theta}$
we see that $\theta$ factors through $\psi$.
As $\pi$ is tame, so is $\psi$.
By Theorem~\ref{T-metric tame minimal to highly proximal} $\psi$ is proximal.
If $X$ is metrizable, then by Theorem~\ref{T-metric tame minimal to highly proximal} $\psi$ is
highly proximal. If $\pi$ is open and $Y$ is metrizable, then by Lemma~\ref{L-open tame to HPI}
$\pi$ is HPI. Using van der Woude's characterization of HPI extensions one sees immediately
that $\psi$ is also HPI. By Lemma~\ref{L-PI to proximal}, $\psi$ is highly proximal.
\end{proof}

\begin{corollary}\label{C-HP over eq}
Suppose that $G$ is Abelian. Then every tame minimal system $(X, G)$ is a
highly proximal extension of an equicontinuous system.
\end{corollary}

Corollary~\ref{C-HP over eq} answers a question of 
Glasner \cite[Problem 2.5]{tame}, who asked whether every metrizable 
tame minimal system $(X, G)$ with $G$ Abelian is a proximal extension of
an equicontinuous system. 
It also generalizes \cite[Theorem 4.3]{NSSEP} in which the conclusion
is established for metrizable null
minimal systems $(X, \Zb)$.

Recall that a subset $H\subseteq G$ is called {\it thick} if for any
finite $F\subseteq G$, one has $H\supseteq sF$ for some $s\in G$.
We say that $H\subseteq G$ is {\it Poincar{\'e}} if for any measure preserving
action of $G$ on a finite measure space $(Y, \mathscr{B}, \mu)$ and
any $A\in \mathscr{B}$ with $\mu(A)>0$, one has $\mu(sA\cap A)>0$
for some $s\in H$. The argument on page~74 of \cite{RCNT} shows that
every thick set is Poincar{\'e}. 



For a dynamical system $(X, G)$ 
and a Borel subset $U\subseteq X$, denote by
$N(U, U)$ the set $\{s\in G: sU\cap U\neq \emptyset\}$.
If $\mu(U)>0$ for some $G$-invariant Borel probability measure $\mu$ on $X$,
then $N(U, U)$ has nonempty intersection with every Poincar{\'e} set. 
In particular, in this case $N(U, U)$ has nonempty intersection with 
every thick set, or, equivalently, $N(U, U)$ is syndetic 
\cite[page 16]{ETJ}.
Using this fact and Lemma~\ref{L-semiopen to IT}, one
sees that the proof of case 1 in \cite[Theorem 3.1]{NSSEP} leads to:

\begin{lemma}\label{L-E to untame}
Suppose that $G$ is Abelian.
Suppose that a metrizable system $(X, G)$
is nonminimal, has a unique minimal subsystem, and has a 
$G$-invariant Borel probability measure with
full support. Then $(X, G)$ is untame.
\end{lemma}

Using Corollary~\ref{C-HP over eq}
and Lemmas~\ref{L-semiopen to IT} and \ref{L-E to untame}
one also sees that the proof of \cite[Theorem 4.4]{NSSEP} works in our context, so that we
obtain:

\begin{lemma}\label{L-uniquely ergodic}
Suppose that $G$ is Abelian. 
Then any metrizable tame minimal system $(X, G)$ is uniquely ergodic.
\end{lemma}

\begin{theorem}\label{T-uniquely ergodic}
Suppose that $G$ is Abelian. 
Then any tame minimal system $(X, G)$ is uniquely ergodic.
\end{theorem}

\begin{proof} Let $\mu_1$ and $\mu_2$ be two 
$G$-invariant Borel probability measures on $X$.
We use the notation established after Lemma~\ref{L-open tame to HPI}.
Denote by $I$ the set of all pairs $(\rho, H)$ such that $\rho\in \CP(X)$, $H$ is a countable subgroup
of $G$, and $(X^*_{\rho, H}, H)$ is minimal.
Take $(Y, G)$ to be the trivial system in Lemma~\ref{L-open to quotient}.
For every $(\rho, H)\in I$
we have, by Lemma~\ref{L-uniquely ergodic}, 
$\sigma_{\rho, H, *}(\mu_1)= \sigma_{\rho, H, *}(\mu_2)$,
where $\sigma_{\rho, H, *}:\M (X)\rightarrow \M (X^*_{\rho, H})$
is the map between spaces of Borel probability measures induced by $\sigma_{\rho, H}$.
Define a partial order on $I$
by $(\rho_1, H_1)\le (\rho_2, H_2)$ if $\rho_1\le \rho_2$ and $H_1\subseteq H_2$.
As mentioned right after Lemma~\ref{L-non-minimal}, when $(\rho_1, H_1)\le (\rho_2, H_2)$
there exists a unique map 
$\sigma_{21}:X^*_{\rho_2, H_2}\rightarrow X^*_{\rho_1, H_1}$
such that
$\sigma_{21}\circ \sigma_{\rho_2, H_2}=\sigma_{\rho_1, H_1}$.
By \cite[Lemma V.3.9]{Vri}, $I$ is directed.
It is easily checked that $X=\varprojlim_{(\rho, H)\in I} X^*_{\rho, H}$.
Thus $\M (X)=\varprojlim_{(\rho, H)\in I} \ M (X^*_{\rho, H})$.
Therefore $\mu_1=\mu_2$.
\end{proof}

Theorem~\ref{T-uniquely ergodic} generalizes
\cite[Theorem 4.4]{NSSEP} in which the conclusion is established for metrizable
null minimal systems $(X, \Zb)$.


\section{I-independence}\label{S-I-indep}

Here we tie together several properties via the notion of
I-independence, which, as Theorem~\ref{T-equiv-indep-prod} suggests,
can be thought of as an analogue of measure-theoretic
weak mixing for $C^*$-dynamical systems (compare also Theorem~\ref{T-even-WM}).

\begin{definition}\label{D-I-indep}
A $C^*$-dynamical system $(A,G,\alpha )$ is said to be
{\it I-independent} if for every finite-dimensional operator subsystem
$V\subseteq A$ and $\varepsilon > 0$ there is a sequence $\{ s_k \}_{k=1}^\infty$ in $G$ 
such that $(s_1, \dots , s_k )$ is a $(1+\varepsilon )$-independence tuple for $V$ for each
$k\geq 1$.
\end{definition}

Note that I-independence is to be distinguished from $\Inf$-independence,
although the two turn out to be equivalent, as the next
proposition demonstrates.

\begin{proposition}\label{P-C-equiv-indep}
Let $(A,G,\alpha )$ be $C^*$-dynamical system. Let
$\mathfrak{S}$ be a collection of finite-dimensional operator subsystems of
$A$ with the property that for every finite set $\Omega\subseteq A$ and
$\varepsilon > 0$ there is a $V\in\mathfrak{S}$ such that
$\Omega\subseteq_\varepsilon V$. Then
the following are equivalent:
\begin{enumerate}
\item $\alpha$ is I-independent,

\item for every $V\in\mathfrak{S}$
and $\varepsilon > 0$ the set $\Ind (\alpha , V, \varepsilon )$ is infinite,

\item for every $V\in\mathfrak{S}$
and $\varepsilon > 0$ the set $\Ind (\alpha , V, \varepsilon )$ is nonempty,

\item $\alpha$ is $\Inf$-independent,

\item $\alpha$ is $\NE$-independent.
\end{enumerate}
\end{proposition}

\begin{proof}
(1)$\Rightarrow$(2)$\Rightarrow$(3) and (4)$\Rightarrow$(5). Trivial.

(2)$\Rightarrow$(4) and (3)$\Rightarrow$(5).
Apply Proposition~\ref{P-dense}.

(5)$\Rightarrow$(1). Let $V$ be a finite-dimensional operator
subsystem of $A$ and let $\varepsilon > 0$. With the aim of verifying
I-independence, we may assume that $V$ does not equal the scalars, so that
$V$ has linear dimension at least two.
By recursion we will construct a sequence $\{ s_1 , s_2 ,\dots \}$ of distinct
elements of $G$ such that for each $k\geq 1$ the linear map
$\varphi_k : V^{\otimes [1,k]} \to V_k :=
[ \alpha_{s_1} (V) \cdots \alpha_{s_k} (V) ] \subseteq A$
determined on elementary tensors by
$a_1 \otimes\cdots\otimes a_k \mapsto \alpha_{s_1} (a_1 )
\cdots \alpha_{s_k} (a_k )$ is a
$(1+\varepsilon )^{1-2^{-k+1}}$-c.b.-isomorphism.

To begin with set $s_1 = e$. Now let $k\geq 1$ and suppose that
$s_1 , s_2 , \dots , s_k$ have been defined so that $\varphi_k$ is a
$(1+\varepsilon )^{1-2^{-k+1}}$-c.b.-isomorphism onto its image. By (3)
there is an $s_{k+1}\in G$ such that the linear map
$\psi : V_k \otimes V_k \to [V_k \alpha_{s_{k+1}} (V_k )]$
determined by $a\otimes b \mapsto a\alpha_{s_{k+1}} (b)$
is a $(1+\varepsilon )^{2^{-k}}$-c.b.-isomorphism. Since $V$ has linear
dimension at least two, we must have $s_{k+1} \notin\{ s_1 , \dots , s_k \}$,
for otherwise $\psi$ would not be injective.
Set $\gamma = \varphi_k \otimes \id_V : V^{\otimes [1,k]} \otimes V =
V^{\otimes [1,k+1]} \to V_k \otimes V$.
By the injectivity of the minimal tensor product, we may
view $V_k \otimes V$ as a subspace of $V_k \otimes V_k$,
in which case we have $\varphi_{k+1} = \psi\circ\gamma$. Then
\begin{align*}
\| \varphi_{k+1} \|_\cb \| \varphi_{k+1}^{-1} \|_\cb &\leq
\| \psi \|_\cb \| \varphi_k \otimes \id_V \|_\cb \| \psi^{-1} \|_\cb
\| \varphi_k^{-1} \otimes \id_V \|_\cb \\
&= \| \psi \|_\cb \| \psi^{-1} \|_\cb \| \varphi_k \|_\cb
\| \varphi_k^{-1} \|_\cb \\
&\leq (1+\varepsilon )^{2^{-k}} (1+\varepsilon )^{1-2^{-k+1}} \\
&= (1+\varepsilon )^{1-2^{-k}} ,
\end{align*}
so that $\varphi_{k+1}$ is a
$(1+\varepsilon )^{1-2^{-k}}$-c.b.-isomorphism, as desired.
Since for each $k\geq 1$ the map $\varphi_k$ is a 
$(1+\varepsilon )$-c.b.-isomorphism, we obtain (1).
\end{proof}

\begin{proposition}\label{P-I-indep-powers}
A $C^*$-dynamical system $(A,G,\alpha )$ is I-independent if and only if
the product system $(A^{\otimes [1,m]},G,\alpha^{\otimes [1,m]} )$ is
I-independent for every $m\in\Nb$.
\end{proposition}

\begin{proof}
For the nontrivial direction we can apply
Proposition~\ref{P-C-equiv-indep} using the following two observations:
(i) the collection of operator subsystems of $A^{\otimes [1,m]}$ of
the form $\bigotimes_{i=1}^m V_i$ for finite-dimensional operator
subsystems $V_1 , \dots , V_m \subseteq A$ satisfies the property
required of $\mathfrak{S}$ in the statement of
Proposition~\ref{P-C-equiv-indep},
and (ii) given finite-dimensional operator subsystems 
$V_1 , \dots , V_m \subseteq A$ and a $\lambda\geq 1$, if a tuple in $G$ is
a $\lambda$-independence tuple for the linear span of $\bigcup_{i=1}^m V_i$ then it 
is a $\lambda^m$-independence tuple for $V_1 \otimes\cdots\otimes V_m$,
as follows from the fact that for any c.b.\ isomorphisms $\varphi_i : E_i \to F_i$
between operator spaces for $i=1,\dots ,m$ the tensor product
$\varphi = \bigotimes_{i=1}^m \varphi_i$ is a c.b.\ isomorphism
with $\| \varphi \|_\cb \| \varphi^{-1} \|_\cb =
\prod_{i=1}^m \| \varphi_i \|_\cb \| \varphi_i^{-1} \|_\cb$.
\end{proof}

\begin{proposition}\label{P-equiv-indep}
Let $(X,G)$ be a dynamical system. Let
$\cB$ be a basis for the topology on $X$ which does not contain
the empty set. Then the following are equivalent:
\begin{enumerate}
\item $(X,G)$ is I-independent,


\item for every finite-dimensional operator subsystem $V\subseteq C(X)$ 
there are a sequence $\{ s_k \}_{k=1}^\infty$ in $G$ and a $\lambda\geq 1$ 
such that $(s_1, \dots , s_k )$ is a $\lambda$-independence tuple for $V$ 
for each $k\geq 1$,

\item every finite-dimensional operator subsystem of $C(X)$ has arbitrarily
long $\lambda$-inde-\linebreak pendence tuples for some $\lambda\geq 1$,

\item every nonempty finite subcollection of $\cB$ has
an infinite independence set,

\item every nonempty finite subcollection of $\cB$ has
arbitrarily large finite independence sets,

\item for every nonempty finite collection $\{ U_1 , \dots , U_m \}
\subseteq\cB$ there is an $s\in G$ such that
$U_i \cap s^{-1} U_j \neq\emptyset$ for all $i,j=1,\dots ,m$.
\end{enumerate}
\end{proposition}

\begin{proof}
(1)$\Rightarrow$(2)$\Rightarrow$(3). Trivial.

(2)$\Rightarrow$(4). Let $\{ U_1 , \dots , U_m \}$ be a nonempty finite
subcollection of $\cB$. By shrinking the $U_i$ if
necessary, we may assume for the purpose of establishing (4) that
$U_i \cap U_j = \emptyset$ for $i\neq j$. For each $i=1,\dots ,m$ take
a nonempty closed set $W_i \subseteq U_i$. Let
$\Theta = \{ g_0 , g_1 , \dots , g_m \}$ be a partition of unity of $X$ such
that $\supp (g_i )\subseteq U_i$ for each $i$ and $\supp (g_0 ) \subseteq
X\setminus \bigcup_{i=1}^m W_i$. Let $V$ be the linear span
of $\Theta$. By (2) there are a sequence $\{ s_1 , s_2 , \dots \}$ in $G$ and a
$\lambda\geq 1$ such that for each $k\geq 1$
the contractive linear map $V^{\otimes [1,k]} \to C(X)$
determined on elementary tensors by
$f_1 \otimes\cdots\otimes f_k \mapsto \alpha_{s_1} (f_1 )
\cdots \alpha_{s_k} (f_k )$ is a $\lambda$-isomorphism
onto its image. Now suppose we are given a $k\in\Nb$ and a 
$\sigma\in {\{ 1,\dots ,m \}}^{\{ 1 , \dots , k \}}$.
Then $\| \alpha_{s_1} (g_{\sigma (1)} )\cdots\alpha_{s_k}
(g_{\sigma (k)} ) \| \geq \lambda^{-1}$. Choose $x\in X$ such that
$g_{\sigma (1)} (s_1 x)\cdots g_{\sigma (k)} (s_k x) =
\| \alpha_{s_1} (g_{\sigma (1)} )\cdots\alpha_{s_k}
(g_{\sigma (k)} ) \|$. Then for each
$i=1,\dots ,k$ we must have $g_{\sigma (i)} (s_i x) > 0$,
which implies that $x\in s_i^{-1} U_{\sigma (i)}$. Hence
$\bigcap_{i=1}^k s_i^{-1} U_{\sigma (i)} \neq\emptyset$, and so we
obtain (4).

(4)$\Rightarrow$(5)$\Rightarrow$(6). Trivial.

(3)$\Rightarrow$(5). This can be proved along the lines of the
argument used for (2)$\Rightarrow$(4).

(6)$\Rightarrow$(1). Let $\Theta = \{ g_1 , \dots , g_m \}$ be a
partition of unity of $X$ for which there are elements
$U_1 , \dots , U_m$ of $\cB$ such that $g_i |_{U_j} =
\delta_{ij} \chi_{U_j}$, where $\chi_{U_j}$ is the characteristic
function of $U_j$. Note that the collection of subspaces
of $C(X)$ spanned by such partitions of unity has the property
required of $\mathfrak{S}$ in the statement of
Proposition~\ref{P-C-equiv-indep}.
By (6) there is an $s\in G$ such that $U_i \cap s^{-1} U_j
\neq\emptyset$ for all $i,j=1,\dots ,m$.
Let $V$ be the subspace of $C(X)$ spanned by
$\Theta$ and let $\varphi : V\otimes V \to A$ be the linear
map determined on elementary tensors by
$f\otimes g \mapsto f\alpha_s (g)$. Then
$\{ g_i \alpha_s (g_j ) \}_{1\leq i,j\leq m}$ is an effective partition of
unity of $X$, and hence is isometrically equivalent to the standard basis of
$\ell_\infty^{m^2}$. Since the subset
$\{ g_i \otimes g_j \}_{1\leq i,j\leq m}$ of
$V\otimes V$ is also isometrically equivalent to standard basis
of $\ell_\infty^{m^2}$, we conclude that $\varphi$ is an isometric
isomorphism. In view of Proposition~\ref{P-C-equiv-indep},
this yields (1).
\end{proof}

\begin{remark}
In the case $G=\Zb$, the proof of Proposition~\ref{P-C-equiv-indep}
shows that the sequence in the
definition of I-independence can be taken to be in $\Nb$.
since $\Ind (\alpha , V, \varepsilon )$ is closed in general
under taking inverses (given a c.b.\ isomorphism 
$\varphi : V\otimes V \to [V\alpha_s (V)]$
with $\varphi (v\otimes w) = v\alpha_s (w)$, the map
$V\otimes V \mapsto [V\alpha_{s^{-1}} (V)]$ defined by
$v\otimes w \mapsto v\alpha_{s^{-1}} (w) = 
\alpha_{s^{-1}} (\varphi (w^* \otimes v^* )^* )$
has the same c.b.\ isomorphism constant). 
Thus in Proposition~\ref{P-equiv-indep} the sequence in each of
conditions (2) and (4) can be taken to be in $\Nb$.
\end{remark}

The next theorem extends \cite[Theorem 2.1]{NSSEP}.
Recall that $(X, G)$ is said to be {\it (topologically) transitive}
if every nonempty open invariant subset of $X$ is dense,
and {\it (topologically) weakly mixing} if the product system
$(X\times X, G)$ is transitive.
For the definitions of
uniform nonnullness and untameness of all orders see Sections~\ref{S-null}
and \ref{S-tame}, respectively.

\begin{theorem}\label{T-equiv-indep-prod}
Let $(X,G)$ be a dynamical system and consider the following conditions:
\begin{enumerate}
\item $(X,G)$ is I-independent,

\item $(X,G)$ is uniformly untame of all orders,

\item $(X,G)$ is uniformly nonnull of all orders,

\item for every $n\in\Nb$ the product system $(X^n ,G)$ is weakly mixing,

\item for every $n\in\Nb$ the product system $(X^n ,G)$ is transitive,

\item $(X,G)$ is uniformly untame,

\item $(X,G)$ is uniformly nonnull,

\item $(X,G)$ is weakly mixing.
\end{enumerate}
Then conditions (1) to (5) are equivalent, and
(5)$\Rightarrow$(6)$\Rightarrow$(7).
When $G$ is Abelian, conditions (1) to (8) are all equivalent.
\end{theorem}

\begin{proof}
The implications
(2)$\Rightarrow$(6)$\Rightarrow$(7) are trivial.
In the case that $G$ is
Abelian, (7)$\Rightarrow$(8) follows as in the proof of
\cite[Theorem 2.1]{NSSEP}
using the lemma in \cite{WM}, while (8)$\Rightarrow$(5) is 
\cite[Proposition II.3]{Fur}.
The implication (1)$\Rightarrow$(2) follows from
Proposition~\ref{P-equiv-indep},
(2)$\Rightarrow$(3) and (4)$\Rightarrow$(5) are trivial, and
(3)$\Rightarrow$(4) is readily seen.

Finally, to show (5)$\Rightarrow$(1), let $\{ U_1 , \dots , U_n \}$ be a nonempty
finite collection of nonempty open subsets of $X$. Let $\Lambda$ be the
set of all integer pairs $(i,j)$ with $1\leq i,j\leq n$.
Set $V_0 = \prod_{(i,j)\in\Lambda} U_i \subseteq X^\Lambda$
and $V_1 = \prod_{(i,j)\in\Lambda} U_j \subseteq X^\Lambda$.
By (5) there is an $s\in G$ such that
$V_0 \cap s^{-1} V_1 \neq\emptyset$, so that $U_i \cap s^{-1} U_j \neq\emptyset$ for
all $i,j=1, \dots ,n$. It then follows by Proposition~\ref{P-equiv-indep}
that $(X,G)$ is I-independent.
\end{proof}

By the above theorem, we can consider I-independence to be the analogue
for noncommutative $C^*$-dynamical systems of, among other properties,
uniform untameness of all orders
(compare the discussion in the last paragraph of Section~\ref{S-entropy tensor}).
On the other hand, the notions
of tameness (resp.\ untameness) and complete untameness make sense for any $C^*$-dynamical
system, the former meaning that no element (resp.\ some element) has an infinite $\ell_1$-isomorphism
set and the latter that every nonscalar element has an infinite
$\ell_1$-isomorphism set. We will end this section by observing
that, in the general $C^*$-dynamical context, I-independence (in fact
a weaker independence condition) implies complete untameness.

%

To verify complete untameness it suffices to check the existence
of an infinite $\ell_1$-isomorphism set over $\Rb$ for every
self-adjoint nonscalar element, as the following lemma
based on Rosenthal-Dor arguments demonstrates.

\begin{lemma}\label{L-sa}
Let $(A,G,\alpha )$ be a $C^*$-dynamical system.
Let $a\in A$, and suppose that at least one of $\re (a)$ and $\im (a)$
has an infinite $\ell_1$-isomorphism set over $\Rb$. Then $a$ has an infinite
$\ell_1$-isomorphism set.
\end{lemma}

\begin{proof}
Suppose that $\{ \alpha_{s_n} (a) \}_{n\in\Nb}$ is a sequence in the orbit of
$a$ which converges
weakly to some $b\in A$. Then $\lim_{n\to\infty} \sigma (\re (\alpha_{s_n} (a)))
= \sigma (\re (b))$ and $\lim_{n\to\infty} \sigma (\im (\alpha_{s_n} (a)) =
\sigma (\im (b))$ for all self-adjoint $\sigma\in A^*$, so that both
$\{ \alpha_{s_n} (\re (a)) \}_{n\in\Nb}$ and
$\{ \alpha_{s_n} (\im (a)) \}_{n\in\Nb}$
are weakly convergent over $\Rb$. It follows from this observation and
\cite{Dor} that $a$ has an infinite $\ell_1$-isomorphism set.
\end{proof}

For a unital $C^*$-algebra $A$, we denote by $\cS_2 (A)$ the collection of 
$2$-dimensional operator subsystems of $A$ equipped with the metric given
by Hausdorff distance between unit balls.

\begin{proposition}\label{P-I-indep-cnr}
Let $(A,G,\alpha )$ be a $C^*$-dynamical system. Let $\lambda\geq 1$,
and suppose that in $\cS_2 (A)$ there is a dense
collection of $V$ for which
there is a sequence $\{ s_1 , s_2 , \dots \}$ in $G$ such that $(s_1 , \dots , s_k )$
is a $\lambda$-independence tuple for $V$ for each $k\geq 1$. Then $\alpha$ is completely 
untame. In particular, an I-independent $C^*$-dynamical system is completely 
untame.
\end{proposition}

\begin{proof}
Let $a$ be a nonscalar self-adjoint element of $A$, and denote by $V$ the 
$2$-dimensional operator system $\spn \{ 1,a \}$. Suppose that 
there is a sequence $\{ s_1 , s_2 , \dots \}$ in $G$ such that for each $k$
the linear map $V^{\otimes [1,k]} \to A$ determined by
$a_1 \otimes\cdots\otimes a_k \mapsto \alpha_{s_1} (a_1 )
\cdots \alpha_{s_k} (a_k )$ is a
$\lambda$-c.b.-isomorphism onto its image. It then follows by
Lemma~\ref{L-tensor-indep} that the set 
$\{ \alpha_{s_k} (a) \}_{k\in\Nb}$
is $\lambda'$-equivalent to the standard basis of
$\ell_1^\Nb$ for some $\lambda'$ depending only on $\lambda$ and the 
spectral diameter of $a$. By a straightforward perturbation argument and
Lemma~\ref{L-sa} we conclude that $\alpha$ is completely untame.
\end{proof}


\section{Independence, Abelianness, and weak mixing in
$C^*$-dynamical systems}\label{S-NC}

Once we express independence in terms of minimal tensor products and
move into the noncommutative realm, a close connection to Abelianness
reveals itself. Indeed one of the
goals of this section is to show that, in simple unital nuclear $C^*$-algebras,
independence and Abelianness in the dynamical context essentially
amount to the same thing (Theorem~\ref{T-nuclear}). This
provides a conceptual basis for the sense gained from examples that concepts
like hyperbolicity and topological K-ness should be
interpreted in the simple nuclear case as certain types of asymptotic Abelianness.
See for instance \cite{Narn,Nesh} and the discussion at the end of
Section~\ref{S-entropy tensor}, and compare also \cite{Dis}, where a relationship 
between tensor product structure and asymptotic Abelianness is established as
a tool in the study of derivations and dissipations. The second
main result of this section concerns the implications for independence
of the existence of weakly mixing invariant states (Theorem~\ref{T-WMcontrInd}).

\begin{definition}\label{D-Abelian}
Let $(A,G,\alpha )$ be a $C^*$-dynamical system. Let $V$ be a
finite-dimensional operator subsystem of $A$ and let $\varepsilon > 0$.
We define $\Ab (\alpha , V, \varepsilon )$ to be the set of all
$s\in G_0$ such that $\| [v,\alpha_s (w) ] \| \leq\varepsilon
\| v \| \| w \|$ for all $v,w\in V$. For a collection $\cC$ of subsets of $G_0$
which is closed under taking supersets,
we say that the system $(A,G,\alpha )$ or the action $\alpha$ is
{\it $\cC$-Abelian} if for every finite-dimensional operator subsystem
$V\subseteq A$ and $\varepsilon > 0$
the set $\Ab (\alpha , V, \varepsilon )$ is a member of $\cC$.
\end{definition}

When $G$ is an infinite group, $\Inf$-Abelianness can be characterized 
as follows.

\begin{proposition}\label{P-Abeliandense}
Let $(A,G,\alpha )$ be a $C^*$-dynamical system with $G$ an infinite group.
Then the following are equivalent:
\begin{enumerate}
\item $\alpha$ is $\Inf$-Abelian,

\item $\alpha$ is $\NE$-Abelian,

\item there is a net $\{ s_\gamma \}_\gamma$ in $G$ (which can be taken to be
a sequence if $A$ is separable) such that
$\lim_\gamma \| [a,\alpha_{s_\gamma} (b) ] \| = 0$ for all $a,b\in A$.
\end{enumerate}
\end{proposition}

\begin{proof}
The implication (1)$\Rightarrow$(2) is trivial. Given a net $\{ s_\gamma \}_\gamma$
as in (3), if it has a subnet tending to infinity then we obviously obtain (1), and if 
not then it has convergent subnet and hence $A$ is commutative, so that
we again obtain (1). Suppose then that (2) holds and let us prove (3). Let
$\mathfrak{D}$ be a set of finite-dimensional operator subsystems
of $A$ directed by inclusion such that $\bigcup\mathfrak{D}$ is dense in $A$.
Let $\Gamma$ be the directed set of all pairs $(V,\varepsilon )$ such that
$V\in\mathfrak{D}$ and $\varepsilon > 0$, where $(V' ,\varepsilon' )\succ
(V,\varepsilon )$ means that $V' \supseteq V$ and $\varepsilon' \leq\varepsilon$.
For every $\gamma = (V,\varepsilon )\in\Gamma$ we can find by (2) an
$s_\gamma \in G$ such that $\| [v,\alpha_{s_\gamma} (w) ] \|
\leq\varepsilon \| v \| \| w \|$ for all $v,w\in V$. Then
$\lim_\gamma \| [a,\alpha_{s_\gamma} (b) ] \| = 0$ for
all $a,b\in A$, as desired. In the case that $A$ is separable we can use instead
a sequence of pairs $(V_n ,1/n)$ where $V_1 \subseteq V_2\subseteq \dots$ is
an increasing sequence of finite-dimensional operator subsystems of $A$ such that
$\bigcup_{n\in\Nb} V_n$ is dense in $A$.
\end{proof}

Condition (3) in Proposition~\ref{P-Abeliandense} in the case where $G$ is the
entire $^*$-automorphism group of $A$ is of importance in $C^*$-algebra
structure and classification theory. See for example Lemma~5.2.3 in \cite{Ror}.

\begin{lemma}\label{L-cp}
For every $\varepsilon > 0$ there is a $\delta > 0$ such that whenever
$V$ is an operator system, $\cH$ is a Hilbert space, and
$\varphi : V\to\cB (\cH )$ is a unital c.b.\ map with
$\| \varphi \|_\cb \leq 1 + \delta$, there is a complete positive map
$\psi : V\to\cB (\cH )$ with $\| \psi - \varphi \| \leq \varepsilon$.
\end{lemma}

\begin{proof}
Let $V$ be an operator system, $\cH$ a Hilbert space, and
$\varphi : V\to\cB (\cH )$ a unital c.b.\ map. We may assume $V$ to be an
operator subsystem of a unital $C^*$-algebra $A$. Then by the Arveson-Wittstock
extension theorem we can extend $\varphi$ to a c.b.\ map
$\tilde{\varphi} : A\to\cB (\cH )$ with
$\| \tilde{\varphi} \|_\cb = \| \varphi \|_\cb$. Using the representation
theorem for c.b.\ maps as in the proof of Lemma~3.3 in \cite{PS}, we
then can produce a completely positive map $\rho : A\to\cB (\cH )$
such that $\| \tilde{\varphi} - \rho \|_\cb
\leq \sqrt{2 \| \tilde{\varphi} \|_\cb (\| \tilde{\varphi} \|_\cb - 1)}$,
from which the lemma follows.
\end{proof}

\begin{lemma}\label{L-contrAb}
For every $\varepsilon > 0$ there is a $\delta > 0$ such
that whenever $A$ is a unital $C^*$-algebra, $\alpha$ is a unital
$^*$-endomorphism of $A$, and
$V$ is a finite-dimensional operator subsystem of $A$ for which the linear map
$V\otimes V\to [V\alpha (V)]$ determined on elementary tensors by
$v_1 \otimes v_2 \mapsto v_1 \alpha (v_2 )$ has c.b.\ norm
at most $1+\delta$, we have
$\| [v_1 , \alpha (v_2 ) ] \| \leq \varepsilon \| v_1 \| \| v_2 \|$
for all $v_1 , v_2 \in V$.
\end{lemma}

\begin{proof}
Let $\varepsilon > 0$, and take $\delta$ as given by Lemma~\ref{L-cp}
with respect to $\varepsilon /2$.
Let $A$ be a unital $C^*$-algebra, $\alpha$ a unital $^*$-endomorphism of $A$,
and $V$ an operator subsystem of $A$, and suppose that the linear
map $\varphi : V\otimes V\to [V\alpha (V)]$ determined on elementary
tensors by $v_1 \otimes v_2 \mapsto v_1 \alpha (v_2 )$ has c.b.\ norm
at most $1+\delta$.
We regard $A$ as a unital $C^*$-algebra of $\cB (\cH )$ for some
Hilbert space $\cH$. Then by our choice of $\delta$
there is a completely positive map
$\psi : V\otimes V \to\cB (\cH )$ with $\| \psi - \varphi \| \leq \varepsilon /2$.
It follows that, for all $v_1 , v_2 \in V$,
\begin{align*}
\| v_1 \alpha (v_2 ) - \alpha (v_2 )v_1 \| &\leq
\| \varphi (v_1 \otimes v_2 ) - \psi (v_1 \otimes v_2 ) \|
+ \| \psi (v_1^* \otimes v_2^* )^* - \varphi (v_1^* \otimes v_2^* )^* \| \\
&\leq 2 \| \varphi - \psi \| \| v_1 \| \| v_2 \| \\
&\leq \varepsilon \| v_1 \| \| v_2 \| ,
\end{align*}
yielding the lemma.
\end{proof}

Lemma~\ref{L-contrAb} immediately yields:

\begin{proposition}\label{P-contrAb}
Let $\cC$ a collection of subsets of $G_0$ which is closed under
taking supersets. Then for a $C^*$-dynamical system with acting semigroup $G$,
$\cC$-contractivity implies $\cC$-Abelianness.
\end{proposition}

\begin{theorem}\label{T-nuclear}
Let $(A,G,\alpha )$ be a $C^*$-dynamical system with $A$ 
nuclear. Let $\cC$ be a collection of subsets of $G_0$ which is
closed under taking supersets. Consider the following conditions:
\begin{enumerate}
\item $\alpha$ is $\cC$-independent,

\item $\alpha$ is $\cC$-contractive,

\item $\alpha$ is $\cC$-Abelian.
\end{enumerate}
Then (1)$\Rightarrow$(2)$\Leftrightarrow$(3), and if $A$ is simple then all three
conditions are equivalent.
\end{theorem}

\begin{proof}
The implication (1)$\Rightarrow$(2) is trivial, while (2)$\Rightarrow$(3) is
a special case of Proposition~\ref{P-contrAb}.

Assume that (3) holds and let us prove (2). Suppose that $\alpha$ is not
$\cC$-contractive. Then there are a finite-dimensional operator subsystem
$V\subseteq A$ and an $\varepsilon > 0$ such that no
$(1+\varepsilon )$-contraction set for $V$ is a member of $\cC$. Denote by
$\Lambda$ the set of pairs $(E,\delta )$ where $E$ is a finite-dimensional
operator subsystem of $A$ containing $V$ and $\delta > 0$. For every 
$\lambda = (E,\delta )\in\Lambda$ let $H_\lambda$ be the set of all 
$s\in G$ such that $\| [v,\alpha_s (w)] \| \leq \delta \| v \| \| w \|$ for all
$v,w\in E$ and the linear map
$\varphi_s : V\otimes V \to [V\alpha_s (V)]$
defined on elementary tensors by $\varphi_s (v\otimes w) = v\alpha_s (w)$
has c.b.\ norm greater than $1+\varepsilon$. Since $\alpha$ is $\cC$-Abelian,
$H_\lambda$ is nonempty for every $\lambda\in\Lambda$, and for all 
$\lambda_1 = (E_1 , \delta_1 )$ and $\lambda_2 = (E_2 , \delta_2 )$ in $\Lambda$
we have $H_\lambda \subseteq H_{\lambda_1} \cap H_{\lambda_2}$
for $\lambda = (E_1 + E_2 , \min (\delta_1 , \delta_2 ))$. Thus the collection  
$\{ H_\lambda : \lambda\in\Lambda \}$ forms a filter base over $G$. Let $\omega$ 
be an ultrafilter containing this collection. We write
$\ell_\infty^G (A) / c_\omega (A)$ for the ultrapower of $A$ with respect to $\omega$
(see \cite[Sect.\ 6.2]{Ror}).

Denote by $\pi$ the quotient map $\ell_\infty^G (A) \to
\ell_\infty^G (A) / c_\omega (A)$.
Let $\Phi_1 , \Phi_2 : A \to \ell_\infty^G (A) / c_\omega (A)$ be the
$^*$-homomorphisms given by $\Phi_1 (a) = \pi ((a)_{s\in G} )$ and
$\Phi_2 (a) = \pi ((\alpha_s (a))_{s\in G} )$ for all $a\in A$.
Then the images of $\Phi_1$ and $\Phi_2$ commute, and thus, since
$A \otimes A = A \otimes_\maxtensor A$ by the nuclearity
of $A$, we obtain
a $^*$-homomorphism $\Phi : A\otimes A \to\ell_\infty^G (A) / c_\omega (A)$
such that $\Phi (a_1 \otimes a_2 ) = \Phi_1 (a_1 ) \Phi_2 (a_2 )$
for all $a_1 , a_2 \in A$ \cite[Prop.\ IV.4.7]{Tak1}.

Using the description of nuclearity as a completely positive approximation
property we can recursively construct a separable nuclear 
operator subsystem $W$ of $A$ containing $V$. 
By the Choi-Effros lifting theorem \cite{CE} there is a unital
completely positive map
$\psi : W\otimes W\to\ell_\infty^G (A)$ such that $\pi\circ\psi = \Phi |_{W\otimes W}$,
viewing $W\otimes W$ as a subspace of $A\otimes A$.
Since the unit ball of $V\otimes V$ is compact, $\omega$ contains the
set $H$ of all $s\in G$
for which $\| \varphi_s - \pi_s \circ\psi |_{V\otimes V} \| \leq
\dim (V)^{-2} \varepsilon$, where $\pi_s$ denotes the coordinate projection
$\ell_\infty^G (A) \to A$ associated to $s$ and $V\otimes V$ is viewed as a 
subspace of $W\otimes W$. Appealing to \cite[Cor.\ 2.2.4]{ER}, for $s\in H$ we have
\begin{align*}
\| \varphi_s \|_\cb &\leq \| \pi_s \circ\psi |_{V\otimes V} \|_\cb +
\| \varphi_s - \pi_s \circ\psi |_{V\otimes V} \|_\cb \\
&\leq 1 + \dim (V)^2 \| \varphi_s - \pi_s \circ\psi |_{V\otimes V} \| \\
&\leq 1+\varepsilon ,
\end{align*}
which yields a contradiction since $H$ intersects every $H_\lambda$. 
We thus obtain (2).

Now suppose that $A$ is simple and let us show (3)$\Rightarrow$(1). Suppose that
$\alpha$ is not $\cC$-independent. Then there are a finite-dimensional operator 
subsystem $V\subseteq A$ and an $\varepsilon > 0$ such that no
$(1+\varepsilon )^2$-independence set for $V$ is a member of $\cC$. As before define
$\Lambda$ to be the set of pairs $(E,\delta )$ where $E$ is a finite-dimensional
operator subsystem of $A$ containing $V$ and $\delta > 0$. For every 
$\lambda = (E,\delta )\in\Lambda$ let $H_\lambda$ be the set of all 
$s\in G$ such that 
$\| [v,\alpha_s (w)] \| \leq \delta \| v \| \| w \|$ for all
$v,w\in E$ and the linear map
$\varphi_s : V\otimes V \to [V\alpha_s (V)]$
defined on elementary tensors by $\varphi_s (v\otimes w) = v\alpha_s (w)$
either has c.b.\ norm greater than $1+\varepsilon$ or is not invertible
or has an inverse with c.b.\ norm greater than $1+\varepsilon$. As in the previous
paragraph we construct a separable nuclear operator subsystem $W$ of $A$ 
containing $V$ and apply the Choi-Effros lifting theorem to obtain 
a unital completely positive map
$\psi : W\otimes W\to\ell_\infty^G (A)$ such that $\pi\circ\psi = \Phi |_{W\otimes W}$.
Since $A$ is simple so is $A\otimes A$, and hence $\Phi$ is faithful.
Consequently $\psi$ is a complete order isomorphism.
Let $\delta$ be a positive number to be specified below. 
Since $A$ is nuclear,
$V$ is $1$-exact, and thus $V\otimes V$ is $1$-exact (see \cite{ER,Pis}).
By a result of Smith \cite[Prop.\ 2.2.2]{ER} and the characterization of
$1$-exactness in terms of almost completely isometric embeddings into
matrix algebras (see \cite[Lemma 17.8]{Pis}), it follows that
there exists a $k\in\Nb$
such that whenever $E$ is an operator space and $\rho : E\to V\otimes V$ is
a bounded linear map we have
$\| \rho \|_\cb \leq \sqrt{1+\delta} \| \id_{M_k} \otimes \rho \|$.
Under the canonical
identification $(M_k \otimes \ell_\infty^G (A)) / (M_k \otimes c_\omega (A))
\cong M_k \otimes (\ell_\infty^G (A) / c_\omega (A))$, we regard
the complete order embedding $\id_{M_k} \otimes \psi :
M_k \otimes W\otimes W \to M_k \otimes \ell_\infty^G (A)$ as a lift of the
restriction to $M_k \otimes W\otimes W$ of the $^*$-homomorphism 
$\id_{M_k} \otimes\Phi : M_k \otimes A \otimes A\to
M_k \otimes (\ell_\infty^G (A) / c_\omega (A))$ with respect to the
quotient map $M_k \otimes \ell_\infty^G (A) \to
(M_k \otimes \ell_\infty^G (A)) / (M_k \otimes c_\omega (A))$.
Since the unit
ball of $M_k \otimes V\otimes V$ is compact, $\omega$ contains the set $H$ of 
all $s\in G$ for which (i) 
$\| \varphi_s - \pi_s \circ\psi |_{V\otimes V} \| \leq
\dim (V)^{-2} \varepsilon$ and (ii) $\pi_s \circ\psi |_{V\otimes V} $ is
invertible with
$\| \id_{M_k} \otimes (\pi_s \circ\psi |_{V\otimes V} )^{-1} \| <
\sqrt{1 + \delta}$
and $\| \varphi_s - \pi_s \circ\psi |_{V\otimes V} \| < \delta$.
From (i) we obtain $\| \varphi_s \|_\cb \leq 1 + \varepsilon$ as before,
while from (ii) we obtain 
$\| (\pi_s \circ\psi |_{V\otimes V} )^{-1} \|_\cb < 1 + \delta$ so that,
by Lemma~\ref{L-perturb}, if
$\delta$ is small enough as a function of $\dim (V)$ and $\varepsilon$ then
$\varphi_s$ is invertible with $\| \varphi_s^{-1} \|_\cb \leq 1+\varepsilon$.
This produces a contradiction since $H$ intersects every $H_\lambda$. 
Thus (3)$\Rightarrow$(1) in the simple case.
\end{proof}

For the remainder of this section $G$ will be a group.
Let $(A,G,\alpha )$ be a $C^*$-dynamical system and $\sigma$ a
$G$-invariant state on $A$. We write $\cH_\sigma$ for the GNS Hilbert space of
$\sigma$ and $\Omega_\sigma$ for the vector in $\cH_\sigma$ associated to the
unit of $A$. We denote by $\cH_{\sigma ,0}$ the orthogonal complement in
$\cH_\sigma$ of the one-dimensional subspace spanned by $\Omega_\sigma$. Note
that $\cH_{\sigma ,0}$ is invariant under the action of $G$. We denote by
$\mathfrak{m}$ the unique invariant mean on the space $\WAP (G)$ of weakly
almost periodic bounded uniformly continuous functions on $G$.
A {\it flight function} is a function $f\in\WAP (G)$ such that
$\mathfrak{m} (|f|) = 0$, which is equivalent to the condition that for every
$\varepsilon > 0$ the set $\{ s\in G : | f(s) | < \varepsilon \}$ is
thickly syndetic (see \cite{BerRos,ETJ}).
We write $C_\sigma$ for the norm closure of
the linear span of the functions $s\mapsto \langle U_s \xi , \zeta \rangle$ on
$G$ for all $\xi , \zeta\in\cH_{\sigma ,0}$, which is a subspace of $\WAP (G)$.
We can alternatively describe $C_\sigma$ as the norm
closure of the linear span of the functions $s\mapsto\sigma (b\alpha_s (a))$
on $G$ for all $a,b\in A$ such that $\sigma (a) = \sigma (b) = 0$.

Recall that a strongly continuous unitary representation $\pi$ of $G$ on a Hilbert
space $\cH$ is said to be {\it weakly mixing} if for all $\xi , \zeta\in\cH$ the
function $s\to |\langle\pi (s)\xi , \zeta\rangle |$ on $G$ is a flight function.

\begin{definition}
Let $(A,G,\alpha )$ be a $C^*$-dynamical system. A $G$-invariant state
$\sigma$ on $A$ is said to be {\it weakly mixing} if the
representation of $G$ on $\cH_{\sigma ,0}$
is weakly mixing, i.e., if $f$ is a flight function for every $f\in C_\sigma$.
\end{definition}

\begin{definition}
Let $(A,G,\alpha )$ be a $C^*$-dynamical system. A $G$-invariant state
$\sigma$ on $A$ is said to be {\it T-mixing} if for all
$a_1 , a_2 , b_1 , b_2 \in A$ and $\varepsilon > 0$ the set of all $s\in G$ such that
$| \sigma (a_1 \alpha_s (b_1 ) a_2 \alpha_s (b_2 )) -
\sigma (a_1 a_2 ) \sigma (b_1 b_2 ) | < \varepsilon$
is thickly syndetic.
\end{definition}


\noindent Since a finite intersection of thickly
syndetic sets is thickly syndetic, the $G$-invariant state $\sigma$ is T-mixing
if and only if for every finite set $\Omega\subseteq A$ and $\varepsilon > 0$ the
set of all $s\in G$ such that $| \sigma (a_1 \alpha_s (b_1 ) a_2 \alpha_s (b_2 )) -
\sigma (a_1 a_2 ) \sigma (b_1 b_2 ) | < \varepsilon$ for all
$a_1 , a_2 , b_1 , b_2 \in\Omega$ is thickly syndetic. Similarly, $\sigma$ is
weakly mixing if and only if for every finite set $\Omega\subseteq A$ and
$\varepsilon > 0$ the set of all $s\in G$ such that
$| \sigma (a \alpha_s (b)) - \sigma (a) \sigma (b) | < \varepsilon$ for all
$a,b \in\Omega$ is thickly syndetic.

T-mixing is a weak form of a clustering property that has been studied 
in the context of quantum statistical mechanics (see for example \cite{BN}). Note that 
it implies weak mixing. If the system $(A,G,\alpha )$ is $\TSyn$-Abelian in the sense 
of Definition~\ref{D-Abelian} (in particular if it is commutative) then
T-mixing is equivalent to weak mixing (use Lemma~\ref{L-WM char}).
This is not the case
in general however, for if we take a $C^*$-probability space $(A,\sigma )$ and
consider the shift $^*$-automorphism $\alpha$ on the infinite reduced free product
$(B,\omega ) = (A,\sigma )^{*\Zb}$, then $\omega$ is weakly mixing but not
T-mixing. See also Example~\ref{E-Bog}.


\begin{lemma}\label{L-TSexpansive}
Let $(A,G,\alpha )$ be $C^*$-dynamical system with $A$ exact. Let
$\tau$ be a faithful T-mixing $G$-invariant tracial state on $A$.
Then $\alpha$ is $\TSyn$-expansive.
\end{lemma}

\begin{proof}
Let $V$ be a finite-dimensional operator subsystem of $A$ and let
$\varepsilon > 0$. Let $(\pi , \cH , \xi )$ be the GNS triple associated to
$\tau$. Since $\tau$ is faithful, $\pi$ is a faithful representation,
and so via $\pi$ we can view $A$ as acting on $\cH$ and $A\otimes A$ as
acting on the Hilbertian tensor product $\cH\otimes\cH$. By the injectivity of the
minimal operator space tensor product we can view $V\otimes V$ as an
operator subsystem of $A\otimes A$.

Let $\{ v_i \}_{i=1}^r$ be an Auerbach basis for $V$. Then
$\cS = \{ v_i \otimes v_j \}_{i,j=1}^r$ is an Auerbach basis for
$V\otimes V$.
By Lemma~\ref{L-perturb} there exists a $\delta > 0$ such that whenever
$W$ is an operator space, $\rho : V\otimes V \to W$ is a linear
isomorphism onto its image with $\max (\| \rho \|_\cb , \| \rho^{-1} \|_\cb )
\leq 1+\delta$, and $\{ w_{ij} \}_{i,j=1}^r$ is a subset of $W$ with
$\| \rho (v_i \otimes v_j ) - w_{ij} \| < 4\delta$ for all
$1\leq i,j \leq r$, the linear map $\psi : V\otimes V \to W$
determined on $\cS$ by $\psi (v_i \otimes v_j ) = w_{ij}$ is
an isomorphism onto its image with
$\max (\| \psi \|_\cb , \| \psi^{-1} \|_\cb ) \leq 1 + \varepsilon$.

Since $A$ is exact, $V$ is $1$-exact as an operator space, and so
by Lemma~17.8 of \cite{Pis} we can find a finite-dimensional
unital subspace $E\subseteq A$ such that, with $p$ denoting the
orthogonal projection of $\cH$ onto $E\xi$,
the (completely contractive) compression map
$\gamma : V\to pVp$ given by $v\mapsto pvp$ has an inverse
with c.b.\ norm less than $\sqrt[4]{1 + \delta}$.
Then the (completely contractive) compression map
$\theta = \gamma\otimes\gamma : V\otimes V \to (p\otimes p)(V\otimes V)
(p\otimes p)$ given by $z\mapsto (p\otimes p)z(p\otimes p)$ is invertible
and 
$\| \theta^{-1} \|_\cb =
\| \gamma^{-1} \otimes\gamma^{-1} \|_\cb\leq \| \gamma^{-1} \|^2_\cb <
\sqrt{1 + \delta}$.

Let $\{ a_i \xi \}_{i=1}^m$ be an orthonormal basis for $E\xi$. Then
$\cT = \{ a_i \xi \otimes a_j \xi \}_{i,j=1}^m$ is an orthonormal basis for
$E\xi\otimes E\xi$. Set $\delta' = (\dim (E))^{-4} \delta$, and
let $\eta$ be a positive real number to be further specified below. Define
$K$ to be the set of all $s\in G\setminus \{ e \}$ such that
\begin{enumerate}
\item $| \tau (a_k^* a_i \alpha_s (a_j a_l^* )) -
\tau (a_k^* a_i ) \tau (a_j a_l^* ) | < \eta$ for all
$1\leq i,j,k,l \leq m$, and

\item $| \tau (c^* v\alpha_s (w) a\alpha_s (bd^* )) -
\tau (c^* va) \tau (wbd^* ) | < \delta'$ for all
$v\otimes w\in\cS$ and $a\xi\otimes b\xi , c\xi\otimes d\xi \in \cT$.
\end{enumerate}
Since $\tau$ is T-mixing, the set $K$ is thickly syndetic.

Let $s\in K$.
Define $S : E\xi \otimes E\xi \to [E\alpha_s (E)]\xi$ to be the surjective
linear map determined on the basis $\cT$ by $S(a_i \xi \otimes a_j \xi ) =
a_i \alpha_s (a_j ) \xi$.
For $1\leq i,j,k,l\leq m$ we have
\begin{align*}
\langle S(a_i \xi \otimes a_j \xi ) , S(a_k \xi\otimes a_l \xi ) \rangle
&= \langle a_i \alpha_s (a_j ) \xi , a_k \alpha_s (a_l ) \xi \rangle \\
&= \tau ((a_k \alpha_s (a_l ))^* a_i \alpha_s (a_j )) 
= \tau (a_k^* a_i \alpha_s (a_j a_l^* )) \\
&\approx_\eta \tau (a_k^* a_i ) \tau (a_j a_l^* ) 
= \langle a_i \xi \otimes a_j \xi , a_k \xi \otimes a_l \xi \rangle .
\end{align*}
We thus see by a simple perturbation argument that by taking $\eta$ small
enough we can ensure that $S$ is invertible with
$\| S \| \| S^{-1} \| \leq \min (\sqrt{1 + \delta} , 2)$ and
$\| S^{-1} - S^* \| < \delta'$.

Denote by $q$ the orthogonal projection of $\cH$ onto $[E\alpha_s (E)]\xi$, and
let $\rho : V\otimes V \to \cB (q\cH )$ be the linear map
given by $\rho (z) = S\theta (z)S^{-1}$ for all $z\in V\otimes V$. Then $\rho$
is an isomorphism onto its image with inverse given by
$\rho^{-1} (z) = \theta^{-1} (S^{-1} zS)$,
and $\| \rho \|_\cb \leq \| S \| \| S^{-1} \| \leq \sqrt{1 + \delta}$ while
$\| \rho^{-1} \|_\cb \leq \| \theta^{-1} \|_\cb \| S^{-1} \| \| S \|
\leq 1+\delta$.

Now suppose we are given $v\otimes w\in\cS$.
For all $a\xi \otimes b\xi , c\xi \otimes d\xi \in\cT$, we have,
using the bound $\| S^{-1} - S^* \| < \delta'$ and
the definition of the set $K$,
\begin{align*}
\lefteqn{\langle S^{-1} qv\alpha_s (w)q S (a\xi \otimes b\xi ), c\xi \otimes d\xi \rangle
\approx_{\delta'} \langle v\alpha_s (w) a\alpha_s (b) \xi , c\alpha_s (d) \xi
\rangle} \hspace*{35mm} \\
&= \tau (\alpha_s (d^* ) c^* v\alpha_s (w) a\alpha_s (b)) 
= \tau (c^* v\alpha_s (w) a\alpha_s (bd^* )) \\
&\approx_{\delta'} \tau (c^* va) \tau (wbd^* ) 
= \tau (c^* va) \tau (d^* wb) \\
&= \langle va\xi , c\xi \rangle \langle wb\xi , d\xi \rangle 
= \langle va\xi \otimes wb\xi , c\xi \otimes d\xi \rangle \\
&= \langle \theta (v\otimes w) (a\xi \otimes b\xi ), c\xi \otimes d\xi \rangle .
\end{align*}
Thus $\| \theta (v\otimes w) - S^{-1} qv\alpha_s (w)q S \| < 2\delta'
(\dim (E))^4 = 2\delta$, and so
\[ \| \rho (v\otimes w ) - qv\alpha_s (w)q \| \leq
\| S \| \| \theta (v\otimes w) - S^{-1} qv\alpha_s (w)q S \| \| S^{-1} \| <
4\delta . \]
We conclude by our choice of $\delta$ that the linear map
$\psi : V\otimes V \to q[V\alpha_s (V)]q$ determined on elementary tensors
by $\psi (v \otimes w) = qv \alpha_s (w)q$
is an isomorphism onto its image with
$\max (\| \psi \|_\cb , \| \psi^{-1} \|_\cb ) \leq 1 + \varepsilon$.

Letting $\kappa : [V\alpha_s (V)] \to q[V\alpha_s (V)]q$ be the
cut-down map $z\mapsto qzq$, we have a commuting diagram
\begin{gather*}
\xymatrix{
V\otimes V \ar[r]^(0.46){\varphi} \ar[dr]_(0.43){\psi}
&[V\alpha_s (V)] \ar[d]^{\kappa} \\
& q[V\alpha_s (V)]q.}
\end{gather*}
where $\varphi$ is the 
linear map determined on elementary tensors by $v\otimes w \mapsto v\alpha_s (w)$.
Since $\kappa$ is completely contractive, the map $\varphi$ is invertible and
$\| \varphi^{-1} \|_\cb = \| \psi^{-1} \circ\kappa \|_\cb \leq
\| \psi^{-1} \|_\cb \leq 1 + \varepsilon$. Thus $K$ is a
$(1+\varepsilon )$-expansion set for $V$, and we obtain the result.
\end{proof}

\begin{theorem}\label{T-WMcontrInd}
Let $(A,G,\alpha )$ be a $\TSyn$-contractive $C^*$-dynamical system with $A$ exact.
Let $\tau$ be a faithful weakly mixing $G$-invariant tracial
state on $A$. Then $\alpha$ is $\TSyn$-independent.
\end{theorem}

\begin{proof}
By Proposition~\ref{P-contrAb}, $\alpha$ is $\TSyn$-Abelian. Hence
$\tau$ is T-mixing, in which case we can apply Lemma~\ref{L-TSexpansive} to conclude
that $\alpha$ is $\TSyn$-independent.
\end{proof}

Theorem~\ref{T-WMcontrInd} applies in particular to the commutative
situation, where the $\TSyn$-contrac-\linebreak tivity and exactness
hypotheses are automatic. In this case one can also obtain
$\TSyn$-independence in a combinatorial fashion
from the characterization of weak mixing
in terms of thickly syndetic sets and a partition of unity
argument (cf.\ the proof of Proposition~\ref{P-equiv-indep}(6)$\Rightarrow$(1)).
It follows, for example, that
if $X$ is a compact manifold, possibly with
boundary, of dimension at least $2$ and $\mu$ is a nonatomic Borel
probability measure on $X$ with full support which is zero on the boundary,
then a generic member of the set $\mathfrak{H}_\mu (X)$ of $\mu$-preserving homeomorphisms
from $X$ to itself equipped with the uniform topology is $\TSyn$-independent, since the
elements of $\mathfrak{H}_\mu (X)$ which
are weakly mixing for $\mu$ form a dense $G_\delta$ subset \cite{KaSt}
(see also \cite{AP}).

The assumption of $\TSyn$-contractivity in Theorem~\ref{T-WMcontrInd} cannot
be dropped in general. Indeed certain Bogoliubov automorphisms of the CAR algebra
are strongly mixing with respect to the unique tracial state but fail
to be $\Inf$-contractive, as the following example demonstrates.

\begin{example}\label{E-Bog}
Let $\cH$ be a separable infinite-dimensional Hilbert space over the
complex numbers.
We write $A(\cH )$ for the CAR algebra over $\cH$.
This is the unique, up to $^*$-isomorphism, unital $C^*$-algebra
generated by the image of an antilinear map $\xi\mapsto a(\xi )$
from $\cH$ to $A(\cH )$ for which the anticommutation relations
\begin{align*}
a(\xi ) a(\zeta )^* + a(\zeta )^* a(\xi ) &= \langle \xi ,\zeta \rangle
1_{A(\cH )} ,\\
a(\xi ) a(\zeta ) + a(\zeta ) a(\xi ) &= 0,
\end{align*}
hold for all $\xi ,\zeta \in\cH$ (see \cite{BR2}).
Every $U$ in the unitary group $\cU (\cH )$ gives
rise to a $^*$-automorphism $\alpha_U$ of $A(\cH )$, called a Bogoliubov
automorphism, by setting $\alpha_U (a(\xi )) = a(U\xi )$ for every $\xi\in\cH$.
The unique tracial state $\tau$ on $A(\cH )$ is given on products of the
form $a(\zeta_n )^* \cdots a(\zeta_1 )^* a(\xi_1 ) \cdots a(\xi_m )$ by
\[ \tau (a(\zeta_n )^* \cdots a(\zeta_1 )^* a(\xi_1 ) \cdots a(\xi_m ))
= \delta_{nm} \det [ \langle {\textstyle\frac12} \xi_i , \zeta_j \rangle ] . \]

Let $U$ be a unitary operator on $\cH$
such that $\lim_{|n| \to\infty} \langle U^n \xi,\zeta\rangle = 0$ for all
$\xi , \zeta\in\cH$ (for example, the bilateral shift with respect to some
orthonormal basis).
It is well known that the corresponding Bogoliubov automorphism $\alpha_U$
of $A(\cH )$ is strongly mixing for $\tau$, i.e.,
$\lim_{n\to\infty} | \tau (a\alpha_U^n (b)) - \tau (a)\tau (b) | = 0$
for all $a,b\in A(\cH )$ (see Example~5.2.21 in \cite{BR2}).
On the other hand, for every $\xi\in\cH$ we have
\begin{align*}
\| [ a(U^n \xi ) , a(\xi )] \| = 2 \| a(U^n \xi )a(\xi ) \| &=
2 \| a(\xi)^* a(U^n \xi )^* a(U^n \xi ) a(\xi )) \|^{1/2} \\
&\geq 2| \tau (a(\xi)^* a(U^n \xi )^* a(U^n \xi ) a(\xi )) |^{1/2} \\
&= \sqrt{\| \xi \|^4 - | \langle U^n \xi ,\xi \rangle |^2} ,
\end{align*}
and this last quantity converges to $\| \xi \|^2$ as $|n| \to \infty$.
This shows that $\alpha_U$ fails to be $\Inf$-Abelian and hence by
Proposition~\ref{P-contrAb} fails to be $\Inf$-contractive.

We also remark that $\tau$ fails to be T-mixing with respect to $\alpha_U$.
Indeed for $\xi\in\cH$ the quantity
$\tau (a(\xi )^* a(U^n \xi )^* a(\xi )a(U^n \xi ))$ is equal to
$\frac14 (|\langle U^n \xi ,\xi \rangle |^2 - \| \xi \|^4 )$, which
converges to $-\frac14 \| \xi \|^4 = -\tau (a(\xi )^* a(\xi ))^2$
as $n\to\infty$.
\end{example}


\section{Independence and weak mixing in UHF algebras}\label{S-UHF}


In the previous section we showed that, under certain conditions, the existence of a
weakly mixing faithful invariant state implies $\TSyn$-independence. What can be said
in the reverse direction? The transfer in dynamics from topology to measure theory
is a subtle one in general; for example, in \cite{LVRA} Huang and Ye exhibited
a u.p.e.\ $\Zb$-system which lacks an ergodic invariant measure of full
support. On the other hand, the presence of noncommutativity at the $C^*$-algebra
level can give rise to a kind of rigidity that in the commutative setting is
more characteristic of measure-theoretic structure. We will illustrate this phenomenon
in one of its extreme forms by showing that, for actions on a UHF algebra,
I-independence implies weak mixing for the unique tracial
state (Theorem~\ref{T-UHF}) and, for single $^*$-automorphisms in the type $d^\infty$ case,
is point-norm generic (Theorem~\ref{T-UHFdenseGdelta}).
We also obtain characterizations of I-independence for Bogoliubov actions
on the even CAR algebra in terms of measure-theoretic weak mixing
(Theorem~\ref{T-even-WM}).

We begin by recording a lemma which reexpresses Corollary~1.6 of \cite{BerRos} in
our $C^*$-dynamical context (see the discussion after Theorem~\ref{T-nuclear}
for information on weak mixing). Note that although $G$ is generally taken to be 
$\sigma$-compact and locally compact in \cite{BerRos}, the 
arguments on which Corollary~1.6 of \cite{BerRos} are based do not require these 
assumptions.

\begin{lemma}\label{L-WM char}
Let $(A,G,\alpha )$ be a $C^*$-dynamical system and let $\sigma$ be a
$G$-invariant state on $A$. Then $\sigma$ is weakly mixing if and only if
for every finite set $\Omega\subseteq A$ and $\varepsilon > 0$ there is
an $s\in G$ such that
$| \sigma (a^* \alpha_s (a)) - \sigma (a^* )\sigma (a) | < \varepsilon$ for all
$a\in\Omega$.
\end{lemma}

When $G$ is Abelian, it suffices in the above characterization of weak mixing
to quantify over singletons in $A$ instead of finite subsets.

\begin{lemma}\label{L-fdnear}
For every $d\in\Nb$ and $\varepsilon > 0$ there is a $\delta > 0$ such that
if $A_1$ and $A_2$ are finite-dimensional unital $C^*$-subalgebras of a
unital $C^*$-algebra $A$ with $\dim (A_2 ) \leq d$ and
$\| [a_1 ,a_2 ] \| \leq \delta \| a_1 \| \| a_2 \|$
for all $a_1 \in A_1$ and $a_2 \in A_2$, then there is a faithful
$^*$-homomorphism $\gamma$ from $A_2$ to the commutant $A_1'$ such that
$\| \gamma - \id_{A_2} \| < \varepsilon$.
\end{lemma}

\begin{proof}
Let $d\in\Nb$ and $\varepsilon > 0$, and
suppose that $A_1$ and $A_2$ are finite-dimensional unital
$C^*$-subalgebras of a unital $C^*$-algebra $A$. By Lemma~2.1 of
\cite{Bra} we can find a $\delta > 0$ depending on $\varepsilon$ and $d$
such that if the unit ball of $A_2$ is approximately included to within
$\delta$ in $A_1'$ then there exists a faithful
$^*$-homomorphism $\gamma : A_2 \to A_1'$ such that
$\| \gamma - \id_{A_2} \| < \varepsilon$. Let $\mu$ be the
normalized Haar measure on the unitary group $\cU (A_1 )$, which is compact by
finite-dimensionality. Then setting
\[ E(a) = \int_{\cU (A_1 )} uau^*\, d\mu (u) \]
for $a\in A$ we obtain a conditional expectation $E$ from $A$ onto
$A_1'$. Moreover, if $\| [a_1 ,a_2 ] \| \leq \delta
\| a_1 \| \| a_2 \|$ for all $a_1 \in A_1$ and $a_2 \in A_2$,
then $\| E(a) - a \| \leq \delta \| a \|$ for all $a\in A_2$, so that
the unit ball of $A_2$ is approximately included to within
$\delta$ in $A_1'$, yielding the lemma.
\end{proof}

Actually, by a result of Christensen \cite{near}, in the above
lemma a $\delta$ not depending on $d$ can be found.

\begin{theorem}\label{T-UHF}
Let $(A,G,\alpha )$ be an I-independent $C^*$-dynamical system with
$A$ a UHF algebra. Then $\alpha$ is weakly mixing
for the unique (and hence $\alpha$-invariant) tracial state $\tau$ on $A$.
\end{theorem}

\begin{proof}
By Propositions~\ref{P-C-equiv-indep} and \ref{P-contrAb}, $\alpha$ is
$\Inf$-Abelian. Let $B$ be a simple finite-dimensional unital
$C^*$-subalgebra of $A$, and let $\varepsilon > 0$.
Suppose we are given an $s\in G$ and a $\delta > 0$ such that
$\| [b_1 ,\alpha_s (b_2 )] \| \leq \delta
\| b_1 \| \| b_2 \|$ for all $b_1 , b_2 \in B$. By Lemma~\ref{L-fdnear},
if we assume $\delta$ to be sufficiently small we can find a
$^*$-homomorphism $\gamma : B\to A$ such that $\| \gamma - \alpha_s |_B \| <
\varepsilon$ and the $C^*$-subalgebras
$B$ and $\gamma (B)$ of $A$ commute. Denote by $\Phi$ the
$^*$-homomorphism from $B \otimes B$ to the $C^*$-algebra generated by
$B$ and $\gamma (B)$ determined on the factors by
$b\otimes 1\mapsto b$ and $1\otimes b \mapsto \gamma (b)$.
Since $B\otimes B$ is simple, the tracial state $\tau\circ\Phi$ is
unique and equal to $\tau |_B \otimes \tau |_B$. For all $a$ and $b$
in the unit ball of $B$, we then have
\[ \tau (a \alpha_s (b))
\approx_\varepsilon \tau (a \gamma (b)) \\
= (\tau\circ\Phi )(a \otimes b) \\
= \tau (a)\tau (b) , \]
so that the set of all $s\in G$ with
$| \tau (a \alpha_s (b)) - \tau (a)\tau (b) | < \varepsilon$ is nonempty.
Thus, since $A$ is the closure of an increasing sequence of simple
finite-dimensional unital $C^*$-subalgebras of $A$, we conclude by
Lemma~\ref{L-WM char} that $\tau$ is weakly mixing.
\end{proof}

\noindent The converse of Theorem~\ref{T-UHF} is false, as Example~\ref{E-Bog}
demonstrates.

We will next examine Bogoliubov actions on the even
CAR algebra. We refer the reader to Example~\ref{E-Bog} for notation
and an outline of the context (see \cite{BR2} for a general reference).
Let $\cH$ be a separable Hilbert space.
The even CAR algebra $A(\cH )_\even$ is the unital
$C^*$-subalgebra of the CAR algebra $A(\cH )$ consisting of those elements
which are fixed by the Bogoliubov automorphism arising from scalar multiplication
by $-1$ on $\cH$, and it is generated by even products of creation and
annihilation operators.
Both the CAR algebra and the even CAR algebra are $^*$-isomorphic to the
type $2^\infty$ UHF algebra (see \cite[Thm.\ 5.2.5]{BR2} and \cite{even}).
Every Bogoliubov automorphism $\alpha_U$ of $A(\cH )$
restricts to a $^*$-automorphism of $A(\cH )_\even$.
A strongly continuous unitary representation
$s\mapsto U_s$ of the group $G$ on
$\cH$ gives rise via Bogoliubov automorphisms to a $C^*$-dynamical system on
each of $A(\cH )$ and $A(\cH )_\even$, in which case we will speak of a
Bogoliubov action and write $\alpha_U$.

\begin{theorem}\label{T-even-WM}
Let $\cH$ be a separable Hilbert space. Let $U$ be a strongly continuous
unitary representation of $G$ on $\cH$,
and consider the corresponding Bogoliubov action $\alpha_U$ on
$A(\cH )_\even$. The following are equivalent:
\begin{enumerate}
\item $U$ is a weakly mixing representation of $G$,

\item the unique tracial state $\tau$ on $A(\cH )_\even$ is weakly mixing
for $\alpha_U$,

\item $\alpha_U$ is $\TSyn$-independent,

\item $\alpha_U$ is I-independent,

\item $\alpha_U$ is $\TSyn$-Abelian,

\item $\alpha_U$ is $\NE$-Abelian.
\end{enumerate}
\end{theorem}

\begin{proof}
(1)$\Rightarrow$(5). Since the representation $U$ of $G$ is weakly mixing,
for every finite set $\Theta\subseteq\cH$ and every $\varepsilon > 0$ the set of
all $s\in G$ for which $| \langle U_s \xi , \zeta \rangle | < \varepsilon$
for all $\xi ,\zeta\in\Theta$ is thickly
syndetic. It is then straightforward to check using the anticommutation relations
and the fact that the even CAR
algebra is generated by even products of creation and annihilation operators
that for every finite set $\Omega\subseteq A(\cH )_\even$ and $\varepsilon > 0$
the set of all $s\in G$ for which $\| [ a,\alpha_{U_s} (b) ] \| < \varepsilon$
for all $a,b\in\Omega$ is
thickly syndetic, i.e., $\alpha_U$ is $\TSyn$-Abelian.

(5)$\Rightarrow$(3). Apply Theorem~\ref{T-nuclear}.

(3)$\Rightarrow$(4). Apply Proposition~\ref{P-C-equiv-indep}.

(4)$\Rightarrow$(2). Apply Theorem~\ref{T-UHF}.

(2)$\Rightarrow$(1). Since $\tau$ is weakly mixing for $\alpha_U$, for all
$\xi , \zeta\in\cH$ the function
\[ s\mapsto \tau (a(\xi )^* a(U_s\zeta )) - \tau (a(\xi )^* ) \tau (a(\zeta ))
= \frac12 \langle U_s \zeta , \xi \rangle \]
on $G$ is a flight function, i.e., $U$ is weakly mixing.

(4)$\Leftrightarrow$(6). This is a consequence of
Theorem~\ref{T-nuclear} and Proposition~\ref{P-C-equiv-indep}.
\end{proof}


We will now restrict our attention to UHF algebras of the form
$M_d^{\otimes\Zb}$ for some $d\geq 2$ and establish generic I-independence
for $^*$-automorphisms.
For a $C^*$-algebra $A$ we denote by $\Aut (A)$ its $^*$-automorphism group
with the point-norm topology, which is a Polish space when $A$ is separable.

\begin{lemma}\label{L-I-indGdelta}
Let $A$ be a separable unital $C^*$-algebra. Then the
I-independent $^*$-automor-\linebreak phisms of $A$ form a $G_\delta$
subset of $\Aut (A)$.
\end{lemma}

\begin{proof}
For every finite-dimensional operator subsystem $V\subseteq A$ and
$\varepsilon > 0$ define
\begin{gather*}
\Gamma (V , \varepsilon ) =
\big\{ \alpha\in\Aut (A)
: \Ind (\alpha,V,\varepsilon' ) \neq\emptyset \text{ for some } \varepsilon'
\in (0,\varepsilon ) \big\} ,
\end{gather*}
which is an open subset of $\Aut (A)$, as can be seen by a straightforward
perturbation argument using Lemma~2.13.2 of \cite{Pis}.
Take an increasing sequence $V_1 \subseteq V_2 \subseteq
\dots$ of finite-dimensional operator subsystems of $A$ whose union is dense
in $A$. Then by Proposition~\ref{P-C-equiv-indep} the
set of I-independent $^*$-automorphisms of $A$ is equal to
$\bigcap_{j=1}^\infty \bigcap_{k=1}^\infty \Gamma (V_j ,1/k)$, which is a
$G_\delta$ set.
\end{proof}

For the definition of the Rokhlin property in the sense that we use it in
the next lemma, see \cite{Kis}.

\begin{lemma}\label{L-R-denseclass}
Let $\gamma$ be a $^*$-automorphism of $M_d^{\otimes\Zb}$ with the
Rokhlin property. Then $\gamma$ has dense conjugacy class in
$\Aut (M_d^{\otimes\Zb} )$.
\end{lemma}

\begin{proof}
Let $\alpha\in\Aut (M_d^{\otimes\Zb} )$, and let $\Omega$ be a finite subset
of $M_d^{\otimes\Zb}$ and $\varepsilon > 0$. To show that $\alpha$ can be
norm approximated to within $\varepsilon$ on $\Omega$ by a conjugate of
$\gamma$, we may assume that $\Omega\subseteq
M_d^{\otimes I}$ for some finite set
$I\subseteq\Zb$. By the classification theory for AF algebras, $\alpha$
is approximately inner. Thus, by enlarging $I$ if necessary, we can find
a unitary $u\in M_d^{\otimes I}$ such that the $^*$-automorphism
$\Ad u \otimes \id$ of $M_d^{\otimes I} \otimes M_d^{\otimes \Zb \setminus I}
= M_d^{\otimes\Zb}$ satisfies
$\| \alpha (a) - (\Ad u \otimes \id )(a) \| < \varepsilon$ for all
$a\in\Omega$. Pick a $^*$-automorphism $\gamma'$ of
$M_d^{\otimes \Zb \setminus I}$ conjugate to $\gamma$. Since
$\Omega\subseteq M_d^{\otimes I}$, we have
$\| \alpha (a) - (\Ad u \otimes \gamma' )(a) \| < \varepsilon$ for all
$a\in\Omega$. Since $\gamma$ has the Rokhlin property, so does
$\Ad u \otimes \gamma'$. Hence by Theorem~1.4 of \cite{Kis} there is a
$^*$-automorphism $\beta$ of $M_d^{\otimes\Zb}$ such that
$\| \Ad u \otimes \gamma' - \beta\circ\gamma\circ\beta^{-1} \| <
\varepsilon$. It follows that $\| \alpha (a) -
(\beta\circ\gamma\circ\beta^{-1} )(a) \| < \varepsilon$ for all
$a\in\Omega$, which establishes the lemma.
\end{proof}


\begin{theorem}\label{T-UHFdenseGdelta}
The I-independent $^*$-automorphisms of $M_d^{\otimes\Zb}$ form a
dense $G_\delta$ subset of $\Aut (M_d^{\otimes\Zb} )$.
\end{theorem}

\begin{proof}
The two-sided tensor product shift $\alpha$ on $M_d^{\otimes\Zb}$ satisfies 
the Rokhlin property \cite{BSKR,RPS} (note that the Rokhlin property for the
one-sided shift implies the Rokhlin property for the two-sided shift)
and thus has dense conjugacy class in
$\Aut (M_d^{\otimes\Zb} )$ by Lemma~\ref{L-R-denseclass}.
It is also clearly I-independent, and so with an appeal to
Lemma~\ref{L-I-indGdelta} we obtain the result.
\end{proof}


For systems on UHF algebras, I-independence is equivalent to
$\Inf$-Abelianness by Proposition~\ref{P-C-equiv-indep} and Theorem~\ref{T-nuclear}.
The following corollary thus ensues by applying Corollary~4.3.11 of \cite{BR1}.

\begin{corollary}
The invariant state space of a generic $^*$-automorphism in
$\Aut (M_d^{\otimes\Zb} )$ is a simplex.
\end{corollary}

Proposition~\ref{P-I-indep-cnr} yields the next corollary.

\begin{corollary}\label{C-generic-cnr}
A generic $^*$-automorphism of $M_d^{\otimes\Zb}$ is
completely untame.
\end{corollary}

\begin{corollary}
The set of inner $^*$-automorphisms of $M_d^{\otimes\Zb}$ is of first
category in $\Aut (M_d^{\otimes\Zb} )$.
\end{corollary}

\begin{proof}
By Corollary~\ref{C-generic-cnr} it suffices to show that no inner
$^*$-isomorphism of $M_d^{\otimes\Zb}$ is completely untame.
Let $u$ be a unitary in
$M_d^{\otimes\Zb}$. If $u$ is scalar then $\Ad u$ is the identity map,
which is obviously tame. If $u$ is not scalar
then $\Ad u$ fails to be completely untame because $u$, being fixed
by $\Ad u$, does not have an infinite $\ell_1$-isomorphism set.
\end{proof}

It is interesting to compare the generic I-independence in Theorem~\ref{T-UHFdenseGdelta} 
with generic behaviour for homeomorphisms
of the Cantor set. Kechris and Rosendal showed by model-theoretic means that the
Polish group of homeomorphisms of the Cantor set has a dense $G_\delta$
conjugacy class \cite{KR}, and a description of the generic homeomorphism has
been given by Akin, Glasner, and Weiss in \cite{AGW}. This generic homeomorphism 
can be seen to be null by examining as follows the construction of Section~1 in 
\cite{AGW}, to which we refer the reader for notation. 
The ``special'' homeomorphism $T(D,C)$ of the Cantor set $X(D,C)$ is
defined at the end of Section~1 in \cite{AGW} and represents 
the generic conjugacy class. It is a product of the homeomorphism
$\tau_{(D,C)}$ of the Cantor set $Z(D,C)$ and an identity map, and so it suffices
to show that $\tau_{(D,C)}$ is null. Now $Z(D,C)$ is a closed subset of
$q(Z(D, C)) \times \Theta$, and $\tau_{(D,C)}$ is the restriction of the product
of an obviously null homeomorphism of $q(Z(D,C))$
and the universal adding machine on $\Theta$, which is an inverse limit of finite systems
and hence is also null. It follows that $\tau_{(D,C)}$ is null, as desired.

As a final remark, we point out that the $^*$-automorphisms of $M_d^{\otimes\Zb}$
which are weakly mixing for $\tau$ form a dense $G_\delta$ subset of
$\Aut (M_d^{\otimes\Zb} )$. This is easily seen using Lemma~\ref{L-WM char} and the
fact that tensor product shift on $M_d^{\otimes\Zb}$ is weakly mixing for $\tau$
and has dense conjugacy class in $\Aut (M_d^{\otimes\Zb} )$
(cf.\ the proof of Theorem~\ref{T-UHFdenseGdelta}).


\section{A tame nonnull Toeplitz subshift}\label{S-Toeplitz}

We construct in this section a tame nonnull Toeplitz subshift.
Toeplitz subshifts were introduced by Jacobs and Keane in \cite{JK}.
An element $x\in \Omega_m:=\{0, 1, \dots, m-1\}^{\Zb}$ is called a
{\it Toeplitz sequence} if for any $j\in \Zb$ there exists an
$n\in \Nb$ such that $x(j+kn)=x(j)$ for all $k\in \Zb$. The subshift
generated by $x$ is called a {\it Toeplitz subshift}. Note that
every Toeplitz subshift is minimal.

Set $m=2$. To construct our Toeplitz sequence $x$ in $\Omega_2$, we
choose an increasing sequence $n_1<n_2<\dots$ in $\Nb$ with
$n_j|n_{j+1}$ for each $j\in \Nb$ and $j\cdot2^j+1$ distinct
elements $y_{j, 0}, y_{j, 1}, \dots, y_{j, j\cdot2^j}$ in
$\Zb/n_j\Zb$ for each $j\in \Nb$ with the following properties:
\begin{enumerate}
\item[(i)] $y_{j+1, k}\equiv y_{j, 0} \!\mod{n_j}$ for all $j\in \Nb$ and all
$0\le k\le (j+1)2^{j+1}$,

\item[(ii)] for each $j\in \Nb$ and $1\le t\le 2^j$, setting $Y_{j,
t}:=\{y_{j, k}: (t-1)j<k\le tj\}$, we have, for all $1\le t<2^j$, $Y_{j,
t+1}=Y_{j, t}+z_{j, t}$ for some $z_{j, t}\in \Zb/n_j\Zb$,

\item[(iii)] there does not exist any $z\in \Zb$ with
$z\equiv y_{j, 0} \!\mod{n_j}$ for all $j\in \Nb$,

\item[(iv)] for each $j\in \Nb$ and any $1\le k_1, k_2\le j\cdot 2^j$ and any
$0\le k_3\le j\cdot 2^j$
we have $y_{j, k_1}-y_{j, 0}\neq y_{j, k_3}-y_{j, k_2}$.
\end{enumerate}
Take a map $f:\bigcup_{j\in \Nb, 1\le t\le 2^j}Y_{j, t}\rightarrow \{0,
1\}$ such that, for each $j\in \Nb$, the maps $\{1, 2, \dots,
j\}\rightarrow \{0, 1\}$ given by $p\mapsto f(y_{j, (t-1)j+p})$
yield exactly all of the elements in $\{0, 1\}^{\{1, 2, \dots, j\}}$
as $t$ runs through $1, \dots, 2^j$.
Now we define our $x$ by
\begin{gather*}
x(s)=
\begin{cases}
f(y_{j, k}), \mbox{ if } s\equiv y_{j, k} \!\!\mod{n_j} \mbox{ for some }
j\in \Nb \mbox{ and some } 1\le k\le 2^j ,\\
0, \mbox{ otherwise}.
\end{cases}
\end{gather*}
Property (i) guarantees that $x$ is well defined. Property
(iii) implies that $x$ is a Toeplitz sequence. Denote by $X$ the
subshift generated by $x$. Also denote by
$A$ (resp.\ $B$) 
the set of elements in $X$ taking value $1$ (resp.\ $0$)
at $0$.
Property (ii) and the condition on $f$
imply that $\tilde{y}_{j, 1}, \dots, \tilde{y}_{j, j}$ is an
independence set for $(A, B)$ for each $j\in \Nb$,
where $\tilde{y}_{j, i}$ is any element in $\Zb$ with $\tilde{y}_{j,
i}\equiv y_{j, i} \!\mod n_j$. Thus $(A, B)$ has
arbitrarily large finite independence sets and hence the subshift
$X$ is nonnull by Proposition~\ref{P-basic N}.

Since $\IT_2(X, \Zb)$ is $\Zb$-invariant,
to show that $X$ is tame, by Proposition~\ref{P-basic R}(2)
it suffices to show that $(A, B)$ has no infinite
independence set.
For each $s\in \Zb$ with $x(s)=1$ let $J(s)$ and $K(s)$ denote the
positive integers such that $1\leq K(s) \leq J(s)2^{J(s)}$
and $s\equiv y_{J(s), K(s)} \!\mod n_{J(s)}$.

\begin{lemma} \label{L-Toeplitz 1}
Suppose that $x(s_1)=x(s_2)=1$ for some $s_1$ and $s_2$ in $\Zb$
with $J(s_1)<J(s_2)$.
Also suppose that $x(s_1+a)=0$ and $x(s_2+a)=1$ for some $a\in \Zb$.
Then $J(s_2+a)=J(s_1)$.
\end{lemma}

\begin{proof}
If $J(s_2+a)>J(s_1)$, then
\[s_1-s_2\equiv (s_1+a)-(s_2+a)\equiv
(s_1+a)-s_2 \!\!\mod n_{J(s_1)}\]
by property (i) and hence $s_1\equiv s_1+a \!\mod n_{J(s_1)}$.
Consequently $x(s_1+a)=x(s_1)=1$, in contradiction to the fact that
$x(s_1+a)=0$. Therefore $J(s_2+a)\le J(s_1)$.

If $J(s_2+a)<J(s_1)$, then $s_1\equiv s_2 \!\mod n_{J(s_2+a)}$ by
property (i) and hence $s_1+a\equiv s_2+a \!\mod n_{J(s_2+a)}$.
Consequently, $x(s_1+a)=x(s_2+a)=1$, in contradiction to the fact that
$x(s_1+a)=0$. Therefore $J(s_2+a)\ge J(s_1)$ and hence $J(s_2+a)=J(s_1)$, as
desired.
\end{proof}

\begin{lemma}\label{L-Toeplitz 2}
Suppose that $\{s_1, s_2, s_3\}$ is an independence set for 
$(A,B)$ and $x(s_i)=1$ for all $1\le i\le 3$. Then
$J(s_1)=J(s_2)=J(s_3)$.
\end{lemma}

\begin{proof} Suppose that the $J(s_i)$ are not all the same.
Without loss of generality, we may assume that $J(s_1)<
\min (J(s_2),J(s_3))$ or $J(s_1)=J(s_3)<J(s_2)$.

Consider the case $J(s_1)<\min (J(s_2), J(s_3))$ first. Since $\{s_1, s_2,
s_3\}$ is an independence set for $(A, B)$, we have
$x(s_1+a)=x(s_3+a)=0$ and $x(s_2+a)=1$ for some $a\in \Zb$. By
Lemma~\ref{L-Toeplitz 1} we have $J(s_2+a)=J(s_1)$. Note that
$s_3\equiv s_2\mod n_{J(s_1)}$ by property (i). Consequently,
$s_3+a\equiv s_2+a \mod n_{J(s_1)}$ and hence $x(s_3+a)=x(s_2+a)=1$,
contradicting the fact that $x(s_3+a)=0$. This rules out the case
$J(s_1)<\min (J(s_2),J(s_3))$.

Consider now the case $J(s_1)=J(s_3)<J(s_2)$. Since $\{s_1, s_2,
s_3\}$ is an independence set for $(A, B)$, we have
$x(s_1+a)=0$ and $x(s_2+a)=x(s_3+a)=1$ for some $a\in \Zb$. By
Lemma~\ref{L-Toeplitz 1} we have $J(s_2+a)=J(s_1)$. Note that
$s_3\equiv s_2 \!\mod n_{j}$ for any $j<J(s_1)\le \min (J(s_2), J(s_3))$ by
property (i). If $J(s_3+a)<J(s_1)=J(s_2+a)$, then $s_3+a\not
\equiv s_2+a \!\mod n_{J(s_3+a)}$ leading to a contradiction. Thus
$J(s_3+a)\ge J(s_1)$. We then have
\begin{align*}
y_{J(s_3),K(s_3)}-y_{J(s_1), 0} &\equiv y_{J(s_3),
K(s_3)}-y_{J(s_2), K(s_2)}\\
&\equiv y_{J(s_3+a), K(s_3+a)}-y_{J(s_2+a), K(s_2+a)}\!\!\mod n_{J(s_1)}
\end{align*}
in contradiction to property (iv). Therefore the case
$J(s_1)=J(s_3)<J(s_2)$ is also ruled out, and we obtain the lemma.
\end{proof}

\begin{lemma}\label{L-Toeplitz 3}
Suppose that $H$ and $H'$ are disjoint nonempty subsets of $\Zb$ and
that $H\cup H'$ is an independence set for
$(A, B)$. Suppose also that $|H|\ge j\cdot (2^j+1)$. Then
there exist an $H_1\subseteq H$ and an $a\in \Zb$ such that
$x(s+a)=1$ and $J(s+a)>j$ for all $s\in H'\cup (H\setminus H_1)$ and $|H_1|\le j$.
\end{lemma}

\begin{proof} We shall prove the assertion via induction on $j$. The case
$j=0$ is trivial since $H'\cup H$ is an independence set for
$(A, B)$.

Assume that the assertion holds for $j=i$. Suppose that
$|H|\ge (i+1)\cdot (2^{i+1}+1)$.
By the assumption there exist an $H_1\subseteq H$ and an $a\in \Zb$ such that
$x(s+a)=1$ and $J(s+a)>i$ for all $s\in H'\cup (H\setminus H_1)$ and $|H_1|\le i$.
Note that $|H\setminus H_1|\ge (i+1)2^{i+1}+1\ge 3$.
By Lemma~\ref{L-Toeplitz 2}, $J(s+a)$ does not
depend on $s\in H'\cup (H\setminus H_1)$. If $J(s+a)>i+1$ for $s\in H'\cup (H\setminus H_1)$, we are done with the induction step.
Suppose then that $J(s+a)=i+1$ for all $s\in H'\cup (H\setminus H_1)$. Then $K(s_1+a)=K(s_2+a)$ for some distinct
$s_1$ and $s_2$ in $H\setminus H_1$. 
In other words, $s_1\equiv s_2 \mod n_{i+1}$. Set $H_2:=H_1\cup \{s_1\}$.
Since $H'\cup H$ is an independence set
for $(A, B)$, there exists a $b\in \Zb$ with $x(s+b)=1$
for all $s\in H'\cup (H\setminus H_2)$ and
$x(s_1+b)=0$. Using property (i) one sees easily that
$J(s_2+a)<J(s_2+b)$. By Lemma~\ref{L-Toeplitz 2} we
have $J(s+b)=J(s_2+b)>i+1$ for all $s\in H'\cup (H\setminus H_2)$.
This finishes the induction step and proves the lemma.
\end{proof}

Suppose that $H\subseteq \Zb$ is an infinite independence set for
$(A, B)$.
Choose a nonempty finite subset $H'\subset H$. By Lemma~\ref{L-Toeplitz 3},
for any $j\in \Nb$ there exists an $a_j\in \Zb$ such that
$x(s+a_j)=1$ and $J(s+a_j)>j$ for all $s\in H'$.
%
%
%
%
By property (i) we have $s+a_j\equiv s'+a_j \!\mod
n_j$ for all $s, s'\in H'$. Thus $n_j | s-s'$. Since
$n_j\to \infty$ as $j\to \infty$, we see that $H'$ cannot contain more than
one element, which is a contradiction. Therefore
$(A, B)$ has no infinite independence sets and
we conclude that $X$ is tame.

\section{Convergence groups}\label{S-convergence}

Boundary-type actions of free groups are often described as exhibiting
chaotic behaviour. The geometric features which
account for this description are tied to complexity within the group
structure and are suggestive of the free-probabilistic notion
of independence. We will see here on the other hand that any convergence group
(in particular, any hyperbolic group acting on its
Gromov boundary) displays a strong lack of dynamical independence in our sense.

Let $(X,G)$ be a dynamical system with $G$ a group and 
$X$ containing at least $3$ points.
We call a net $\{s_i\}_{i\in I}$ in $G$ {\it wandering} if
for any $s\in G$, $s_i \neq s$ for sufficiently large $i\in I$.
Recall that $G$ is said to act as a {\it (discrete) convergence group} 
on $X$ if
for any wandering net $\{s_i\}_{i\in I}$ in $G$ there exist
$x, y\in X$ and a subnet $\{s_j\}_{j\in J}$ of $\{s_i\}_{i\in I}$
such that $s_j K\to y$ as $j\to\infty$ for
every compact set $K\subseteq X\setminus \{x\}$ \cite{Bowditch}.
When $X$ is metrizable, one can replace ``wandering net'' and
``subnet'' by ``sequence of distinct elements'' and ``subsequence'', 
respectively, in the above definition. 

\begin{lemma}\label{L-bound}
Let $G$ act as a convergence group on $X$.
Let $Z_1$ and $Z_2$ be two disjoint closed subsets of $X$.
Then there exists a finite subset $F\subseteq G$ such that
$HH^{-1}\subseteq F$ for every independence set $H$ of $Z_1$ and $Z_2$.
\end{lemma}

\begin{proof}
Suppose that such an $F$ does not exist. Then we can find independence
sets $H_1, H_2, \dots $ for the pair $(Z_1 , Z_2 )$ such that
$H_nH^{-1}_n\not \subseteq \bigcup^{n-1}_{k=1}H_kH^{-1}_k$ for all $n\geq 2$.
Choose an $s_n\in  H_nH^{-1}_n\setminus \bigcup^{n-1}_{k=1}H_kH^{-1}_k$
for each $n$.
Since independence sets are right translation invariant, $\{e, s_n\}$
is an independence set for $(Z_1 , Z_2 )$ when $n\ge 2$.
As $G$ acts as a convergence group on $X$
we can find $x, y\in X$ and a subnet $\{s_j\}_{j\in J}$ of $\{s_n\}_{n\in \Nb}$
such that $s_j K\rightarrow y$ as $j\rightarrow \infty$ for
every compact set $K\subseteq X\setminus \{x\}$. Without loss of generality we
may assume that $x\notin Z_1$. Then $s_j Z_1 \rightarrow y$ as
$j\rightarrow \infty$.
We separate the cases $y\notin Z_1$ and $y\in Z_1$. If $y\notin Z_1$,
then $s_j Z_1 \cap Z_1 = \emptyset$ when $j$ is large enough, contradicting
the fact that $\{e, s_j\}$
is an independence set for $(Z_1 , Z_2 )$. If $y\in Z_1$, then $s_j Z_1 \cap Z_2
=\emptyset$ when $j$ is large enough, again contradicting
the fact that $\{e, s_j\}$
is an independence set for $(Z_1 , Z_2 )$.
Therefore there does exist an $F$ with the desired property.
\end{proof}

\begin{theorem}\label{T-convergence}
Let $G$ act as a convergence group on $X$. Then $(X, G)$ is null.
\end{theorem}

\begin{proof}
Given disjoint closed subsets $Z_1$ and $Z_2$ of $X$ let $F$ be as in Lemma~\ref{L-bound}.
Then $|H|\le |HH^{-1}|\le |F|$ for any independence set $H$ for
the pair $(Z_1 , Z_2 )$.
Thus $Z_1$ and $Z_2$ do not have arbitrarily large finite independence sets.
It follows by Proposition~\ref{P-basic N}(2) that $(X, G)$ is null.
\end{proof}

An action of $G$ on a locally compact 
Hausdorff space $Y$ is said to be {\it properly discontinuous}
if, given any compact subset $Z\subseteq Y$, $sZ\cap Z\neq \emptyset$
for only finitely many $s\in G$.  
For a general reference on hyperbolic spaces
and hyperbolic groups, see \cite{BH}. 
It is a theorem of Tukia that if a group $G$ acts properly discontinuously
as isometries
on a proper geodesic hyperbolic space $Y$, then the induced action
of $G$ on the Gromov compactification $\overline{Y}=Y\cup \partial Y$ is
a convergence group action \cite[Theorem 3A]{Tukia}. 
Thus we get:

\begin{corollary}\label{C-boundary}
If a group $G$ acts properly discontinuously as isometries
on a proper geodesic hyperbolic space $Y$, then the induced action
of $G$ on the Gromov compactification $\overline{Y}$ is null. In particular,
the action of a hyperbolic group on its Gromov boundary is null.
\end{corollary}

It happens on the other hand that, as the following remark notes, 
nonelementary convergence group actions
fail to be $\HNS$ (hereditarily nonsensitive \cite[Definition 9.1]{GM}), 
although $\HNS$ dynamical
systems are tame (see the proof of \cite[Corollary 5.7]{DEBS}).

\begin{remark}\label{R-convergence vs equicontinuous}
For any dynamical 
system $(X, G)$, denote by $L(X, G)$ its {\it limit set}, defined as
the set of all $x\in X$ such that $| \{ s\in G : U\cap sU \neq\emptyset\} |
=\infty$ for every neighbourhood $U$ of $x$.
Clearly $L(X, G)$ is a closed $G$-invariant subset of $X$.
A convergence group action of $G$ on $X$ is 
{\it nonelementary} if $|L(X,G)|>2$, or, equivalently, there is no
one- or two-point subset of $X$ fixed setwise by $G$ \cite[Theorem 2T]{Tukia}.
In this event, $L(X, G)$ is an infinite perfect set, the action of $G$ on
$L(X, G)$ is minimal, and there exists a {\it loxodromic} 
element $s\in G$ \cite[Theorem 2S, Lemmas 2Q, 2D]{Tukia},
i.e., $s$ has exactly two distinct fixed points $x$ and $y$, 
and $s^n K\to y$ as
$n\to \infty$ for any compact set $K\subseteq X\setminus\{x\}$.
Notice that $\{x, y\}\subseteq L(X, G)$. 
Clearly the action of $G$ on $L(X, G)$ is
nonequicontinuous in this event. Since minimal $\HNS$
dynamical systems are 
equicontinuous and subsystems of $\HNS$ dynamical systems are $\HNS$, 
we see that no nonelementary convergence group action is $\HNS$.
\end{remark}

\end{document}